\documentclass{lmcs}
\pdfoutput=1
\usepackage[utf8]{inputenc}

\usepackage{lastpage}
\lmcsdoi{21}{3}{32}
\lmcsheading{}{\pageref{LastPage}}{}{}%
{Jul.~10,~2023}{Sep.~26,~2025}{}


\usepackage{hyperref}
\theoremstyle{plain}

\usepackage{float}
\usepackage{mathpartir}
\usepackage{tikz-cd}
\usepackage{calc}

\begin{document}

\definecolor{green}{rgb}{0.54, 0.66, 0.36}
\newcommand{\judge}[2]{\boldsymbol{\lfloor}#1\boldsymbol{\rfloor}\;#2}
\newcommand{\judgectx}[1]{\boldsymbol{\lfloor}#1\boldsymbol{\rfloor}}
\newcommand{\judgext}[2]{\boldsymbol{\lfloor}#1\boldsymbol{\rfloor}_{\textnormal{\texttt{ext}}}\;#2}
\newcommand{\type}{\textsc{Type}}
\newcommand{\tp}{\textsc{tp}}
\newcommand{\op}{\text{op}}
\newcommand{\discr}{\textsc{Discr}}
\newcommand{\grpd}{\textsc{Grpd}}
\newcommand{\gr}{\int}
\newcommand{\cattypes}{(\mathcal{C},\tp,-.-,P)}
\newcommand{\cattypesc}{(\mathcal{C},\tp_{\mathcal{C}},-.-,P_{\mathcal{C}})}
\newcommand{\cattypesd}{(\mathcal{D},\tp_{\mathcal{D}},-.-,P_{\mathcal{D}})}
\newcommand{\equivpair}[2]{\langle\;#1\;|\;#2\;\rangle}
\newcommand{\splitt}{\mathsf{split}}
\newcommand{\pair}{\mathsf{p}}
\newcommand{\pairing}{\mathsf{pair}}
\newcommand{\ev}{\mathsf{ev}}
\newcommand{\f}{\textnormal{\texttt{f}}}
\newcommand{\g}{\textnormal{\texttt{g}}}
\newcommand{\p}{\textnormal{\texttt{p}}}
\newcommand{\q}{\textnormal{\texttt{q}}}
\newcommand{\ci}{\textnormal{\texttt{c}}}
\newcommand{\di}{\textnormal{\texttt{d}}}
\newcommand{\choice}{\textsc{q}}
\renewcommand\subset{\subseteq}
\renewcommand\phi{\varphi}
\newcommand{\id}{\mathsf{id}}
\newcommand{\pid}{\mathbf{{\Pi}}}
\newcommand{\sigmad}{\mathbf{{\Sigma}}}
\renewcommand\subset{\subseteq}
\renewcommand\phi{\varphi}
\newcommand{\C}{\boldsymbol{\mathfrak{C}}}
\newcommand{\ettm}{\mathbf{eTT}}
\newcommand{\pttm}{\mathbf{pTT}}
\newcommand{\spttm}{\textnormal{\texttt{h}}\mathbf{pTT}}
\newcommand{\ett}{\textnormal{ETT}}
\newcommand{\ptt}{\textnormal{PTT}}
\newcommand{\sptt}{\textnormal{hPTT}}
\newcommand{\mltt}{\textnormal{MLTT}}
\newcommand{\uip}{\textnormal{UIP}}
\newcommand{\1}{\textnormal{\texttt{1}}}
\newcommand{\T}{\textsc{T}}
\newcommand{\refl}{\mathsf{r}}
\renewcommand{\j}{\mathsf{J}}
\newcommand{\h}{\mathsf{H}}
\newcommand{\funext}{\mathsf{funext}}
\renewcommand{\sectionautorefname}{Section}
\renewcommand{\subsectionautorefname}{Subsection}

\title[Relating homotopy equivalences to conservativity]{Relating homotopy equivalences to conservativity\texorpdfstring{\\}{} in dependent type theories with computation axioms}
\titlecomment{Earlier preprint versions of this article were titled \textit{A conservativity result for homotopy elementary types in dependent type theory} and \textit{Relating homotopy equivalences to conservativity in dependent type theories with propositional computation}.}
\thanks{Research supported by a School of Mathematics full-time EPSRC Doctoral Training Partnership Studentship 2019/2020, by an E-COST-MEETING-CA20111 WG6 2022, and by an E-COST-CA20111 Short-Term Scientific Mission Grant.}

\author[M.~Spadetto]{Matteo Spadetto\lmcsorcid{0000-0002-6495-7405}}
\address{University of Nantes, France}
\email{matteo.spadetto.42@gmail.com}

\begin{abstract}
  \noindent We prove a conservativity result for extensional type theories over propositional ones, i.e. dependent type theories with propositional computation rules, or computation axioms, using insights from homotopy type theory. The argument exploits a notion of canonical homotopy equivalence between contexts, and uses the notion of a category with attributes to phrase the semantics of theories of dependent types. Informally, our main result asserts that, for judgements essentially concerning h-sets, reasoning with extensional or propositional type theories is equivalent.
\end{abstract}

\maketitle

\section{Introduction}\label{section2.1}

In recent decades, dependent type theory has emerged as a powerful tool in the foundations of mathematics. Dependent types, i.e. types varying over the terms of other types, allow for precise and expressive specifications of mathematical statements and proofs, and provide a formal language for reasoning about them \cite{nordstrom2000martin,COQUAND198895}. This is not just of theoretical significance, as dependent types also have practical applications in fields such as software verification, where proof assistants like Coq and Agda have been successfully employed. More recently, the field has seen significant developments, based on the numerous insights provided by the foundational work of
Vladimir Voevodsky in univalent foundations. This includes the emergence of Homotopy Type Theory \cite{hottbook} and the Univalent Foundations program \cite{vlad}, which provide a new understanding of the concept of equality. The focus of this paper is on the notion of conservativity between two different theories of dependent types, exploiting the insights coming from Homotopy Type Theory.

Suppose that we are given two dependent type theories $\T_1$ and $\T_2$ such that the inference rules, hence the type constructors, of $\T_1$ are contained (or inferable) in $\T_2$. In other words, let us assume that $\T_2$ extends $\T_1$. The first question that one can ask is whether $\T_1$ and $\T_2$ can deduce the same judgements, in which case the theories actually coincide (or are logically equivalent, respectively). However, even when this is not the case, one can still ask a question which is perhaps even more interesting: whether the two theories prove the same \textit{statements}. By following the general interpretation of \textit{statements as types} and \textit{proofs as terms} (see \cite{Howard1980-HOWTFN-2, martin1984intuitionistic, hottbook} for more details), this property corresponds to the following: the theory $\T_1$ considers any given type inhabited, provided it is inhabited according to $\T_2$, i.e. whenever $\T_2$ proves a term judgement ${t:T}$, then $\T_1$ proves a term judgement of the form ${\tilde{t}:T}$. When this happens, we say that $\T_2$ is \textit{conservative} over $\T_1$, meaning that \textit{by ``reinforcing'' the theory $\T_1$ with the rules of $\T_2$ one does not risk to modify---namely increase---the deductive power of $\T_1$}.

The property of conservativity was studied by Martin Hofmann \cite{phdhofmann, hofmannconservativity} for an intensional type theory with \textit{extensional concepts} as $\T_1$ and the extensional type theory for $\T_2$. In this case a full conservativity result actually follows, but the additional extensional requirements on $\T_1$---e.g. the identity proof irrelevance: $$\judge{x:A,\;p:x=x}{p= \refl(x)}$$ for paths in between a term and itself---are fundamental. In fact, a dependent type theory with fully intensional identity types does not make the type: $$\judge{x,y:A,\;p,q:x=y}{p=q:\type}$$ inhabited, while the extensional type theory does. This is the main topic of our work.

\medskip

\noindent \textbf{In this paper} we consider dependent type theory $\T$, that we may call \textit{propositional}, seen as a further weakening of a fully intensional dependent type theory, and we look for a concrete family of type judgements such that the extensional type theory is conservative over $\T$ \textit{relative} to that family, i.e. whenever the extensional type theory makes one of these types inhabited, so does $\T$---see \autoref{main result}. In detail, the theory $\T$ is going to be endowed with \textit{propositional}, or \textit{axiomatic} \cite{ottenspadetto2025,spadetto25}, \textit{identity types}.

A dependent type theory is said to have \textit{propositional identity types} if it is endowed with a type constructor consisting of the usual formation, introduction, and elimination rules of intensional identity type---see \autoref{intensional id}---except for the judgemental equality of its computation rule. The latter is only required to hold in a weakened form, called \textit{propositional form}, i.e. it is \textit{replaced by a propositional equality}: whenever we are given judgements: \begin{center} $\judge{x,y:A;\;p:x=y}{C(x,y,p):\type}$ \;\;\;\; and \;\;\;\; $\judge{x:A}{q(x):C(x,x,r(x))}$ \end{center} and hence a term judgement $\judge{x,y:A;\;p:x=y}{\j(q,x,y,p):C(x,y,p)}$ by the elimination rule, then, in place of asking that the term equality judgement $\judge{x:A}{\j(q,x,x,r(x))\equiv q(x)}$ holds, we only ask that it holds \textit{propositionally}, i.e. that an additional term judgement of the form: $$\judge{x:A}{\h(q,x):\j(q,x,x,r(x))=q(x)}$$ holds.

Cohen, Coquand, Danielsson, Huber, and M\"ortberg \cite{MR3124820,MR3781068} conducted initial analyses related to propositional identity types. Another work \cite{MR3281415} introduces a univalent model of $\mltt$ where the propositional computation rule for identity types is validated, although its judgmental version is not. Following this, the type constructor has been thoroughly examined by van den Berg and Moerdijk \cite{MR3828037, MR3795638}, who introduced and explored a notion of semantics for dependent type theories with propositional identity types, using the notion of a \textit{path category}. One may consider the same form of weakening for the computation rule of dependent sum types and dependent product types: these type constructors satisfying a propositional computation rule are called \textit{propositional} (or \textit{axiomatic}) \textit{dependent sum types} and \textit{propositional} (or \textit{axiomatic}) \textit{dependent product types} respectively \cite{ottenspadetto2025,spadetto25}. The dependent type theory $\T$ that we are going to consider is therefore a dependent type theory having propositional identity types, propositional dependent sum types, and propositional dependent product types with \textit{function extensionality} in propositional form, together with an arbitrary family of atomic types and atomic terms. We call such a theory \textit{propositional dependent type theory}. The aim of this paper is to provide a semantic proof that the corresponding extensional type theory is conservative over the propositional type theory $\T$ relative to the family of type judgements of the latter obtained by inductively applying the type formation rules to the atomic types that are provably h-sets in the latter.

The property of conservativity between various weakenings of extensional and intensional theories has been studied by several authors, in addition to the already cited Martin Hofmann \cite{hofmannconservativity}. Bocquet \cite{bocquet} examines the property of \textit{Morita equivalence} between such a propositional theory and its extension obtained by strictifying the computation rule for identity types. Kapulkin and Li \cite{zbMATH07978139} reformulate and prove Hofmann's result in terms of such a Morita equivalence. Winterhalter, Sozeau, and Tabareau \cite{winterhalter:hal-01849166}, building on Hofmann's point of view and on an approach devised by Oury \cite{MR2197014}, define a translation from extensional theories to intensional ones with uniqueness of identity proofs ($\uip$) and function extensionality. This translation is then adapted by Boulier and Winterhalter \cite{boulierwinterhalter,Winterhalter2020} to obtain one from extensional theories to propositional ones, the latter again with uniqueness of identity proofs and function extensionality. This adaptation is used to achieve a completely syntactical proof of the conservativity of the former over the latter.

The present paper can thus be seen as an alternative, semantic, approach to the same problem addressed by Boulier and Winterhalter, building on Hofmann's argument. Our approach involves considering a propositional theory of dependent types, where dependent product types include function extensionality in propositional form, without initially requiring the uniqueness of identity proofs. The argument we follow gradually leads us to restrict the theory to contexts, called \textit{h-elementary}, generated only from h-sets. Consequently, we end up with a result formulated more abstractly than the one contained in \cite{boulierwinterhalter,Winterhalter2020}, with a specific concrete case that essentially matches the latter. This approach also enables us to better understand both the strengths and limitations of Hofmann's semantic argument by starting from a broader perspective. We also emphasise that the inductive notion of h-elementary type---see \autoref{simple type}---is not merely a technical tool that allows us to continue applying Hofmann's argument, but it is also a natural and expressive family of statements to consider in a theory. For instance, all statements in Heyting arithmetic are h-elementary statements.

\medskip

\noindent\textbf{Structure of the paper.} In \autoref{section2.2} we recall the notion of category with attributes, that constitutes a notion of semantics for dependent type theories. We recall a well-known notion of \textit{sound} semantics for extensional dependent type theories based on categories with attributes and we present an analogous one for propositional ones.

In \autoref{section2.3} we consider the notions of homotopy equivalence between types (namely \textit{homotopy equivalence}) and of homotopy equivalence between contexts (\textit{context homotopy equivalence}) in a propositional dependent type theory. These notions identify maps (between type and contexts, respectively) that are invertible up to homotopy (i.e. pointwise propositional equality and pointwise context propositional equality \cite{MR2469279} respectively). We recall how to extend a homotopy equivalence between two contexts via a homotopy equivalence between two types in the given contexts, and how to obtain new homotopy equivalences between types, by starting from other given ones and applying the type constructors. Finally, we inductively define a family of \textit{canonical} homotopy equivalences and a family of \textit{canonical} context homotopy equivalences, notions appearing in other works \cite{phdhofmann,hofmannconservativity,MAIETTI2009319,contentemaietti,maiettisabelli} mentioned in \autoref{section2.3.VI}. This notion is used in \autoref{section2.5} in order to make the syntax of a ``sub-theory'' of a propositional type theory into a model of an extensional type theory.

\autoref{section2.4} is devoted to studying properties of the family of the canonical context homotopy equivalences. These properties are deduced by induction on the complexity of such an equivalence: as the family of the canonical context homotopy equivalences depends on the family of the canonical homotopy equivalences between types, a result for the former usually follows by induction starting from a result for the latter, which is obtained by induction as well. The properties of reflexivity, symmetry and transitivity hold for the general family of the canonical context homotopy equivalences. However, the fundamental property that any two parallel canonical equivalences are homotopic, i.e. pointwise context propositionally equal, does not hold, hence we need to restrict ourselves to the smaller family of canonical equivalences between those contexts that we call \textit{contexts with h-propositional identities} (see \autoref{concrete context}). Having this property is fundamental to ensure the well-definedness of the composition in the category with attributes that we define in \autoref{section2.5}, starting from the syntax of the given propositional type theory. In detail, we use the property of \textit{having h-propositional identities} at the end of the proof of \autoref{only one canonical homotopy equivalence}.

\autoref{section2.5} defines a category with attributes $\C$ starting from the \textit{h-elementary} contexts of the given propositional type theory, identified up to canonical context homotopy equivalence. The h-elementary contexts are a special case of contexts with h-propositional identities that we define at the beginning of the section (see \autoref{simple context}). This restriction makes the natural family of display maps that we define for $\C$ into a \textit{cartesian} natural transformation, so that we obtain a model of the strict substitution. Specifically, the h-elementariness property is used to deduce the uniqueness of the factorisation of a given square against a naturality square. We conclude our analysis by proving in \autoref{section2.6} that $\C$ is actually a model of the extensional type theory. This allows us to infer, by soundness, the conservativity of the extensional theory over the h-elementary contexts of the propositional one.

\subsection{Preliminary conventions and terminology}\label{primitiveterminology}\label{section2.1.I}

In this paper, a dependent type theory is intended to be \textit{extensional} if it has extensional identity types, dependent product types, and dependent sum types, i.e. if it satisfies the rules of \autoref{extensional id}, \autoref{dependent pi}, and \autoref{dependent sigma} (see \autoref{inferencerules}) respectively---or, for dependent sums, their equivalent formulation in \autoref{dependent sigma char}. Analogously, a dependent type theory is said to be \textit{propositional} if it has propositional identity types, propositional dependent product types with the propositional function extensionality rule, and propositional dependent sum types, i.e. if it satisfies the rules of \autoref{propositional id}, \autoref{propositional pi}, and \autoref{propositional sigma} respectively---or, for propositional dependent sums, their equivalent formulation in \autoref{propositional sigma char}.

However, for sake of simplicity, for the remainder of the paper we consider \textit{one} propositional type theory: i.e. we assume that we are given a dependent type theory, which is propositional, together with a family of atomic types and atomic terms and subject to the usual structural rules---for an enumeration see \cite[Chapter III---context formation, context equality, judgement formation, judgement equality]{MR1134134}, \cite{MR1674451}, \cite{Hofmann1997}. We refer to this specific dependent type theory as Propositional Type Theory ($\ptt$). Analogously, as there is only one extensional type theory that we actually consider in this paper (relative to the given $\ptt$), i.e. the one whose atomic types are the atomic types of $\ptt$ that are provably \textit{h-sets} in $\ptt$ (and whose atomic terms are the atomic terms of these atomic h-sets in $\ptt$), we refer to it as Extensional Type Theory ($\ett$).

We use the symbology $\judge{-}{-}$ to indicate judgements of $\ptt$ and the symbology $\judgext{-}{-}$ to indicate the ones of $\ett$. We use the symbol $\equiv$ to indicate judgmental equalities between contexts, types in the same context and terms in the same context and of the same type. We use the symbol $=$ to indicate propositional equalities i.e. identity types. Sometimes (especially in diagrams) we adopt the notations $x\Rightarrow y$ and $x\Leftarrow y$, where $x\Rightarrow y$ indicates the type $x=y$ and $x\Leftarrow y$ indicates the type $y=x$.

Sometimes, we use path induction on paths that may not appear to be general, i.e. paths of the form $x,y:A, p:x=y$ for some type $A$. However, in these instances, we mean that the specific type on which we are performing path induction can be generalised so that the specific path $p'$ is replaced by a general path as $p$. This allows us to use path induction, after which we substitute $p'$ for the general path $p$ in order to obtain the desired result. This argument is indicated as \textit{generalised path induction} or \textit{generalised path elimination}.

\section{Recap on categories with attributes}\label{semantica via categorie con attributi}\label{section2.2}

This section is devoted to the notion of a \textit{category with attributes} and its use to define an opportune notion of semantics for several kinds of dependent type theories. For further details, we refer the reader to \cite{phdcartmell,moggi_1991, KAPULKIN2021106563}.

\begin{defi}[Category with attributes]\label{cwa}
Suppose that we are given:

\begin{itemize}
    \item A category $\mathcal{C}$ with terminal object $\1$, whose objects are called \textit{(semantic) contexts}.
    \item A functor $\mathcal{C}^{\op}\xrightarrow{\tp} \discr$, that we call \textit{presheaf of (semantics) types}. If $\Gamma$ is a semantic context, then an object $A$ of the category $\tp \Gamma$ is said to be a \textit{semantic type} in semantic context $\Gamma$. Here $\discr$ denotes the category of (small) discrete categories.
    \item A functor $\gr \tp \xrightarrow{-.-}\mathcal{C}$, that we call \textit{(semantic) context extension}. Here $\gr\tp$ denotes the Grothendieck construction associated to the presheaf $\tp$.
    \item A cartesian natural transformation: $$\Big(\gr\tp \xrightarrow{-.-}\mathcal{C}\Big)\xrightarrow{P} \Big(\gr\tp \xrightarrow{\pi}\mathcal{C}\Big)$$ where $\pi$ denotes the projection on the first component (we remind that a natural transformation is said to be \textit{cartesian} when its naturality squares are pullbacks). The $(\Gamma,A)$-component: $$\Gamma.A \xrightarrow{P_{A}} \Gamma$$ of $P$ is called \textit{display map} of $A$.
\end{itemize} Then we say that the quadruple $(\mathcal{C},\tp,-.-,P)$ is a \textit{category with attributes}.
\end{defi}

The notion of category with attributes constitutes a notion of semantics for dependent type theories with the usual structural rules \cite{MR1674451,Hofmann1997}\cite[Chapter III---context formation, context equality, judgement formation, judgement equality]{MR1134134}, hence with the usual \textit{strict} notion of substitution. In this case, semantic terms of a category with attributes are presented as sections of the display maps. A modification of  this concept, producing an equivalent notion of semantics, is used by other authors \cite{clairambault_dybjer_2014,dybjer1996internal,phdhofmann,Hofmann1997}, where one requires that, for every semantic type $A$, a set of \textit{semantic terms} is given. Other requirements imply that this set is in bijection with the sections of the corresponding display map $P_A$. In our case (see \autoref{notation and terminology}), the semantic terms of $A$ are defined to be the the sections of $P_A$ themselves: it turns out that our notion of category with attributes is the one of \cite{Hofmann1997} in the special case where the bijections between the semantic terms (of given semantic types) and the sections (of the corresponding display maps) are identities.

For an actually more general notion of semantics for dependent type theories we refer the reader to the notion of \textit{comprehension category}, see \cite{MR1201808}, or to the one of \textit{display map category}, see \cite{taylor1999practical, MR1674451,mossvonglehn2018}. In this case, in fact, a dependent type theory modelled by such a structure does not necessarily enjoy a strictly functorial notion of substitution, analogously to the case of locally cartesian closed categories (see \cite{MR0727078}) as pointed out by Hofmann \cite{hofmannlccc}. For more details on a pseudo-functorial notion of substitution in dependent type theories, we refer the reader to \cite{curien1993substitution}.

\begin{rem}[Notation and terminology]\label{notation and terminology}
Let $(\mathcal{C},\tp,-.-,P)$ be a category with attributes.

\begin{itemize}
    \item For the $(\Gamma,A)$-component of $P$, we do not explicitly write the dependence on $\Gamma$. This does not generate ambiguity because, whenever we consider a semantic type, its semantic context is always specified.
    
    \item If $f$ is a morphism of semantic contexts $\Delta\to\Gamma$, then we denote as $-[f]$ the action: $$\tp \Gamma\to \tp \Delta$$ of $\tp$ on $f$. Hence, if $A$ is a semantic type in $\Gamma$, then $A[f]$ denotes the action of $-[f]$ on $A$. However, for sake of clarity we sometimes use the general notation $\tp f$ in place of $-[f]$, and hence we write $(\tp f)A$ in place of $A[f]$. E.g. this happens in \autoref{section2.5} and \autoref{section2.6}.
    
    \item We remind that the objects of $\gr \tp$ are pairs $(\Gamma,A)$, where $\Gamma$ is a semantic context and $A$ is semantic type in $\Gamma$. An arrow $(\Delta,B)\to(\Gamma,A)$ of $\gr \tp$ is a pair $(f,\phi)$, where $f$ is an arrow $\Delta \to \Gamma$ and $\phi$ an arrow $B \to A[f]$. We refer to \cite{MR3093880} for further details on this notion. The category $\tp \Delta$ is discrete, therefore an arrow $(\Delta,B)\to(\Gamma,A)$ is nothing but a morphism of contexts $\Delta \to \Gamma$ such that $B$ is the type $A[f]$ in $\Delta$. As we got rid of the dependence on the $\phi$ component, we denote the image of such an arrow via the functor $-.-$ as: $$\Delta.A[f]\xrightarrow{f.A}\Gamma.A$$ and any naturality (pullback, since $P$ is cartesian) square of $P$ is of the form: \[\begin{tikzcd} \Delta.A[f] \arrow[d, "P_{A[f]}"'] \arrow[r, "f.A"] & \Gamma.A \arrow[d, "P_A"] \\ \Delta \arrow[r, "f"] & \Gamma \end{tikzcd}\] for some semantic context morphism $f$.
    
    \item Let $A$ be a semantic type in semantic context $\Gamma$. The sections $\Gamma \xrightarrow{a} \Gamma.A$ of its display map $\Gamma.A \xrightarrow{P_A}\Gamma$ are called \textit{(semantic) terms} of $A$ in $\Gamma$.
\end{itemize}
\end{rem}

In a given dependent type theory, whenever we are given contexts $\gamma : \Gamma$ and $\delta : \Delta$, a type $\judge{\delta : \Delta}{B(\delta) : \type}$ and a morphism of contexts $\judge{\gamma : \Gamma}{f(\gamma) : \Delta,B}$ (see the beginning of \autoref{section2.3} for both the notion of \textit{morphism of contexts} and the meaning of this notation), then we implicitly mean that we are given a judgement of the form: $$\judge{\gamma : \Gamma}{\overline{f}(\gamma) : B(P_{B(\delta)}f(\gamma))}$$ such that the substitution $f(\gamma): \Delta,B(\delta)$ decomposes as: $$P_{B(\delta)}f(\gamma) : \Delta,\;\overline{f}(\gamma) : B(P_{B(\delta)}f(\gamma))$$ where $P_{B(\delta)}$ denotes the substitution $\judge{\delta :\Delta,y:B(\delta)}{\delta:\Delta}$. Every category with attributes models this phenomenon:

\begin{prop}\label{morphism decomposition}
Let $\cattypes$ be a category with attributes. Let $f$ be a semantic context morphism $\Gamma \to \Delta.B$, where $B$ is a semantic type in $\Delta$. Then there is a unique semantic term $\overline{f}$ of $B[P_Bf]$ in $\Gamma$ such that the diagram: \[\begin{tikzcd} \Gamma \arrow[drrr, bend left, "f"] \arrow[dr, "\overline{f}"'] \\ & \Gamma.B[P_Bf] \arrow[rr, "P_Bf.B"'] && \Delta.B \end{tikzcd}\] commutes. Moreover, if $f$ is of the form: $$\Gamma \xrightarrow{b}\Gamma.B[g]\xrightarrow{g.B}\Delta.B$$ for some semantic context morphism $\Gamma \xrightarrow{g}\Delta$ and some semantic term $b$ of  $B[g]$, then the semantic term $\overline{f}$ is $b$ itself (observe that: $$P_Bf=P_B(g.B)b=gP_{B[g]}b=g$$ hence in this case $\Gamma.B[P_Bf]=\Gamma.B[g]$, so that the equality $\overline{f}=b$ typechecks). \proof By cartesianity of $P$.
\endproof
\end{prop}

\subsection{Extending substitution to terms}\label{section2.2.I}

In a category with attributes $\cattypes$ the notion of substitution for semantic types (given by the presheaf $\tp$) naturally extends to semantic terms. Suppose that we are given a morphism of semantic contexts $\Delta \xrightarrow{f}\Gamma$, a semantic type $A$ in $\Gamma$ and a semantic term $a$ of $A$. As the square: \[\begin{tikzcd} \Delta.A[f] \arrow[d, "P_{A[f]}"'] \arrow[r, "f.A"] & \Gamma.A \arrow[d, "P_A"] \\ \Delta \arrow[r, "f"] & \Gamma \end{tikzcd}\] is a pullback and $P_A a f = f$, there is a unique section $a[f]$ of $P_{A[f]}$ such that: \[\begin{tikzcd} \Delta \arrow[r, "f"] \arrow[d, "a\text{[}f\text{]}"'] & \Gamma \arrow[d, "a"]\\ \Delta.A[f] \arrow[r, "f.A"] & \Gamma.A \end{tikzcd}\] commutes. We define the action of the substitution $-[f]$ on $a$ as the semantic term $a[f]$ of $A[f]$. As expected:

\begin{prop}
For a semantic term $a$ of type $A$ in context $\Gamma$ and a semantic morphism $\Delta\xrightarrow{f}\Gamma$, the semantic terms $a[f]$ and $\overline{af}$ coincide.

\proof
It holds that $af=(f.A)a[f]$. Hence by \autoref{morphism decomposition} $\overline{af}$ is $a[f]$.
\endproof
\end{prop}

\noindent
and additionally:

\begin{prop}\label{functoriality of a}
For a semantic term $a$, the semantic term $a[f]$ is functorial in $f$.

\proof The term $a[1_\Gamma]$ of $A[1_\Gamma]=A$ makes the diagram: 
\[\begin{tikzcd}
	\Gamma && \Gamma \\
	\\
	{\Gamma.A} && {\Gamma.A}
	\arrow["{1_\Gamma}", from=1-1, to=1-3]
	\arrow["{a[1_{\Gamma}]}"{description}, from=1-1, to=3-1]
	\arrow["a"{description}, from=1-3, to=3-3]
	\arrow["{1_\Gamma.A=1_{\Gamma.A}}", from=3-1, to=3-3]
\end{tikzcd}\] commute, hence $a[1_\Gamma]=a$.

If we are given morphisms of contexts $\Omega \xrightarrow{g} \Delta \xrightarrow{f} \Gamma$ then the minimal squares in the diagram: 
\[\begin{tikzcd}
	\Omega && \Delta && \Gamma \\
	\\
	{\Omega.A[f][g]} && {\Delta.A[f]} && {\Gamma.A}
	\arrow["g", from=1-1, to=1-3]
	\arrow["{a[f][g]}"{description}, from=1-1, to=3-1]
	\arrow["f", from=1-3, to=1-5]
	\arrow["{a[f]}"{description}, from=1-3, to=3-3]
	\arrow["a"{description}, from=1-5, to=3-5]
	\arrow["{g.A[f]}", from=3-1, to=3-3]
	\arrow["{f.A}", from=3-3, to=3-5]
\end{tikzcd}\] commute, hence the outer rectangle---whose lower side is $(f.A)(g.A[f])=(fg).A$---commutes as well. As $a[fg]$ is the unique term of $A[fg]=A[f][g]$ making the outer rectangle commute, we conclude that $a[f][g]=a[fg]$. \endproof
\end{prop}

We end the current subsection recalling the following:

\begin{lem}\label{encode two terms in one morphism}
Suppose that $a$ and $b$ are semantic terms of some semantic type $A$ in some semantic context $\Gamma$. Then: $$a;b:=(\;\Gamma\xrightarrow{b}\Gamma.A \xrightarrow{a.A[P_A]}\Gamma.A.A[P_A]\;)=(\;\Gamma\xrightarrow{a}\Gamma.A\xrightarrow{b[P_A]}\Gamma.A.A[P_A]\;)$$ and $a;b$ is the unique arrow $\Gamma\to\Gamma.A.A[P_A]$ whose postcompositions via $P_{A[P_A]}$ and via $P_A.A$ are $a$ and $b$ respectively. Moreover, if $\Delta \xrightarrow{f}\Gamma$ is a morphism of contexts, then: \[\begin{tikzcd}
	\Delta &&& \Gamma \\
	\\
	{\Delta.A[f].A[f][P_{A[f]}]} & {\Delta.A[f].A[P_A(f.A)]} && {\Gamma.A.A[P_A]}
	\arrow["f", from=1-1, to=1-4]
	\arrow["{a[f];b[f]}"{description}, from=1-1, to=3-1]
	\arrow["{a;b}"{description}, from=1-4, to=3-4]
	\arrow["{f.A.A[P_A]}", from=3-2, to=3-4]
	\arrow[Rightarrow, no head, from=3-1, to=3-2]
\end{tikzcd}\] commutes.

\proof By building the term $b=b[P_Aa]$ according to the definition of $b[P_Aa]$ itself:
\[\begin{tikzcd}
	\Gamma && {\Gamma.A} && \Gamma \\
	\\
	{\Gamma.A=\Gamma.A[P_Aa]} && {\Gamma.A.A[P_A]} && {\Gamma.A} \\
	{} \\
	\Gamma && {\Gamma.A} && \Gamma
	\arrow["a", from=5-1, to=5-3]
	\arrow["{P_A}", from=5-3, to=5-5]
	\arrow["a", from=1-1, to=1-3]
	\arrow["{P_A}", from=1-3, to=1-5]
	\arrow["{P_A}"{description}, from=3-5, to=5-5]
	\arrow["{P_A=P_{A[P_Aa]}}"{description}, from=3-1, to=5-1]
	\arrow["b"{description}, from=1-5, to=3-5]
	\arrow["{b[P_A]}"{description}, from=1-3, to=3-3]
	\arrow["{b[P_Aa]}"{description}, from=1-1, to=3-1]
	\arrow["{a.A[P_A]}", from=3-1, to=3-3]
	\arrow["{P_A.A}", from=3-3, to=3-5]
	\arrow["{P_{A[P_A]}}"{description}, from=3-3, to=5-3]
	\arrow["\lrcorner"{anchor=center, pos=0.125}, draw=none, from=3-3, to=5-5]
	\arrow["\lrcorner"{anchor=center, pos=0.125}, draw=none, from=3-1, to=5-3]
\end{tikzcd}\] the upper left-hand square commutes, hence we are done with the first equality. Since $P_{A[P_A]}b[P_A]a=a$ and $(P_A.A)(a.A[P_A])b=((P_Aa).A)b=b$, we conclude that $a;b$ is the unique arrow $\Gamma\to\Gamma.A.A[P_A]$ whose postcompositions via $P_{A[P_A]}$ and via $P_A.A$ are $a$ and $b$ respectively.

In order to verify that the diagram of the statement commutes, we use the first presentation: $$\Delta\xrightarrow{b[f]}\Delta.A[f]\xrightarrow{a[f].A[f][P_{A[f]}]}\Delta.A[f].A[f][P_{A[f]}]$$ for $a[f];b[f]$, and the second one: $$\Gamma\xrightarrow{a}\Gamma.A\xrightarrow{b[P_A]}\Gamma.A.A[P_A]$$ for $a;b$. We are left to verify that: \[\begin{tikzcd}[column sep=scriptsize]
	\Delta &&&& \Gamma && {\Gamma.A} \\
	\\
	{\Delta.A[f]} &&& {\Delta.A[f].A[f][P_{A[f]}]} & {\Delta.A[f].A[P_A(f.A)]} && {\Gamma.A.A[P_A]}
	\arrow["{b[P_A]}"{description}, from=1-7, to=3-7]
	\arrow["{f.A.A[P_A]}", from=3-5, to=3-7]
	\arrow[Rightarrow, no head, from=3-4, to=3-5]
	\arrow["{b[P_A(af)]}"{description}, from=1-1, to=3-1]
	\arrow["{a[f].A[f][P_{A[f]}]}", from=3-1, to=3-4]
	\arrow["f", from=1-1, to=1-5]
	\arrow["a", from=1-5, to=1-7]
\end{tikzcd}\] commutes, since $b[P_A(af)]=b[f]$. Since the diagram: \[\begin{tikzcd}
	\Delta &&& \Gamma &&& {\Gamma.A} \\
	\\
	{\Delta.A[f]} &&& {\Gamma.A} &&& {\Gamma.A.A[P_A]}
	\arrow["{b[P_A]}"{description}, from=1-7, to=3-7]
	\arrow["{b[P_A(af)]}"{description}, from=1-1, to=3-1]
	\arrow["f", from=1-1, to=1-4]
	\arrow["a", from=1-4, to=1-7]
	\arrow["{f.A}", from=3-1, to=3-4]
	\arrow["{a.A[P_A]}", from=3-4, to=3-7]
\end{tikzcd}\] commutes (and this is true because the equality $(a.A[P_A])(f.A)=(af).A[P_A]$ holds) we are done if we verify that: \[\begin{tikzcd}
	{\Delta.A[f]} &&& {\Gamma.A} \\
	\\
	{\Delta.A[f].A[P_A(f.A)]} &&& {\Gamma.A.A[P_A]}
	\arrow["{a.A[P_A]}"{description}, from=1-4, to=3-4]
	\arrow["{f.A}", from=1-1, to=1-4]
	\arrow["{f.A.A[P_A]}", from=3-1, to=3-4]
	\arrow["{a[f].A[P_A(f.A)]}"{description}, from=1-1, to=3-1]
\end{tikzcd}\] commutes and use that $a[f].A[f][P_{A[f]}]=a[f].A[P_A(f.A)]$. This is actually true as: $$(\;f.A.A[P_A]\;)(\;a[f].A[P_A(f.A)]\;)=(\;(f.A)a[f]\;).A[P_A]=(\;af\;).A[P_A]=(\;a.A[P_A]\;)(\;f.A\;)$$ hence we are done.
\endproof
\end{lem}

\begin{rem}
If we are given a dependent type theory, if we are given judgements of the form $\judge{\gamma}{a(\gamma):A(\gamma)}$ and $\judge{\gamma}{b(\gamma):A(\gamma)}$, then the morphism of semantic contexts $a;b$ of \autoref{encode two terms in one morphism} corresponds to the morphism of contexts $\judge{\gamma}{\gamma,a(\gamma),b(\gamma)}$.
\end{rem}

\subsection{Variable terms}\label{variable term}\label{section2.2.II}

Let $\cattypes$ be a category with attributes and let $A$ be a semantic type in a semantic context $\Gamma$. Then, we might consider the semantic term $\overline{1_{\Gamma.A}}$ of $A[P_A]$, that is a section: $$\Gamma.A\xrightarrow{\overline{1_{\Gamma.A}}}\Gamma.A.A[P_A]$$ of $P_{A[P_A]}$. We denote this semantic term as $v_A$ and call it \textit{semantic variable} of $A$ and remind that it is characterised, among the sections of $P_{A[P_A]}$, as the one satisfying $(P_A.A)v_A=1_{\Gamma.A}$. We recall, without proof, some important equalities:

\begin{lem}\label{substitution in a variable term}
Let $\Delta\xrightarrow{f}\Gamma$ be an arrow in $\mathcal{C}$ and consider the corresponding extension $\Delta.A[f]\xrightarrow{f.A}\Gamma.A$. Then the equality: $$v_A[f.A]=v_{[A[f]]}$$ between semantic terms of $A[P_A(f.A)]=A[fP_{A[f]}]$ holds.
\end{lem}

\begin{lem}\label{variable&term}
If $a$ is a semantic term of $A$, then: $$(\;\Gamma \xrightarrow{a} \Gamma.A \xrightarrow{v_A} \Gamma.A.A[P_A]\;)=(\;\Gamma \xrightarrow{a} \Gamma.A \xrightarrow{a.A[P_A]} \Gamma.A.A[P_A]\;).$$ In particular, the equality $v_A[a]=a$ between sections of $P_A$ holds.
\end{lem}

\begin{rem}
In a given dependent type theory, if $\gamma:\Gamma$ is a context and $\judge{\gamma}{A(\gamma):\type}$ is a type, then two morphisms of contexts naturally arise: \[\begin{aligned}
\judge{\gamma,x:A(\gamma)&}{\gamma:\Gamma}\\
\judge{\gamma,x:A(\gamma)&}{x:A(\gamma)\equiv A(P_{A(\gamma)}(\gamma,x)).}
\end{aligned}\] The former corresponds to the morphism of contexts $P_A$ in a given category with attributes, while the latter, seen as a term rather then a general morphism of contexts, to the semantic term $v_A$ associated to a given semantic type $A$.
\end{rem}

\subsection{Semantics of extensional and propositional type theories}\label{section2.2.III}

In this subsection we briefly describe the additional structure that a category with attributes needs to be equipped with in order to allow on it an interpretation of an extensional type theory or of a propositional one.

\medskip

The following two notions, \autoref{semantic extensional identity types} and \autoref{semantic propositional identity types}, define when a category with attributes is equipped with semantic identity types (in the extensional and propositional case, respectively). In detail, the following notion is meant to model the inference rules for the extensional identity types, that here we recall in a concise form:
\begin{figure}[H]
\begin{center}
$\begin{alignedat}{2}
\text{Form\,}&\inferrule{\;\;\quad\quad\quad{A : \type}\quad\quad\quad\;\;}{\judge{x,x' : A}{x=x' : \type}}&\text{Extensionality\,}&\inferrule{A : \type}{\makebox[\widthof{$\judge{x, x' : A;\; p : x=x'}{p\equiv \refl(x)}$}][c]{$\judge{x, x' : A;\; p : x=x'}{x\equiv x'}$}}\\
\\
\text{Intro\,}&\inferrule{\;\;\quad\quad\quad{A : \type}\quad\quad\quad\;\;}{\judge{x : A}{\refl(x) : x=x}}&\;\;\;\quad\text{Id proof irrelevence\,}&\inferrule{A : \type}{\judge{x, x' : A;\; p : x=x'}{p\equiv \refl(x)}}\\
\end{alignedat}$
\end{center}
\end{figure}
\noindent referring the reader to \autoref{extensional id} for the extended version.

\begin{defi}[Semantic extensional identity types]\label{semantic extensional identity types} We say that a category with attributes $\cattypes$ is \textbf{equipped with semantic extensional identity types} if, for every semantic context $\Gamma$ and every semantic type $A$ in context $\Gamma$, there is a choice of:
\begin{itemize}
    \item (\textit{Formation}) a semantic type $\id_A$ of semantic context $\Gamma.A.A[P_A]$;
    \item (\textit{Introduction}) a morphism of contexts: $$\Gamma.A\xrightarrow{\refl_A}\Gamma.A.A[P_A].\id_A$$ such that the diagram: \[\begin{tikzcd}
	{\Gamma.A} && {\Gamma.A.A[P_A].\id_A} \\
	\\
	&& {\Gamma.A.A[P_A]}
	\arrow["{\refl_A}", from=1-1, to=1-3]
	\arrow["{P_{\id_A}}", from=1-3, to=3-3]
	\arrow["{v_A}"', from=1-1, to=3-3]
	\end{tikzcd}\] commutes;
\end{itemize} in such a way that the following properties are satisfied: \begin{itemize}
\item (\textit{Extensionality}) For all semantic terms $a$ and $b$ of semantic type $A$ and every semantic term $\Gamma\xrightarrow{p}\Gamma.\id_A(a,b)$ of semantic type: $$\id_A[a;b]=\id_A[\;\Gamma\xrightarrow{b}\Gamma.A \xrightarrow{a.A[P_A]}\Gamma.A.A[P_A]\;]=\id_A[\;\Gamma\xrightarrow{a}\Gamma.A\xrightarrow{b[P_A]}\Gamma.A.A[P_A]\;]$$ in context $\Gamma$, the equality $a=b$ between semantic terms of type $A$ and the equality: $$(\;\Gamma\xrightarrow{p}\Gamma.\id_A[a;b]=\Gamma.\id_A[a;a]\;)=(\;\Gamma\xrightarrow{\refl_A^a}\Gamma.\id_A[a;a]\;)$$ between semantic terms of type $\id_A[a;a]$ hold. Here $\refl_A^a$ is the unique arrow $\Gamma\dashrightarrow\Gamma.\id_A[a;a]$ such that: \[\begin{tikzcd}
	\Gamma &&&&& {\Gamma.A} \\
	\\
	& {\Gamma.\id_A[a;a]} &&&& {\Gamma.A.A[P_A].\id_A} \\
	\\
	& \Gamma && {\Gamma.A} && {\Gamma.A.A[P_A]}
	\arrow["a", from=5-2, to=5-4]
	\arrow["{a.A[P_A]}", from=5-4, to=5-6]
	\arrow["{P_{\id_A(a,a)}}"{description}, from=3-2, to=5-2]
	\arrow["{P_{\id_A}}"{description}, from=3-6, to=5-6]
	\arrow["{((a.A[P_A])a).\id_A}\;=\;a;a.\id_A", from=3-2, to=3-6]
	\arrow["{\refl_A}"{description}, from=1-6, to=3-6]
	\arrow["a", from=1-1, to=1-6]
	\arrow[Rightarrow, no head, from=1-1, to=5-2,  bend right]
	\arrow[dashed, from=1-1, to=3-2]
	\arrow["\lrcorner"{anchor=center, pos=0.125}, draw=none, from=3-2, to=5-4]
\end{tikzcd}\] commutes (observe that $P_{\id_A}\refl_Aa=v_Aa=(a.A[P_A])a$ by \autoref{variable&term}).
\item (\textit{Compatibility with the substitution}) If $f$ is a morphism of semantic contexts $\Delta\to\Gamma$ then the equality: $$\id_A[\;\Delta.A[f].A[fP_{A[f]}]\xrightarrow{f.A.A[P_A]}\Gamma.A.A[P_A]\;]=\id_{A[f]}$$ between semantic types in context $\Delta.A[f].A[fP_{A[f]}]$ holds and: $$\refl_A[f]=\refl_{A[f]}$$ where $\refl_A[f]$ is the unique arrow $\Delta.A[f]\to\Delta.A[f].A[f][P_{A[f]}].\id_{A[f]}$ such that $P_{\id_A}\refl_A[f]=v_A$ and such that the diagram: \[\begin{tikzcd}
	{\Delta.A[f]} &&& {\Gamma.A} \\
	\\
	{\Delta.A[f].A[f][P_{A[f]}].\id_{A[f]}} &&& {\Gamma.A.A[P_A].\id_A}
	\arrow["{\refl_{A}[f]}"{description}, from=1-1, to=3-1]
	\arrow["{f.A.A[P_A].\id_A}", from=3-1, to=3-4]
	\arrow["{\refl_A}"{description}, from=1-4, to=3-4]
	\arrow["{f.A}", from=1-1, to=1-4]
\end{tikzcd}\] commutes.
\end{itemize}
\end{defi}

The following notion is meant to model the inference rules for the propositional identity types, that here we recall in a concise form:
\begin{figure}[H]
\begin{center}
$\begin{alignedat}{2}
\text{Form\,}&\inferrule{\;\;\quad\quad\quad{A : \type}\quad\quad\quad\;\;}{\judge{x,x' : A}{x=x' : \type}}&\text{Elim\,}&\inferrule{A : \type\\\\\judge{x, x';\; p : x=x'}{C(x,x',p) : \type}\\\\\judge{x}{c(x) : C(x,x,\refl(x))}}{\judge{x, x';\; p}{\j(c,x,x',p) : C(x,x',p)}}\\
\\
\text{Intro\,}&\inferrule{\;\;\quad\quad\quad{A : \type}\quad\quad\quad\;\;}{\judge{x : A}{\refl(x) : x=x}}&\;\;\;\quad\text{Prop comp\,}&\inferrule{A : \type\\\\\judge{x, x';\; p : x=x'}{C(x,x',p) : \type}\\\\ \judge{x}{c(x) : C(x,x,\refl(x))}}{\judge{x}{\h(c,x):\j(c,x,x,\refl(x))= c(x)}}\\
\end{alignedat}$
\end{center}
\end{figure}
\noindent referring the reader to Figure \autoref{propositional id} for the extended version.

\begin{defi}[Semantic propositional identity types]\label{semantic propositional identity types}
We say that a category with attributes $\cattypes$ is \textbf{equipped with semantic propositional identity types} if it satisfies \textit{formation}, \textit{introduction}, \textit{compatibility with the substitution} of \autoref{semantic extensional identity types} and moreover:
\begin{itemize}
    \item (\textit{Elimination and propositional computation}) for every semantic context $\Gamma$, every semantic type $A$ in context $\Gamma$, every semantic type $C$ in context $\Gamma.A.A[P_A].\id_A$ and every semantic term $\Gamma.A\xrightarrow{c}\Gamma.A.C[\refl_A]$ of $C[\refl_A]$ in context $\Gamma.A$, there is a choice of a semantic term: $$\Gamma.A.A[P_A].\id_A\xrightarrow{\j_c}\Gamma.A.A[P_A].\id_A.C$$ of type $C$ in context $\Gamma.A.A[P_A].\id_A$ and of a semantic term: $$\Gamma.A\xrightarrow{\h_c}\Gamma.A.\id_{C[\refl_A]}[\j_c[\refl_A];c]$$ of type $\id_{C[\refl_A]}[\j_c[\refl_A];c]$ in context $\Gamma.A$;
    \item (\textit{Additional compatibility with the substitution}) for every semantic context $\Gamma$, every semantic type $A$ in context $\Gamma$, every semantic type $C$ in context $\Gamma.A.A[P_A].\id_A$, every semantic term $\Gamma.A\xrightarrow{c}\Gamma.A.C[\refl_A]$ of $C[\refl_A]$ in context $\Gamma.A$ and every morphism of semantic contexts $\Delta\xrightarrow{f}\Gamma$, the diagram: \[\begin{tikzcd}[cramped, column sep=tiny]
	{\Delta.A[f].A[fP_{A[f]}].\id_{A[f]}} & {\Delta.A[f].A[P_A(f.A)].\id_A[f.A.A[P_A]]} \\
	\\
	\\
	{\Delta.A[f].A[fP_{A[f]}].\id_{A[f]}.C'} & {\Delta.A[f].A[P_A(f.A)].\id_A[f.A.A[P_A]].C'}
	\arrow[Rightarrow, no head, from=1-1, to=1-2]
	\arrow[Rightarrow, no head, from=4-1, to=4-2]
	\arrow["{\j_{\big(\Delta.A[f]\xrightarrow{c[f.A]}\Delta.A[f].C[\refl_A(f.A)]=\Delta.A[f].C[(f.A.A[P_A].\id_A)\refl_{A[f]}]\big)}}", from=1-1, to=4-1, shift right=15]
	\arrow["{\j_c[f.A.A[P_A].\id_A]}", from=1-2, to=4-2, shift left=15]
\end{tikzcd}\] where $C'\equiv C[f.A.A[P_A].\id_A]$, commutes i.e. the equality $\j_{c[f.A]}=\j_c[f.A.A[P_A].\id_A]$ holds, and that the diagram: 
\[\begin{tikzcd}[row sep=small]
	{\Delta.A[f]} && {\Delta.A[f].\id_{C[(f.A.A[P_A].\id_A)\refl_{A[f]}]}[\;\j_{c[f.A]}[\refl_{A[f]}];c[f.A]\;]} \\
	&& {\Delta.A[f].\id_{C[\refl_A][f.A]}[\;\j_{c}[f.A.A[P_A].\id_A][\refl_{A[f]}];c[f.A]\;]} \\
	&& {\Delta.A[f].\id_{C[\refl_A][f.A]}[\;\j_{c}[\refl_A][f.A];c[f.A]\;]} \\
	&& {\Delta.A[f].(\;\id_{C[\refl_A]}[\;f.A.C[\refl_A].C[\refl_AP_{C[\refl_A]}]\;]\;)[\;\j_{c}[\refl_A][f.A];c[f.A]\;]} \\
	{\Delta.A[f]} && {\Delta.A[f].(\;\id_{C[\refl_A]}[\j_c[\refl_A];c]\;)[f.A]}
	\arrow[Rightarrow, no head, from=1-3, to=2-3]
	\arrow[Rightarrow, no head, from=2-3, to=3-3]
	\arrow[Rightarrow, no head, from=3-3, to=4-3]
	\arrow["{\text{ \autoref{encode two terms in one morphism}}}", Rightarrow, no head, from=4-3, to=5-3]
	\arrow["{\h_c[f.A]}"', from=5-1, to=5-3]
	\arrow["{\h_{c[f.A]}}", from=1-1, to=1-3]
	\arrow[Rightarrow, no head, from=1-1, to=5-1]
\end{tikzcd}\] commutes i.e. the equality $\h_{c[f.A]}=\h_c[f.A]$ holds.
\end{itemize}
\end{defi}

The following two notions (\autoref{semantic dependent product types} and \autoref{semantic propositional dependent product types}) define when a category with attributes is equipped with semantic dependent product types (in the extensional and propositional case, respectively). In detail, the following notion is meant to model the inference rules for the dependent product types, that here we recall in a concise form:
\begin{figure}[H]
\begin{center}$\begin{alignedat}{2}
\text{Form\,}&\inferrule{\;\;\;{A : \type}\quad\quad\judge{x:A}{B(x):\type}\;\;\;}{{\Pi_{x:A}B(x) : \type}}&\quad\text{Intro\,}&\inferrule{{A : \type}\quad\quad\judge{x:A}{B(x):\type}\\\\\judge{x:A}{y(x):B(x)}}{{\lambda x.y(x):\Pi_{x:A}B(x)}}\\
\\
\text{Elim\,}&\inferrule{\;\;\;{A : \type}\quad\quad\judge{x:A}{B(x):\type}\;\;\;}{\judge{z:\Pi_{x:A}B(x);\;x:A}{\ev(z,x):B(x)}}&\quad\quad\text{Comp\,}&\inferrule{{A : \type}\quad\quad\judge{x:A}{B(x):\type}\\\\\judge{x:A}{y(x):B(x)}}{\judge{x:A}{\ev(\lambda x.y(x),x)\equiv y(x)}}\\
\end{alignedat}$\end{center}
\end{figure}
\noindent referring the reader to \autoref{dependent pi} for the extended version. Here, we also present the semantic counterpart of the expansion rule, which follows automatically in the presence of extensional identity types.

\begin{defi}[Semantic dependent product types]\label{semantic dependent product types}
We say that a category with attributes $\cattypes$ is \textbf{equipped with semantic dependent product types} if:
\begin{itemize}
    \item (\textit{Formation}) for every semantic context $\Gamma$, every semantic type $A$ in context $\Gamma$ and every semantic type $B$ in context $\Gamma.A$, there is a choice of a semantic type $\pid_A^B$;
    \item (\textit{Introduction}) for every semantic context $\Gamma$, every semantic type $A$ in context $\Gamma$, every semantic type $B$ in context $\Gamma.A$ and every semantic term $\Gamma.A\xrightarrow{b}\Gamma.A.B$ of $B$, there is a choice of a semantic term $\Gamma\xrightarrow{\lambda b}\Gamma.\pid_A^B$ of $\pid_A^B$;
    \item (\textit{Elimination}) for every semantic context $\Gamma$, every semantic type $A$ in context $\Gamma$, every semantic type $B$ in context $\Gamma.A$, every semantic term $\Gamma\xrightarrow{z}\Gamma.\pid_A^B$ of $\pid_A^B$ and every semantic term $\Gamma\xrightarrow{a}\Gamma.A$ of $A$, a choice of a semantic term $\Gamma\xrightarrow{\ev_z^a}\Gamma.B[a]$ of $B[a]$;
\end{itemize} in such a way that the following properties are satisfied:
\begin{itemize}
    \item (\textit{Compatibility with the substitution}) For every semantic context $\Gamma$, every semantic type $A$ in context $\Gamma$ and every semantic type $B$ in context $\Gamma.A$, and for every choice of a semantic term $\Gamma\xrightarrow{z}\Gamma.\pid_A^B$ of $\pid_A^B$, of a semantic term $\Gamma\xrightarrow{a}\Gamma.A$ of $A$ and of a semantic term $\Gamma.A\xrightarrow{b}\Gamma.A.B$ of $B$, if $\Delta\xrightarrow{f}\Gamma$ is a morphism of semantic contexts, then: \begin{itemize}
        \item the equality: $$\pid_A^B[f]=\pid_{A[f]}^{B[f.A]}$$ between semantic types in context $\Gamma$ holds;
        \item the equality: $$\bigg(\ev_{\Gamma\xrightarrow{z}\Gamma.\pid_A^B}^{\Gamma\xrightarrow{a}\Gamma.A}\bigg)[f]=\bigg(\ev_{\Gamma\xrightarrow{z[f]}\Gamma.\pid_{A[f]}^{B[f.A]}}^{\Gamma\xrightarrow{a[f]}\Gamma.A[f]}\bigg)$$ between semantic terms of type $B[\;\Delta\xrightarrow{f}\Gamma\xrightarrow{a}\Gamma.A\;]=B[\;\Delta\xrightarrow{a[f]}\Delta.A[f]\xrightarrow{f.A}\Gamma.A\;]$ holds;
        \item the equality: $$(\;\lambda (\;\Gamma.A\xrightarrow{b}\Gamma.A.B\;)\;)[f] = \lambda (\;\Delta.A[f]\xrightarrow{b[f.A]}\Delta.A[f].B[f.A]\;)$$ between semantic terms of type $\pid_A^B[f]=\pid_{A[f]}^{B[f.A]}$ holds.
    \end{itemize}
    
    \item (\textit{Computation}) If $\Gamma.A\xrightarrow{b}\Gamma.A.B$ is a semantic term of $B$ and if $\Gamma\xrightarrow{a}\Gamma.A$ is a semantic term of $A$, for some semantic context $\Gamma$, some semantic type $A$ in context $\Gamma$ and some semantic type $B$ in context $\Gamma.A$, then the equality: $$(\;\Gamma\xrightarrow{\ev_{\lambda b}^a}\Gamma.B[a]\;)=(\;\Gamma\xrightarrow{b[a]}\Gamma.B[a]\;)$$ between semantic terms of type $B[a]$ holds.
    \item (\textit{Expansion}) If $\Gamma\xrightarrow{z}\Gamma.\pid_A^B$ is a semantic term of type $\pid_A^B$, then the equality: $$\lambda(\;\Gamma.A\xrightarrow{\Bigg(\ev_{\big(\;\Gamma.A\xrightarrow{z[P_A]}\Gamma.A.\pid_A^B[P_A]=\Gamma.A.\pid_{A[P_A]}^{B[P_A.A]}\;\big)}^{\big(\;\Gamma.A\xrightarrow{v_A}\Gamma.A.A[P_A]\;\big)}\Bigg)}\Gamma.A.B[\;(P_A.A)v_A\;]=\Gamma.A.B\;)=z$$ between semantic terms of type $\pid_A^B$ holds.
\end{itemize}
\end{defi}

The following notion is meant to model the inference rules for the propositional dependent product types, that here we recall in a concise form:
\begin{figure}[H]
\begin{center}
$\begin{alignedat}{2}
\text{Form\,}&\inferrule{{A : \type}\\\\\quad\judge{x:A}{B(x):\type}\quad}{{\Pi_{x:A}B(x) : \type}}&\quad\text{Intro\,}&\inferrule{{A : \type}\quad\quad\judge{x:A}{B(x):\type}\\\\\judge{x:A}{y(x):B(x)}}{{\lambda x.y(x):\Pi_{x:A}B(x)}}\\
\\
\text{Elim\,}&\inferrule{{A : \type}\\\\\quad\judge{x:A}{B(x):\type}\quad}{\judgectx{z:\Pi_{x:A}B(x);\;x:A}\\\\{\ev(z,x):B(x)}}&\quad\quad\text{Prop comp\,}&\inferrule{{A : \type}\quad\quad\judge{x:A}{B(x):\type}\\\\\judge{x:A}{y(x):B(x)}}{\judgectx{x:A}\\\\{\beta^\Pi(y,x):\ev(\lambda x.y(x),x)=y(x)}}\\
\\
\\
\\
\\
&&\quad\text{Intro\,}&\inferrule{{A : \type}\quad\quad\judge{x:A}{B(x):\type}}{\judgectx{z,z';\;p:\Pi_{x:A}\ev(z,x)=\ev(z',x)}\\\\{\funext(z,z',p):z=z'}}\\
\\
\text{Prop exp\,}&\inferrule{{A : \type}\\\\\quad\judge{x:A}{B(x):\type}\quad}{\boldsymbol{\lfloor}z,z';\;q:z=z'\boldsymbol{\rfloor}\\\\\eta_\funext(z,z',q):q=\\\\\funext(z,z',\lambda x.\ev(q,x))}&\quad\text{Prop comp\,}&\inferrule{{A : \type}\quad\quad\judge{x:A}{B(x):\type}\\\\\text{ }}{\judgectx{z,z';\;p:\Pi_{x:A}\ev(z,x)=\ev(z',x)}\\\\\beta_\funext(z,z',p):\\\\\lambda x.\ev(\funext(z,z',p),x)=p}
\end{alignedat}$\end{center}
\end{figure}
\noindent referring the reader to \autoref{propositional pi} and for the extended version. We do not explicitly write down the semantics of the propositional extensionality rules, but they can be formulated
analogously---we refer the reader to \cite{spadetto25} for additional details. Again, we present the semantic counterpart of the propositional expansion rule, which follows automatically from propositional function extensionality---syntactically, one defines: $$\eta^\Pi(z)\equiv\funext(\;\lambda x.\ev(z,x),z,\lambda x.(\;\beta^\Pi(\ev(z,-),x)^{-1}\;)\;):z=\lambda x.\ev(z,x)$$ in context $z:\Pi_{x:A}B(x)$. We will also refer to the computation and expansion rules as $\beta$-reduction and $\eta$-expansion, respectively.

\begin{defi}[Semantic propositional dependent product types]\label{semantic propositional dependent product types} We say that a category with attributes $\cattypes$ is \textbf{equipped with semantic propositional dependent product types} if it is equipped with semantic propositional identity types, it satisfies \textit{formation}, \textit{introduction}, \textit{elimination}, \textit{compatibility with the substitution} of \autoref{semantic dependent product types} and moreover:
\begin{itemize}
    \item (\textit{Propositional computation and additional compatibility with the substitution}) for every semantic context $\Gamma$, every semantic type $A$ in context $\Gamma$, every semantic type $B$ in context $\Gamma.A$, every semantic term $\Gamma.A\xrightarrow{b}\Gamma.A.B$ of $B$ and every semantic term $\Gamma\xrightarrow{a}\Gamma.A$ of $A$, there is a choice of a semantic term: $$\Gamma\xrightarrow{\beta_b^a}\Gamma.\id_{B[a]}[\;\ev_{\lambda b}^a;b[a]\;]$$ of type $\id_{B[a]}[\;\ev_{\lambda b}^a;b[a]\;]$, in such a way that the diagram: 
\[\begin{tikzcd}[row sep=small]
	\Delta && {\Delta.\id_{B[a]}[\;\ev_{\lambda b}^a;b[a]\;][f]} \\
	&& {\Delta.\id_{B[a]}[\;f.B[a].B[aP_{B[a]}]\;][\;\ev_{\lambda b}^a[f];b[a][f]\;]} \\
	&& {\Delta.\id_{B[a][f]}[\;\ev_{\lambda b}^a[f];b[a][f]\;]} \\
	&& {\Delta.\id_{B[f.A][a[f]]}[\;\ev_{\lambda b}^a[f];b[f.A][a[f]]\;]} \\
	&& {\Delta.\id_{B[f.A][a[f]]}[\;\ev_{(\lambda b)[f]}^{a[f]};b[f.A][a[f]]\;]} \\
	\Delta && {\Delta.\id_{B[f.A][a[f]]}[\;\ev_{\lambda (b[f.A])}^{a[f]};b[f.A][a[f]]\;]}
	\arrow["{\text{ \autoref{encode two terms in one morphism}}}", Rightarrow, no head, from=1-3, to=2-3]
	\arrow[Rightarrow, no head, from=2-3, to=3-3]
	\arrow[Rightarrow, no head, from=3-3, to=4-3]
	\arrow[Rightarrow, no head, from=4-3, to=5-3]
	\arrow[Rightarrow, no head, from=5-3, to=6-3]
	\arrow["{\beta_b^a[f]}", from=1-1, to=1-3]
	\arrow["{\beta_{b[f.A]}^{a[f]}}", from=6-1, to=6-3]
	\arrow[Rightarrow, no head, from=1-1, to=6-1]
	\end{tikzcd}\] commutes i.e. $\beta_b^a[f]=\beta_{b[f.A]}^{a[f]}$, for every morphism of contexts $\Delta\xrightarrow{f}\Gamma$;
    \item (\textit{Propositional expansion and additional compatibility with the substitution}) for every semantic context $\Gamma$, every semantic type $A$ in context $\Gamma$, every semantic type $B$ in context $\Gamma.A$ and every semantic term $\Gamma\xrightarrow{z}\Gamma.\pid_A^B$ of $\pid_A^B$, there is a choice of a semantic term: $$\Gamma\xrightarrow{\eta_z}\Gamma.\id_{\pid_A^B}[\;\lambda\ev_{z[P_A]}^{v_A};z\;]$$ of type $\id_{\pid_A^B}[\;\lambda\ev_{z[P_A]}^{v_A};z\;]$, in such a way that the diagram: 
    \[\begin{tikzcd}[row sep=small]
	\Delta && {\Delta.\id_{\pid_{A}^{B}}[\;\lambda\ev_{z[P_A]}^{v_A};z\;][f]} \\
	&& {\Delta.\id_{\pid_A^B}[\;f.\pid_A^B.\pid_A^B[P_{\pid_A^B}]\;][\;(\lambda\ev_{z[P_A]}^{v_A})[f];z[f]\;]} \\
	&& {\Delta.\id_{\pid_A^B[f]}[\;(\lambda\ev_{z[P_A]}^{v_A})[f];z[f]\;]} \\
	&& {\Delta.\id_{\pid_{A[f]}^{B[f.A]}}[\;\lambda(\ev_{z[P_A]}^{v_A}[f.A]);z[f]\;]} \\
	&& {\Delta.\id_{\pid_{A[f]}^{B[f.A]}}[\;\lambda(\ev_{z[P_A][f.A]}^{v_A[f.A]});z[f]\;]} \\
	\Delta && {\Delta.\id_{\pid_{A[f]}^{B[f.A]}}[\;\lambda(\ev_{z[f][P_{A[f]}]}^{v_{A[f]}});z[f]\;]}
	\arrow["{\text{ \autoref{encode two terms in one morphism}}}", Rightarrow, no head, from=1-3, to=2-3]
	\arrow[Rightarrow, no head, from=2-3, to=3-3]
	\arrow[Rightarrow, no head, from=3-3, to=4-3]
	\arrow[Rightarrow, no head, from=4-3, to=5-3]
	\arrow["{\text{ \autoref{substitution in a variable term}}}", Rightarrow, no head, from=5-3, to=6-3]
	\arrow["{\eta_z[f]}", from=1-1, to=1-3]
	\arrow["{\eta_{z[f]}}", from=6-1, to=6-3]
	\arrow[Rightarrow, no head, from=1-1, to=6-1]
	\end{tikzcd}\] commutes i.e. $\eta_z[f]=\eta_{z[f]}$, for every morphism of contexts $\Delta\xrightarrow{f}\Gamma$.
\end{itemize}
\end{defi}

The following two notions (\autoref{semantic dependent sum types} and \autoref{semantic propositional dependent sum types}) define when a category with attributes is equipped with semantic dependent sum types (in the extensional and propositional case, respectively). In detail, the following notion is meant to model the inference rules for the dependent sum types, that here we recall in a concise form:
\begin{figure}[H]
\begin{center}
$\begin{alignedat}{2}
\text{Form\,}&\inferrule{{A : \type}\\\\\judge{x:A}{B(x):\type}}{{\Sigma_{x:A}B(x) : \type}}&\quad\text{Elim\,}&\inferrule{{A : \type}\quad\quad\judge{x:A}{B(x):\type}\\\\ \judge{u:\Sigma_{x:A}B(x)}{C(u):\type} \\\\ \judge{x:A;\;y:B(x)}{c(x,y):C(\langle x,y\rangle)}}{\judge{u:\Sigma_{x:A}B(x)}{\splitt(c,u):C(u)}}\\
\\
\text{Intro\,}&\inferrule{{A : \type}\\\\\judge{x:A}{B(x):\type}}{\judgectx{x:A;\;y:B(x)}\\\\{\langle x,y\rangle:\Sigma_{x:A}B(x)}}&\quad\quad\quad\text{Comp\,}&\inferrule{{A : \type}\quad\quad\judge{x:A}{B(x):\type}\\\\ \judge{u:\Sigma_{x:A}B(x)}{C(u):\type} \\\\ \judge{x:A;\;y:B(x)}{c(x,y):C(\langle x,y\rangle)}}{\judgectx{x:A;\;y:B(x)}\\\\{\splitt(c,\langle x,y\rangle)\equiv c(x,y)}}\\
\end{alignedat}$
\end{center}
\end{figure}
\noindent referring the reader to \autoref{dependent sigma} for the extended version.

\begin{defi}[Semantic dependent sum types]\label{semantic dependent sum types}
A category with attributes $\cattypes$ is \textbf{equipped with semantic dependent sum types} if:
\begin{itemize}
    \item (\textit{Formation and introduction}) for every semantic context $\Gamma$, every semantic type $A$ in context $\Gamma$ and every semantic type $B$ in context $\Gamma.A$, there is a choice of a semantic type $\sigmad_A^B$ and of a morphism of contexts $\Gamma.A.B\xrightarrow{\pair_A^B}\Gamma.\sigmad_A^B$ such that: \[\begin{tikzcd}
	{\Gamma.A.B} && {\Gamma.\sigmad_A^B} \\
	\\
	{\Gamma.A} && \Gamma
	\arrow["{P_B}"{description}, from=1-1, to=3-1]
	\arrow["{P_A}", from=3-1, to=3-3]
	\arrow["{P_{\sigmad_A^B}}"{description}, from=1-3, to=3-3]
	\arrow["{\pair_A^B}", from=1-1, to=1-3]
	\end{tikzcd}\]
	commutes;
    \item (\textit{Elimination}) for every semantic context $\Gamma$, every semantic type $A$ in context $\Gamma$, every semantic type $B$ in context $\Gamma.A$, every semantic type $C$ in context $\Gamma.\sigmad_A^B$ and every semantic term $c$ of type $C[\pair_A^B]$ in context $\Gamma.A.B$, there is a choice of a semantic term $\splitt_c$ of type $C$ in context $\Gamma.\sigmad_A^B$;
\end{itemize} in such a way that the following properties are satisfied:
\begin{itemize}
    \item (\textit{Computation}) For every semantic context $\Gamma$, every semantic type $A$ in context $\Gamma$, every semantic type $B$ in context $\Gamma.A$, every semantic type $C$ in context $\Gamma.\sigmad_A^B$ and every semantic term $c$ of type $C[\pair_A^B]$ in context $\Gamma.A.B$, the equality: $$(\;\Gamma.A.B\xrightarrow{\splitt_c[\pair_A^B]}\Gamma.A.B.C[\pair_A^B]\;)=(\;\Gamma.A.B\xrightarrow{c}\Gamma.A.B.C[\pair_A^B]\;)$$ between semantic terms of type $C[\pair_A^B]$ holds.
    \item (\textit{Compatibility with the substitution}) For every semantic context $\Gamma$, every semantic type $A$ in context $\Gamma$, every semantic type $B$ in context $\Gamma.A$, every semantic type $C$ in context $\Gamma.\sigmad_A^B$, every semantic term $c$ of type $C[\pair_A^B]$ in context $\Gamma.A.B$, if $\Delta\xrightarrow{f}\Gamma$ is a morphism of semantic contexts, then:
    \begin{itemize}
        \item the equality: $$\sigmad_A^B[f]=\sigmad_{A[f]}^{B[f.A]}$$ between semantic types in context $\Gamma$ holds;
        \item the diagram: \[\begin{tikzcd}
	{\Delta.A[f].B[f.A]} &&& {\Gamma.A.B} \\
	\\
	{\Delta.\sigmad_{A[f]}^{B[f.A]}} & {\Delta.\sigmad_A^B[f]} && {\Gamma.\sigmad_A^B}
	\arrow["{\pair_A^B}"{description}, from=1-4, to=3-4]
	\arrow["{f.A.B}", from=1-1, to=1-4]
	\arrow[Rightarrow, no head, from=3-1, to=3-2]
	\arrow["{f.\sigmad_A^B}", from=3-2, to=3-4]
	\arrow["{\pair_{A[f]}^{B[f.A]}}"{description}, from=1-1, to=3-1]
\end{tikzcd}\] commutes;
        \item the equality between the semantic terms: $$\begin{aligned} \Delta.\sigmad_A^B[f]&\xrightarrow{\splitt_c[f.\sigmad_A^B]}\Delta.\sigmad_A^B[f].C[f.\sigmad_A^B]\\
        \Delta.\sigmad_{A[f]}^{B[f.A]}&\xrightarrow{\splitt_{c[f.A.B]}} \Delta.\sigmad_{A[f]}^{B[f.A]}.C[f.\sigmad_{A}^{B}] \end{aligned}$$ of semantic type $C[f.\sigmad_A^B]$ in semantic context $\Delta.\sigmad_{A}^{B}[f]=\Delta.\sigmad_{A[f]}^{B[f.A]}$ holds, where we remind that $c[f.A.B]$ is a term: $$\Delta.A[f].B[f.A]\;\to\;\Delta.A[f].B[f.A].\big(\;C\big[\pair_A^B(f.A.B)\big]\;=\;C\big[f.\sigmad_A^B\big]\big[\pair_{A[f]}^{B[f.A]}\big]\;\big).$$
    \end{itemize}
\end{itemize}
\end{defi}

The following notion is meant to model the inference rules for the propositional dependent sum types, that here we recall in a concise form:
\begin{figure}[H]
\begin{center}
$\begin{alignedat}{2}
\text{Form\,}&\inferrule{{A : \type}\\\\\judge{x:A}{B(x):\type}}{{\Sigma_{x:A}B(x) : \type}}&\quad\text{Elim\,}&\inferrule{{A : \type}\quad\quad\judge{x:A}{B(x):\type}\\\\ \judge{u:\Sigma_{x:A}B(x)}{C(u):\type} \\\\ \judge{x:A;\;y:B(x)}{c(x,y):C(\langle x,y\rangle)}}{\judge{u:\Sigma_{x:A}B(x)}{\splitt(c,u):C(u)}}\\
\\
\text{Intro\,}&\inferrule{{A : \type}\\\\\judge{x:A}{B(x):\type}}{\judgectx{x:A;\;y:B(x)}\\\\{\langle x,y\rangle:\Sigma_{x:A}B(x)}}&\quad\quad\text{Prop comp\,}&\inferrule{{A : \type}\quad\quad\judge{x:A}{B(x):\type}\\\\ \judge{u:\Sigma_{x:A}B(x)}{C(u):\type} \\\\ \judge{x:A;\;y:B(x)}{c(x,y):C(\langle x,y\rangle)}}{\judgectx{x:A;\;y:B(x)}\\\\{\sigma(c,x,y):\splitt(c,\langle x,y\rangle)= c(x,y)}}\\
\end{alignedat}$
\end{center}
\end{figure}
\noindent referring the reader to \autoref{propositional sigma} for the extended version.

\begin{defi}[Semantic propositional dependent sum types]\label{semantic propositional dependent sum types}
We say that $\cattypes$ is \textbf{equipped with semantic propositional dependent sum types} if it is equipped with semantic propositional identity types, it satisfies \textit{formation and introduction}, \textit{elimination}, \textit{compatibility with the substitution} of \autoref{semantic dependent sum types} and moreover:
\begin{itemize}
    \item (\textit{Propositional computation}) for every semantic context $\Gamma$, every semantic type $A$ in context $\Gamma$, every semantic type $B$ in context $\Gamma.A$, every semantic type $C$ in context $\Gamma.\sigmad_A^B$ and every semantic term $c$ of type $C[\pair_A^B]$ in context $\Gamma.A.B$, there is a choice of a semantic term: $$\Gamma.A.B\xrightarrow{\sigma_c} \Gamma.A.B.\id_{C[\pair_A^B]}[\;\splitt_c[\pair_A^B];c\;]$$ of type $\id_{C[\pair_A^B]}[\;\splitt_c[\pair_A^B];c\;]$;
    \item (\textit{Additional compatibility with the substitution}) for every semantic context $\Gamma$, every semantic type $A$ in context $\Gamma$, every semantic type $B$ in context $\Gamma.A$, every semantic type $C$ in context $\Gamma.\sigmad_A^B$, every semantic term $c$ of type $C[\pair_A^B]$ in context $\Gamma.A.B$ and every morphism of semantic contexts $\Delta\xrightarrow{f}\Gamma$, the following diagram: 
    \[\begin{tikzcd}[cramped, column sep=tiny, row sep=small]
	{\Delta.A[f].B[f.A]} & {\Delta.A[f].B[f.A].\id_{C[f.\sigmad_A^B][\pair_{A[f]}^{B[f.A]}]}[\;\splitt_{c[f.A.B]}[\pair_{A[f]}^{B[f.A]}];c[f.A.B]\;]} \\
	& {\Delta.A[f].B[f.A].\id_{C[\pair_{A}^{B}][f.A.B]}[\;\splitt_{c[f.A.B]}[\pair_{A[f]}^{B[f.A]}];c[f.A.B]\;]} \\
	& {\Delta.A[f].B[f.A].\id_{C[\pair_{A}^{B}][f.A.B]}[\;\splitt_{c}[f.\sigmad_A^B][\pair_{A[f]}^{B[f.A]}];c[f.A.B]\;]} \\
	& {\Delta.A[f].B[f.A].\id_{C[\pair_{A}^{B}][f.A.B]}[\;\splitt_{c}[\pair_A^B][f.A.B];c[f.A.B]\;]} \\
	& {\Delta.A[f].B[f.A].\id_{C[\pair_{A}^{B}]}[\;f.A.B.C[\pair_{A}^{B}].C'} \\
	{\Delta.A[f].B[f.A]} & {\Delta.A[f].B[f.A].\id_{C[\pair_{A}^{B}]}[\;\splitt_{c}[\pair_A^B];c\;][f.A.B]}
	\arrow["{\text{ \autoref{encode two terms in one morphism}}}", Rightarrow, no head, from=5-2, to=6-2]
	\arrow[Rightarrow, no head, from=4-2, to=5-2]
	\arrow[Rightarrow, no head, from=3-2, to=4-2]
	\arrow[Rightarrow, no head, from=2-2, to=3-2]
	\arrow[Rightarrow, no head, from=1-2, to=2-2]
	\arrow["{\sigma_{c[f.A.B]}}", from=1-1, to=1-2]
	\arrow["{\sigma_c[f.A.B]}"', from=6-1, to=6-2]
	\arrow[Rightarrow, no head, from=1-1, to=6-1]
	\end{tikzcd}\] where $C'\equiv C[\pair_{A}^{B}P_{C[\pair_{A}^{B}]}]\;][\;\splitt_{c}[\pair_A^B][f.A.B];c[f.A.B]\;]$, commutes i.e. $\sigma_{c[f.A.B]}=\sigma_c[f.A.B]$.
\end{itemize}
\end{defi}

We have all the notions that we need to give the following:

\begin{defi}\label{modelofsomeT}
Let $\cattypes$ be a category with attributes. \begin{itemize}
    \item If $\T$ is a given extensional type theory, we say that $\cattypes$ is \textbf{model of $\T$} if it is equipped with semantic extensional identity types, with semantic dependent product types, with semantic dependent sum types and with a choice of a semantic type (and of a  semantic term) in context $\1$ for every atomic type (and every atomic term) of $\T$.
    \item If $\T$ is a given propositional type theory, we say that $\cattypes$ is \textbf{model of $\T$} if it is equipped with semantic propositional identity types, with semantic propositional dependent product types, with semantic propositional dependent sum types and with a choice of a semantic type (and of a  semantic term) in context $\1$ for every atomic type (and every atomic term) of $\T$.
\end{itemize}
\end{defi}

\begin{defi}\label{strict morphism}
Let $\cattypesc$ and $\cattypesd$ be categories with attributes. A \textbf{morphism of categories with attributes}: $$F:\cattypesc\to\cattypesd$$ consists of:
\begin{itemize}
    \item a functor $F: \mathcal{C}\to\mathcal{D}$ preserving the terminal object;
    \item a natural transformation, that we continue calling $F$, from $\tp_{\mathcal{C}}$ to $\tp_{\mathcal{D}}F$, i.e. the equality: $$F(A[f])=(FA)[Ff]$$ holds for every semantic type $A$ in some semantic context $\Gamma$ of $\mathcal{C}$ and every morphism of semantic contexts $\Delta\xrightarrow{f}\Gamma$ of $\mathcal{C}$;
\end{itemize} in such a way that the equality: $$F(-.-)=(F-).(F-)$$ holds and that $FP_A=P_{FA}$ for every semantic type $A$ in some semantic context $\Gamma$ of $\mathcal{C}$.

Moreover, if $\cattypesc$ and $\cattypesd$ are models of a given extensional (propositional, respectively) type theory $\T$, then a \textbf{morphism of $\T$}: $$\cattypesc\to\cattypesd$$ is a morphism of categories with attributes $\cattypesc\to\cattypesd$ preserving the semantic extensional (propositional, respectively) identity types, the semantic (propositional, respectively) dependent product types, the semantic (propositional, respectively) dependent sum types and the choice of the semantic types (and of the semantic term) in context $\1$.
\end{defi}

\begin{rem} The one that we presented in \autoref{strict morphism} happens to be the \textit{strict} notion of morphism between categories with attributes, and it is considered e.g. by Cartmell \cite{phdcartmell} and by Kapulkin and Lumsdaine \cite{KAPULKIN2021106563}. In this paper we only deal with morphisms between categories with attributes in this strict form: they all strictly commute  with the given semantic context extensions and the display maps. The syntactic model of an extensional (propositional) type theory (see \autoref{syntactic model}) enjoys a strict universal initiality property with respect to morphisms of the given existential (propositional, respectively) type theory in this strong form.

However several weakenings of this notion are available in the literature, depending on their strictness in the commutativity with the category with attributes structure. E.g. the notion used by Clairambault and Dybjer \cite{clairambault_dybjer_2014} preserves the semantic context extension only up to natural isomorphism. We refer the reader to \cite{newstead2018algebraic} for more details (regarding in this case the related structure of natural model).
\end{rem}

\noindent \textbf{The propositional type theory $\sptt$.} Now, let us consider the given $\ptt$ (see \autoref{primitiveterminology}). Let us consider the propositional type theory contained in $\ptt$ (meaning that all of its contexts, types, terms and judgements are contexts, types, terms and judgements of $\ptt$) whose atomic types are the ones of $\ptt$ that are provably h-sets in $\ptt$ (and whose atomic terms are the atomic terms of these atomic h-sets in $\ptt$). We indicate this specific propositional type theory as $\sptt$, since the contexts of $\sptt$ are the ones that we call \textit{homotopy elementary} contexts of $\ptt$ (see \autoref{simple context}).

\begin{rem}\label{extensional is propositional}
Since $\ett$ and $\sptt$ have the same atomic types and terms (see \autoref{primitiveterminology}), every model of $\ett$ is canonically a model of $\sptt$: in order to obtain a choice of the terms $\j_c$, $\h_c$, $\beta_b^a$, $\eta_z$ and $\sigma_c$, one defines them as instances of $\refl_T^t$ for opportune semantic types $T$ and semantic terms $t$. Then, all of the additional compatibilities with the substitution follow in fact by the one of $\refl_T$.

With this choice of a structure of model of $\sptt$ for every model of $\ett$, every morphism of $\ett$ is canonically a morphism of $\sptt$.
\end{rem}

\begin{rem}\label{syntactic model}
The category whose objects are the contexts of $\ett$ and whose arrows are the morphisms of contexts (identified up to
renaming their free variables and up to componentwise equality judgement) has for terminal object the empty context and constitutes a category with attributes with the following data:
\begin{itemize}
    \item the semantic types in a given semantic context $\gamma:\Gamma$ are the type judgements $\judgext{\gamma:\Gamma}{A(\gamma):\type}$ of $\ett$ in context $\gamma$; the presheaf of semantic types act on the morphisms of contexts by substitution;
    \item the semantic context extension maps a pair $(\gamma:\Gamma,\judgext{\gamma}{A(\gamma):\type})$ to the context: $$\gamma:\Gamma,x:A(\gamma)$$ and an arrow $\judgext{\delta:\Delta}{f(\delta):\Gamma}$ of target $(\gamma:\Gamma,\judgext{\gamma}{A(\gamma):\type})$ to the morphism of contexts: $$\delta,x'\xrightarrow{\;\;\judgext{\delta:\Delta,x':A(f(\delta))}{f(\delta):\Gamma,x':A(f(\delta))}\;\;}\gamma,x;$$
    \item whenever $\judgext{\gamma}{A(\gamma):\type}$ is a semantic type in semantic context $\gamma$ i.e. a type judgement in context $\gamma$, then the display map $P_{\judgext{\gamma}{T(\gamma):\type}}$ is the morphism of contexts: $$\judgext{\gamma,x:A(\gamma)}{\gamma};$$ its sections---i.e. the semantic terms of $\judgext{\gamma}{A(\gamma):\type}$---are the morphisms of contexts $\gamma\to\gamma,x$ of the form: $$\judgext{\gamma}{\gamma,a(\gamma):A(\gamma)}$$ for some term judgement $\judgext{\gamma}{a(\gamma):A(\gamma)}$.
\end{itemize} We indicate as $\ettm$ this category with attributes. Then $\ettm$ is a model of $\ett$ with the clear choices of the semantic extensional identity types, of the semantic dependent sum types, and of the semantic dependent product types. The choice of a semantic type (and of a  semantic term) in empty context for every atomic type (and every atomic term) is the identity.

Analogously, $\ptt$ and $\sptt$ form models of themselves, which we indicate as $\pttm$ and $\spttm$.
\end{rem}

\begin{thm}[Soundness]\label{soundness i.e. initiality}
For every model $\cattypes$ of $\ett$, there is unique a morphism $\ettm\to\cattypes$ of $\ett$.

Analogously, for every model $\cattypes$ of $\ptt$, there is unique a morphism: $$\pttm\to\cattypes$$ of $\ptt$ and, for every model $\cattypes$ of $\sptt$, there is unique a morphism: $$\spttm\to\cattypes$$ of $\sptt$.
\end{thm}

We conclude this subsection by noticing the following \autoref{canonical interpretation}. For further details on this general notion of semantics (or equivalent ones) for dependent type theories, we refer the reader to \cite{phdcartmell,pitts2000,MR1134134}.

\begin{rem}\label{canonical interpretation}
By \autoref{extensional is propositional} and by \autoref{soundness i.e. initiality}, there is unique a morphism: $$\spttm\to\ettm$$ of $\sptt$. We call it \textbf{canonical interpretation} of $\sptt$ into $\ett$ and denote it as $|\cdot|$. A priori, being a morphism of $\sptt$, the mapping $|\cdot|$ is defined on contexts, morphisms of contexts and type judgements: we now show how to extend the mapping $|\cdot|$ to types and terms in context.

Let $\gamma:\Gamma$ and let $\judge{\gamma}{A(\gamma):\type}$ be a type judgement of $\sptt$ in context $\gamma$ (i.e. an \textit{h-elementary} type judgement of $\ptt$ in \textit{h-elementary} context $\gamma$---see \autoref{simple type} and \autoref{simple context}). As $\judge{\gamma}{A(\gamma):\type}$ is a semantic type in semantic context $\gamma$ in $\spttm$, then $|\judge{\gamma}{A(\gamma):\type}|$ needs to be a semantic type in semantic context $|\gamma|$ in $\ettm$ i.e. a type judgement of $\ett$ in context $|\gamma|$. Hence $|\judge{\gamma}{A(\gamma):\type}|$ is of the form: $$\judgext{|\gamma|}{|A(\gamma)|:\type}$$ where $|A(\gamma)|$ denotes therefore a type of $\ett$ in context $|\gamma|$.

Let $\judge{\gamma}{a(\gamma):A(\gamma)}$ be a term judgement of $\sptt$. Then the morphism $\judge{\gamma}{\gamma:\Gamma,a(\gamma):A(\gamma)}$ is a section of $P_{\judge{\gamma}{A(\gamma):\type}}$ in $\spttm$ and therefore $|\judge{\gamma}{\gamma:\Gamma,a(\gamma):A(\gamma)}|$ is a section $|\gamma|\to|\gamma|,|x|:|A(\gamma)|$ of: $$|P_{\judge{\gamma}{A(\gamma):\type}}|=P_{|\judge{\gamma}{A(\gamma):\type}|}=P_{\judgext{|\gamma|}{|A(\gamma)|:\type}}=\judgext{|\gamma|,|x|:|A(\gamma)|}{|\gamma|}$$ in $\ettm$. Hence $|\judge{\gamma}{\gamma:\Gamma,a(\gamma):A(\gamma)}|$ needs to be of the form: $$\judgext{|\gamma|}{|\gamma|,|a(\gamma)|:|A(\gamma)|}$$ for a term judgement $\judgext{|\gamma|}{|a(\gamma)|:|A(\gamma)|}$, where $|a(\gamma)|$ denotes therefore a term of $\ett$ in context $|\gamma|$ and of type $|A(\gamma)|$.

If $\judge{\delta}{f(\delta):\Gamma}$ is a morphism of contexts of $\spttm$ and if we write $|A|(|\gamma|)$ and $|a|(|\gamma|)$ for $|A(\gamma)|$ and $|a(\gamma)|$---respectively---in context $|\gamma|$, we observe that: \begin{center} $\judgext{|\gamma|}{|A(f(\delta))|\equiv|A|(|f(\delta)|)}$ \, and \, $\judgext{|\gamma|}{|a(f(\delta))|=|a|(|f(\delta)|)}$ \end{center} because $|\cdot|$ needs to commute with the substitution.
\end{rem}

\section{Homotopy equivalences of contexts in Propositional Type Theory}\label{section:homotopy equivalencesofcontexts}\label{section2.3}

This section mostly deals with the notion of morphism of contexts and the concept of generalised identity type of two parallel morphisms of context. Here we briefly recap these notions as long as we need them. For more details, we refer the reader to \cite{MR2469279}.

\medskip

Suppose that we are given two contexts $\gamma:\Gamma$ and $\delta:\Delta$ of $\ptt$, where the former is an abbreviation for the list $\gamma_1 : \Gamma_1,\gamma_2: \Gamma_2(\gamma_1),...,\gamma_n:\Gamma_n(\gamma_1,...,\gamma_{n-1})$ and the latter for $\delta_1 : \Delta_1,\delta_2: \Delta_2(\delta_1),...,\delta_m:\Delta_m(\delta_1,...,\delta_{m-1})$. With the expression $\judge{\gamma:\Gamma}{a(\gamma):\Delta}$ we mean a list of judgements: \begin{itemize}
    \item $\judge{\gamma : \Gamma}{a_1(\gamma):\Delta_1}$
    \item $\judge{\gamma : \Gamma}{a_2(\gamma):\Delta_2(a_1(\gamma))}$
    \item $\judge{\gamma : \Gamma}{a_3(\gamma):\Delta_3(a_1(\gamma),a_2(\gamma))}$
    \item ...
    \item $\judge{\gamma : \Gamma}{a_m(\gamma):\Delta_m(a_1(\gamma),a_2(\gamma),a_3(\gamma),...,a_{m-1}(\gamma))}$
\end{itemize} and we call such a list a \textit{morphism of contexts} $\gamma\to\delta$. Now, suppose that $\delta' : \Delta'$ indicates the context: $$\delta_1:\Delta_1,...,\delta_{m-1}:\Delta_{m-1}(\delta_1,...,\delta_{m-2})$$ and $\judge{\gamma:\Gamma}{a'(\gamma):\Delta'}$ indicates the list: \begin{itemize}
    \item $\judge{\gamma : \Gamma}{a_1(\gamma):\Delta_1}$
    \item $\judge{\gamma : \Gamma}{a_2(\gamma):\Delta_2(a_1(\gamma))}$
    \item $\judge{\gamma : \Gamma}{a_3(\gamma):\Delta_3(a_1(\gamma),a_2(\gamma))}$
    \item ...
    \item $\judge{\gamma : \Gamma}{a_{m-1}(\gamma):\Delta_m(a_1(\gamma),a_2(\gamma),a_3(\gamma),...,a_{m-2}(\gamma))}$
\end{itemize} and $A(\delta')$ indicates the type in context $\Delta_m(\delta_1,...,\delta_{m-1})$. Then we may also write anyone of the following: $$\judge{\gamma:\Gamma}{a(\gamma):\Delta',A}\text{ \, \, \, }\judge{\gamma:\Gamma}{a'(\gamma),a_m(\gamma):\Delta}$$ $$\judge{\gamma:\Gamma}{a'(\gamma),a_m(\gamma):\Delta',A}\text{ \, \, \, }\judge{\gamma:\Gamma}{a'(\gamma):\Delta',a_m(\gamma):A(a'(\gamma))}$$ in order to indicate the same morphism of contexts $\judge{\gamma:\Gamma}{a(\gamma):\Delta}$.

If we are given two parallel morphisms of contexts $\judge{\gamma}{a(\gamma):\Delta}$ and $\judge{\gamma}{b(\gamma):\Delta}$, the expression $\judge{\gamma}{p(\gamma):a(\gamma)=b(\gamma)}$ indicates the list: \begin{itemize}
    \item $\judge{\gamma}{p_1(\gamma):a_1(\gamma)=b_1(\gamma)}$
    \item $\judge{\gamma}{p_2(\gamma):a_2(\gamma)=p_1(\gamma)^*b_2(\gamma)}$
    \item $\judge{\gamma}{p_3(\gamma):a_3(\gamma)=(p_1(\gamma),p_2(\gamma))^*b_3(\gamma)}$
    \item ...
    \item $\judge{\gamma}{p_m(\gamma):a_m(\gamma)=(p_1(\gamma),...,p_{m-1}(\gamma))^*b_m(\gamma)}$
\end{itemize} where the operations $(p_1(\gamma),...,p_{k}(\gamma))^*$ are defined by sequential (generalised) path inductions on $p_{k}$, $p_{k-1}$, ..., $p_2$, and $p_1$ and hence make the identity types: $$(\;\refl(a_1(\gamma)),...,\refl(a_{k}(\gamma))\;)^*d=d$$ inhabited (remind that $\ptt$ has \textit{propositional} identity types) if $d:\Delta_{k+1}(b_1(\gamma),...,b_k(\gamma))$. We call \textit{context propositional equality} this new meaning of the symbol $=$ in between two parallel morphisms of contexts.

As shown by Gambino and Garner \cite{MR2469279,MR2525957}, the expression $\judge{\gamma:\Gamma}{p(\gamma):a(\gamma)=b(\gamma)}$ formally verifies the same rules of \autoref{propositional id} verified by the propositional equality and one can simply prove this by sequential (generalised) path induction. In particular, a context elimination (i.e. path induction) rule, with a corresponding context propositional computation rule, is satisfied by the context propositional equality, and moreover every context homotopy equivalence---see \autoref{context homotopy equivalence}---is also a context half-adjoint equivalence.

\begin{defi} \label{context homotopy equivalence}
A \textbf{context homotopy equivalence} between $\gamma:\Gamma$ and $\delta:\Delta$ is a couple of morphisms of contexts of the form:
\begin{itemize}
    \item $\judge{\gamma : \Gamma}{\f(\gamma) : \Delta}$
    \item $\judge{\delta : \Delta}{\g(\delta): \Gamma}$
\end{itemize} such that \textit{there exist} couples of judgements of the form:
\begin{itemize}
    \item $\judge{\gamma : \Gamma}{\p(\gamma):\g(\f(\gamma))=\gamma}$
    \item $\judge{\delta : \Delta}{\q(\delta):\f(\g(\delta))=\delta}$
\end{itemize} hence, as usual, the expression $\judge{\gamma : \Gamma}{\f(\gamma) : \Delta}$ is an abbreviation for the list of judgements:
\begin{itemize}
    \item $\judge{\gamma : \Gamma}{\f_1(\gamma):\Delta_1}$
    \item $\judge{\gamma : \Gamma}{\f_2(\gamma):\Delta_2(\f_1(\gamma))}$
    \item $\judge{\gamma : \Gamma}{\f_3(\gamma):\Delta_3(\f_1(\gamma),\f_2(\gamma))}$
    \item ...
    \item $\judge{\gamma : \Gamma}{\f_m(\gamma):\Delta_m(\f_1(\gamma),\f_2(\gamma),\f_3(\gamma),...,\f_{m-1}(\gamma))}$
\end{itemize}
and the expression $\judge{\gamma : \Gamma}{\p(\gamma):\g(\f(\gamma))=\gamma}$ is an abbreviation for the list of judgements: \begin{itemize}
    \item $\judge{\gamma :\Gamma}{\p_1(\gamma):\g_1(\f(\gamma))=\gamma_1}$
    \item $\judge{\gamma :\Gamma}{\p_2(\gamma):\g_2(\f(\gamma))=\p_1(\gamma)^*\gamma_2}$
    \item $\judge{\gamma :\Gamma}{\p_3(\gamma):\g_3(\f(\gamma))=(\p_1(\gamma),\p_2(\gamma))^*\gamma_3}$
    \item ...
    \item $\judge{\gamma :\Gamma}{\p_n(\gamma):\g_n(\f(\gamma))=(\p_1(\gamma),\p_2(\gamma),...,\p_{n-1}(\gamma))^*\gamma_n}$.
\end{itemize} \end{defi}

In \autoref{section2.3.I} we show under what hypotheses these equivalences can be extended to wider contexts.

\subsection{Extension of context homotopy equivalences}\label{section2.3.I} In what follows, and throughout the entire paper, we will often adopt the following conventions for naming variables: if we denote two given types by $A$ and $A'$, and denote a variable of type $A$ by $x$, then we will often write $\underline{x}$ to denote a variable of type $A'$. Additionally, we will use $x'$, $x''$, etc. (and similarly $\underline{x}'$, $\underline{x}''$, etc.) to denote variables typed by re-indexings of $A$ (and of $A'$, respectively).

Let $\gamma:\Gamma$ and $\delta:\Delta$ be contexts and let us assume that we are given a context homotopy equivalence $(\f;\g)$ as follows: \begin{itemize}
    \item $\judge{\gamma : \Gamma}{\f(\gamma) : \Delta}$
    \item $\judge{\delta : \Delta}{\g(\delta): \Gamma}$
\end{itemize} between them as before. Suppose that we are given the judgements:
\begin{itemize}
    \item $\judge{\gamma:\Gamma}{A(\gamma):\type}$
    \item $\judge{\delta:\Delta}{A'(\delta):\type}$
    \item $\judge{\gamma:\Gamma,x:A(\gamma)}{f(\gamma,x):A'(\f(\gamma))}$
    \item $\judge{\gamma:\Gamma,\underline{x}':A'(\f(\gamma))}{g(\gamma,\underline{x}'):A(\gamma)}$
\end{itemize} such that there are terms of the form:
\begin{itemize}
    \item $\judge{\gamma:\Gamma,x:A(\gamma)}{p(\gamma,x): x=g(\gamma,f(\gamma,x))}$
    \item $\judge{\gamma:\Gamma,\underline{x}':A'(\f(\gamma))}{q(\gamma,\underline{x}'): \underline{x}'=f(\gamma,g(\gamma,\underline{x}'))}$.
\end{itemize} In other words, we are given a \textit{homotopy equivalence $(f;g)$ between the types in context $A(\gamma)$ and $A'(\delta)$ relative to the context homotopy equivalence $(\f;\g)$}. We can use these data in order to augment the context homotopy equivalence $(\f;\g)$ to a context homotopy equivalence between the contexts: \begin{center} $\gamma:\Gamma,x : A(\gamma)$ and $\delta : \Delta,\underline{x}:A'(\delta)$. \end{center} The remainder of the current subsection is devoted to showing this construction: \textit{the reader who is willing to skip the details may turn to \autoref{context homotopy equivalence expansion}.}

\medskip

We observe that $\judge{\delta:\Delta,\underline{x}'':A'(\f(\g(\delta)))}{g(\g(\delta),\underline{x}''):A(\g(\delta))}$ and that: $$\judge{\delta:\Delta,\underline{x}:A'(\delta)}{\q(\delta)^*\underline{x}:A'(\f(\g(\delta)))}$$ hence $\judge{\delta:\Delta,\underline{x}:A'(\delta)}{g(\g(\delta),\q(\delta)^*\underline{x}):A(g(\delta))}$. Let us rename:
\begin{itemize}
    \item $\judge{\gamma:\Gamma,x:A(\gamma)}{\f_{m+1}(\gamma,x)\equiv f(\gamma,x):A'(\f(\gamma))}$
    \item $\judge{\delta:\Delta,\underline{x}:A'(\delta)}{\g_{n+1}(\delta,\underline{x})\equiv g(\g(\delta),\q(\delta)^*\underline{x}):A(\g(\delta))}$.
\end{itemize}

We observe that: $$\f_{m+1}(\g(\delta),\g_{n+1}(\delta,\underline{x}))\equiv f(\g(\delta),g(\g(\delta),\q(\delta)^*\underline{x}))\overset{q(\g(\delta),\q(\delta)^*\underline{x})}{=}\q(\delta)^*\underline{x}$$ in context $\delta:\Delta,\underline{x}:A'(\delta)$. Moreover: $$\begin{aligned}\g_{n+1}(\f(\gamma),\f_{m+1}(\gamma,x))&\equiv g(\g(\f(\gamma)),\q(\f(\gamma))^*\f_{m+1}(\gamma,x))\\&=g(\g(\f(\gamma)),\f(\p(\gamma))^*\f_{m+1}(\gamma,x))\\&=\p(\gamma)^*x\end{aligned}$$ in context $\gamma:\Gamma,x:A(\gamma)$, where the first identity type is inhabited as one can assume w.l.o.g. that: $$\judge{\gamma :\Gamma}{\q(\f(\gamma))=\f(\p(\gamma))}$$ and the second by applying based (generalised) path induction $n$ times on $\p_1(\gamma)$, $\p_2(\gamma)$, ..., $\p_n(\gamma)$ and finally using that $\judge{\gamma:\Gamma,x:A(\gamma)}{p(\gamma,x): g(\gamma,f(\gamma,x))=x}$.

In conclusion, we saw that we can augment the contexts $\gamma:\Gamma$ and $\delta:\Delta$ in: $$\gamma:\Gamma,x:A(\gamma)\textnormal{ and }\delta:\Delta,\underline{x}:A'(\delta)$$ respectively, in such a way that they continue being homotopy equivalent. In fact the extended context morphisms:
\begin{itemize}
    \item $\judge{\gamma : \Gamma,  x:A(\gamma)}{\f(\gamma) : \Delta,\f_{m+1}(\gamma,x):A'(\f(\gamma))}$ \, i.e.
    
    $\judge{\gamma : \Gamma,  x:A(\gamma)}{(\f,\f_{m+1})(\gamma,x):\Delta,A'}$
    \item $\judge{\delta : \Delta, \underline{x}:A'(\delta)}{\g(\delta): \Gamma,\g_{n+1}(\delta,\underline{x}):A(\g(\delta))}$ \,  i.e.
    
    $\judge{\delta : \Delta, \underline{x}:A'(\delta)}{(\g,\g_{n+1})(\delta,\underline{x}):\Gamma,A}$
\end{itemize}
constitute a homotopy equivalence, since: \begin{center}
    $\judge{\gamma : \Gamma}{\p(\gamma):\g(\f(\gamma))=\gamma}$ \, and \, $\judge{\gamma:\Gamma,x:A(\gamma)}{\g_{n+1}(\f(\gamma),\f_{m+1}(\gamma,x))=\p(\gamma)^*x}$
\end{center} that is: \[\begin{aligned}  \judge{\gamma : \Gamma,x:A(\gamma)}{\g(\f(\gamma)),\g_{n+1}(\f(\gamma),\f_{m+1}(\gamma,x))&=\gamma,x} \textnormal{ \, i.e.} \\
\judge{\gamma : \Gamma,x:A(\gamma)}{(\g,\g_{n+1})((\f,\f_{m+1})(\gamma,x))&=\gamma,x} \end{aligned}\] and since:
\begin{center}
    $\judge{\delta : \Delta}{\q(\delta):\f(\g(\delta))=\delta}$ \, and \, $\judge{\delta :\Delta,\underline{x}: A'(\delta)}{\f_{m+1}(\g(\delta),\g_{n+1}(\delta,\underline{x}))=\q(\delta)^*\underline{x}}$
\end{center} that is: \[\begin{aligned}  \judge{\delta :\Delta,\underline{x}: A'(\delta)}{\f(\g(\delta)),\f_{m+1}(\g(\delta),\g_{n+1}(\delta,\underline{x}))&=\delta,\underline{x}} \textnormal{ \, i.e.} \\
\judge{\delta :\Delta,\underline{x}: A'(\delta)}{(\f,\f_{m+1})((\g,\g_{n+1})(\delta,\underline{x}))&=\delta,\underline{x}}. \end{aligned}\]

Let us summarise this into the following:

\begin{lem}[Extension]\label{context homotopy equivalence expansion}
Let $\gamma : \Gamma$ and $\delta : \Delta$ and let: \begin{center}$\judge{\gamma:\Gamma}{\f(\gamma):\Delta}$\\ $\judge{\delta:\Delta}{\g(\delta):\Gamma}$\end{center} be a context homotopy equivalence $(\f;\g)$. If we are given types: $$\judge{\gamma:\Gamma}{A(\gamma):\type}\text{ and }\judge{\delta:\Delta}{A'(\delta):\type}$$ together with a homotopy equivalence: \begin{center}$\judge{\gamma:\Gamma,x:A(\gamma)}{f(\gamma,x):A'(\f(\gamma))}$\\ $\judge{\gamma:\Gamma,\underline{x}':A'(\f(\gamma))}{g(\gamma,\underline{x}'):A(\gamma)}$\end{center} between $A(\gamma)$ and $A'(\delta)$ relative to $(\f;\g)$ then: \begin{center} $\judge{\gamma : \Gamma,  x:A(\gamma)}{(\f,\f_{m+1})(\gamma,x):\Delta,A'}$\\ $\judge{\delta : \Delta, \underline{x}:A'(\delta)}{(\g,\g_{n+1})(\delta,\underline{x}):\Gamma,A}$\end{center} is a context homotopy equivalence $(\f,\f_{m+1};\g,\g_{n+1})$, where: \begin{center} $\judge{\gamma:\Gamma,x:A(\gamma)}{\f_{m+1}(\gamma,x)\equiv f(\gamma,x):A'(\f(\gamma))}$ \\ $\judge{\delta:\Delta,\underline{x}:A'(\delta)}{\g_{n+1}(\delta,\underline{x})\equiv g(\g(\delta),\q(\delta)^*\underline{x}):A(\g(\delta))}$.
\end{center} We call $(\f,\f_{m+1};\g,\g_{n+1})$ the \textnormal{extension of $(\f;\g)$ via $(f;g)$}.
\end{lem}

In the next subsections we analyse specific shapes of homotopy equivalences between types relative to a given context homotopy equivalence. The first regards the ones coming from an application of the dependent product constructor.

\subsection{Dependent product of homotopy equivalences}\label{section2.3.II}

We start the current subsection by briefly describing the following data:

\medskip

\noindent \textbf{Data.} Let us assume that $\gamma:\Gamma$ and $\delta:\Delta$ are contexts and that $(\f;\g)$ is a context homotopy equivalence between them. Moreover, let us assume that we are given judgements: \begin{itemize}
        \item $\judge{\gamma:\Gamma}{A(\gamma):\type}$ \, and \, $\judge{\gamma:\Gamma,x:A(\gamma)}{B(\gamma,x):\type}$
        \item $\judge{\delta:\Delta}{A'(\delta):\type}$ \, and \, $\judge{\delta:\Delta,\underline{x}:A'(\delta)}{B'(\delta,\underline{x}):\type}$
\end{itemize} together with a homotopy equivalence $(f_1;g_1)$: \begin{center}$\judge{\gamma:\Gamma,x:A(\gamma)}{f_1(\gamma,x):A'(\f(\gamma))}$\\ $\judge{\gamma:\Gamma,\underline{x}':A'(\f(\gamma))}{g_1(\gamma,\underline{x}'):A(\gamma)}$\end{center} between $A(\gamma)$ and $A'(\delta)$ relative to $(\f;\g)$, and a homotopy equivalence $(f_2;g_2)$: \begin{center}$\judge{\gamma:\Gamma,x:A(\gamma),y:B(\gamma,x)}{f_2(\gamma,x,y):B'(\f(\gamma),f_1(\gamma,x))}$ \\ $\judge{\gamma:\Gamma,x:A(\gamma),\underline{y}':B'(\f(\gamma),f_1(\gamma,x))}{g_2(\gamma,x,\underline{y}'):B(\gamma,x)}$\end{center} between $B(\gamma,x)$ and $B'(\delta,\underline{x})$ relative to the extension $(\f,\f_{m+1};\g,\g_{n+1})$ of $(\f;\g)$ via $(f_1;g_1)$ (see \autoref{context homotopy equivalence expansion}).

\medskip

\noindent \textbf{Fact.} We can use these data to construct a homotopy equivalence $(f^\Pi;g^\Pi)$ between the types $\Pi_{x : A(\gamma)}B(\gamma,x)$ and $\Pi_{\underline{x}:A'(\delta)}B'(\delta,\underline{x})$ relative to $(\f;\g)$ as follows, leading to \autoref{pi equiv} below. The reader may go through the proof under the additional assumption that the contexts $\gamma$ and $\delta$ are empty. However, here we present the proof in full generality because there are several (homotopic) ways to define the pair $(f^\Pi;g^\Pi)$ and the one we choose here will determine the style of the subsequent proofs. A similar remark applies to \autoref{section2.3.III} and \autoref{section2.3.IV}.

\medskip

Let us start by considering homotopies: \begin{center}
    $\judge{\gamma:\Gamma,x:A(\gamma)}{p_1(\gamma,x):x=g_1(\gamma,f_1(\gamma,x))}$\\
    $\judge{\gamma:\Gamma,\underline{x}':A'(\f(\gamma))}{q_1(\gamma,\underline{x}'):\underline{x}'=f_1(\gamma,g_1(\gamma,\underline{x}'))}$
\end{center} and homotopies: \begin{center}
    $\judge{\gamma:\Gamma,x:A(\gamma),y:B(\gamma,x)}{p_2(\gamma,x,y):y=g_2(\gamma,x,f_2(\gamma,x,y))}$\\
    $\judge{\gamma:\Gamma,x:A(\gamma),\underline{y}':B'(\f(\gamma),f_1(\gamma,x))}{q_2(\gamma,x,\underline{y}'):\underline{y}'=f_2(\gamma,x,g_2(\gamma,x,\underline{y}'))}$
\end{center} as in our assumptions.

\medskip

Let us fix the context $\gamma :\Gamma,z :\Pi_{x : A(\gamma)}B(\gamma,x)$ and let us observe that $\judge{\gamma,z,x:A}{\ev(z,x):B(\gamma,x)}$, hence: $$\judge{\gamma,z,x}{f_2(\gamma,x,\ev(z,x)):B'(\f(\gamma),f_1(\gamma,x))}.$$ Since $\judge{\gamma,\underline{x}':A'(\f(\gamma))}{g_1(\gamma,\underline{x}'):A(\gamma)}$, then: $$\judge{\gamma,z,\underline{x}'}{f_2(\gamma,g_1(\gamma,\underline{x}'),\ev(z,g_1(\gamma,\underline{x}'))):B'(\f(\gamma),f_1(\gamma,g_1(\gamma,\underline{x}')))}$$ hence: $$\judge{\gamma,z,\underline{x}'}{q_1(\gamma,\underline{x}')^*f_2(\gamma,g_1(\gamma,\underline{x}'),\ev(z,g_1(\gamma,\underline{x}'))):B'(\f(\gamma),\underline{x}')}.$$ We conclude that:  $$\judge{\gamma,z}{f^\Pi(\gamma,z)\equiv\lambda \underline{x}':A'(\f(\gamma))\;.\;q_1(\gamma,\underline{x}')^*f_2(\gamma,g_1(\gamma,\underline{x}'),\ev(z,g_1(\gamma,\underline{x}')))}$$ is a term of type: $$\Pi_{\underline{x}':A'(\f(\gamma))}B(\f(\gamma),\underline{x}')\equiv [\Pi_{\underline{x}:A'(\delta)}B'(\delta,\underline{x})](\f(\gamma)).$$

\medskip

Vice versa, let us fix the context $\gamma :\Gamma, \underline{z}':\Pi_{\underline{x}':A'(\f(\gamma))}B(\f(\gamma),\underline{x}')$ and let us observe that: $$\judge{\gamma,\underline{z}',\underline{x}':A'(\f(\gamma))}{\ev(\underline{z}',\underline{x}'):B'(\f(\gamma),\underline{x}')}$$ hence $\judge{\gamma,\underline{z}',x:A(\gamma)}{\ev(\underline{z}',f_1(\gamma,x)):B'(\f(\gamma),f_1(\gamma,x))}$ and therefore: $$\judge{\gamma,\underline{z}',x:A(\gamma)}{g_2(\gamma,x,\ev(\underline{z}',f_1(\gamma,x))):B(\gamma,x).}$$ We conclude that: $$\judge{\gamma,\underline{z}'}{g^\Pi(\gamma,\underline{z}')\equiv\lambda x:A(\gamma)\;.\; g_2(\gamma,x,\ev(\underline{z}',f_1(\gamma,x))):\Pi_{x:A(\gamma)}B(\gamma,x).}$$

\medskip

We claim that $(f^\Pi;g^\Pi)$ is a homotopy equivalence between:\begin{center}$\Pi_{x:A(\gamma)}B(\gamma,x)$ and $\Pi_{\underline{x}:A'(\delta)}B'(\delta,\underline{x})$\end{center}relative to $(\f;\g)$. We start by verifying that $g^\Pi(\gamma,f^\Pi(\gamma,z))\equiv z$ in context $\gamma,z$. Let $\underline{z}'\equiv f^\Pi(\gamma,z)$ and let us observe that: $$\begin{aligned}\ev(\underline{z}',f_1(\gamma,x))&=q_1(\gamma,f_1(\gamma,x))^*f_2(\gamma,g_1(\gamma,f_1(\gamma,x)),\ev(z,g_1(\gamma,f_1(\gamma,x))))\\&=f_1(\gamma,p_1(\gamma,x))^*f_2(\gamma,g_1(\gamma,f_1(\gamma,x)),\ev(z,g_1(\gamma,f_1(\gamma,x))))\\&=f_2(\gamma,x,\ev(z,x))\end{aligned}$$ where the first equality follows by $\beta$-reduction, the second because without loss of generality $q_1(\gamma,f_1(\gamma,x))=f_1(\gamma,p_1(\gamma,x))$ and the third by based (generalised) path induction on $p_1(\gamma,x)$. By propositional functoriality: $$g_2(\gamma,x,\ev(\underline{z}',f_1(\gamma,x)))=g_2(\gamma,x,f_2(\gamma,x,\ev(z,x)))\overset{p_2(\gamma,x,\ev(z,x))}{=}\ev(z,x)$$ hence by propositional function extensionality: $$g^\Pi(\gamma,f^\Pi(\gamma,z))\equiv g^\Pi(\gamma,\underline{z}')=\lambda x : A(\gamma)\;.\;\ev(z,x)=z$$ where the last identity type is inhabited by $\eta$-expansion.

Vice versa, let us verify that $f^\Pi(\gamma,z)=\underline{z}'$, where $z\equiv g^\Pi(\gamma,\underline{z}')$. At first, we observe that: $$\ev(z,g_1(\gamma,\underline{x}'))=g_2(\gamma,g_1(\gamma,\underline{x}'),\ev(\underline{z}',f_1(\gamma,g_1(\gamma,\underline{x}'))))$$ hence: $$\begin{aligned}f_2(\gamma,g_1(\gamma,\underline{x}'),\ev(z,g_1(\gamma,\underline{x}')))&=f_2(\gamma,g_1(\gamma,\underline{x}'),g_2(\gamma,g_1(\gamma,\underline{x}'),\ev(\underline{z}',f_1(\gamma,g_1(\gamma,\underline{x}')))))\\&=\ev(\underline{z}',f_1(\gamma,g_1(\gamma,\underline{x}')))\end{aligned}$$ where the first identity type is inhabited by propositional functoriality and the second one by the term: $$q_2(\gamma,g_1(\gamma,\underline{x}'),\ev(\underline{z}',f_1(\gamma,g_1(\gamma,\underline{x}')))).$$ Secondly, we observe that: $$q_1(\gamma,\underline{x}')^*f_2(\gamma,g_1(\gamma,\underline{x}'),\ev(z,g_1(\gamma,\underline{x}')))=q_1(\gamma,\underline{x}')^*\ev(\underline{z}',f_1(\gamma,g_1(\gamma,\underline{x}')))=\ev(\underline{z}',\underline{x}')$$ where the former equality follows by propositional functoriality and the latter by based (generalised) path induction on $q_1(\gamma,\underline{x}')$. Therefore: $$f^\Pi(\gamma,g^\Pi(\gamma,\underline{z}'))\equiv f^\Pi(\gamma,z)=\lambda \underline{x}':A'(\f(\gamma))\;.\; \ev(\underline{z}',\underline{x}')=\underline{z}'$$ where the first identity type is inhabited by propositional function extensionality and the second by $\eta$-expansion.

\medskip

Let us summarise this into the following:

\begin{lem}\label{pi equiv}
Let us assume that $\gamma:\Gamma$ and $\delta:\Delta$ are contexts and that $(\f;\g)$ is a context homotopy equivalence between them. Moreover, let us assume that we are given judgements: \begin{itemize}
        \item $\judge{\gamma:\Gamma}{A(\gamma):\type}$ \, and \, $\judge{\gamma:\Gamma,x:A(\gamma)}{B(\gamma,x):\type}$
        \item $\judge{\delta:\Delta}{A'(\delta):\type}$ \, and \, $\judge{\delta:\Delta,\underline{x}:A'(\delta)}{B'(\delta,\underline{x}):\type}$
\end{itemize} together with a homotopy equivalence $(f_1;g_1)$: \begin{center}$\judge{\gamma:\Gamma,x:A(\gamma)}{f_1(\gamma,x):A'(\f(\gamma))}$\\ $\judge{\gamma:\Gamma,\underline{x}':A'(\f(\gamma))}{g_1(\gamma,\underline{x}'):A(\gamma)}$\end{center} between $A(\gamma)$ and $A'(\delta)$ relative to $(\f;\g)$, and a homotopy equivalence $(f_2;g_2)$: \begin{center}$\judge{\gamma:\Gamma,x:A(\gamma),y:B(\gamma,x)}{f_2(\gamma,x,y):B'(\f(\gamma),f_1(\gamma,x))}$ \\ $\judge{\gamma:\Gamma,x:A(\gamma),\underline{y}':B'(\f(\gamma),f_1(\gamma,x))}{g_2(\gamma,x,\underline{y}'):B(\gamma,x)}$\end{center} between $B(\gamma,x)$ and $B'(\delta,\underline{x})$ relative to the extension $(\f,\f_{m+1};\g,\g_{n+1})$ of $(\f;\g)$ via $(f_1;g_1)$ (see \autoref{context homotopy equivalence expansion}).

In context: \begin{center} $\gamma:\Gamma$ and $z : \Pi_{x:A(\gamma)}B(\gamma,x)$ and $\underline{z}':[\Pi_{\underline{x}:A'(\delta)}B'(\delta,\underline{x})](\f(\gamma))$ \end{center} if we name: \begin{center}
    $\judge{\gamma,z}{f^\Pi(\gamma,z)\equiv\lambda \underline{x}':A'(\f(\gamma))\;.\;q_1(\gamma,\underline{x}')^*f_2(\gamma,g_1(\gamma,\underline{x}'),\ev(z,g_1(\gamma,\underline{x}'))):[\Pi_{\underline{x}:A'(\delta)}B'(\delta,\underline{x})](\f(\gamma))}$
    
    \medskip
    
    $\judge{\gamma,\underline{z}'}{g^\Pi(\gamma,\underline{z}')\equiv\lambda x:A(\gamma)\;.\; g_2(\gamma,x,\ev(\underline{z}',f_1(\gamma,x))):\Pi_{x:A(\gamma)}B(\gamma,x)}$
\end{center} then $(f^\Pi;g^\Pi)$ is a homotopy equivalence between $\Pi_{x:A(\gamma)}B(\gamma,x)$ and $\Pi_{\underline{x}:A'(\delta)}B'(\delta,\underline{x})$ relative to $(\f;\g)$.
\end{lem}

\subsection{Dependent sum of homotopy equivalences}\label{section2.3.III}

As for \autoref{section2.3.II}, let us describe some:

\medskip

\noindent \textbf{Data.} Let us assume that we are given a context homotopy equivalence $(\f;\g)$ between $\gamma:\Gamma$ and $\delta:\Delta$ and judgements: \begin{itemize}
        \item $\judge{\gamma:\Gamma}{A(\gamma):\type}$ \, and \, $\judge{\gamma:\Gamma,x:A(\gamma)}{B(\gamma,x):\type}$
        \item $\judge{\delta:\Delta}{A'(\delta):\type}$ \, and \, $\judge{\delta:\Delta,\underline{x}:A'(\delta)}{B'(\delta,\underline{x}):\type}$
\end{itemize} together with a homotopy equivalence $(f_1;g_1)$: \begin{center}$\judge{\gamma:\Gamma,x:A(\gamma)}{f_1(\gamma,x):A'(\f(\gamma))}$\\ $\judge{\gamma:\Gamma,\underline{x}':A'(\f(\gamma))}{g_1(\gamma,\underline{x}'):A(\gamma)}$\end{center} between $A(\gamma)$ and $A'(\delta)$ relative to $(\f;\g)$, and a homotopy equivalence $(f_2;g_2)$: \begin{center}$\judge{\gamma:\Gamma,x:A(\gamma),y:B(\gamma,x)}{f_2(\gamma,x,y):B'(\f(\gamma),f_1(\gamma,x))}$ \\ $\judge{\gamma:\Gamma,x:A(\gamma),\underline{y}':B'(\f(\gamma),f_1(\gamma,x))}{g_2(\gamma,x,\underline{y}'):B(\gamma,x)}$\end{center} between the types $B(\gamma,x)$ and $B'(\delta,\underline{x})$ relative to the extension $(\f,\f_{m+1};\g,\g_{n+1})$ of $(\f;\g)$ via $(f_1;g_1)$ (\autoref{context homotopy equivalence expansion}).

\medskip

\noindent \textbf{Fact.} Again, we can use these data to construct a homotopy equivalence: $$(f^\Sigma;g^\Sigma)$$ between $\Sigma_{x : A(\gamma)}B(\gamma,x)$ and $\Sigma_{\underline{x}:A'(\delta)}B'(\delta,\underline{x})$ relative to $(\f;\g)$ as follows, leading to \autoref{sigma equiv} below.

\medskip

As in \autoref{section2.3.II}, let us consider homotopies: \begin{center}
    $\judge{\gamma:\Gamma,x:A(\gamma)}{p_1(\gamma,x):g_1(\gamma,f_1(\gamma,x))=x}$\\
    $\judge{\gamma:\Gamma,\underline{x}':A'(\f(\gamma))}{q_1(\gamma,\underline{x}'):f_1(\gamma,g_1(\gamma,\underline{x}'))=\underline{x}'}$
\end{center} and homotopies: \begin{center}
    $\judge{\gamma:\Gamma,x:A(\gamma),y:B(\gamma,x)}{p_2(\gamma,x,y):g_2(\gamma,x,f_2(\gamma,x,y))=y}$\\
    $\judge{\gamma:\Gamma,x:A(\gamma),\underline{y}':B'(\f(\gamma),f_1(\gamma,x))}{q_2(\gamma,x,\underline{y}'):f_2(\gamma,x,g_2(\gamma,x,\underline{y}'))=\underline{y}'}$
\end{center} as in our assumptions.

\medskip

Let us fix the context $\gamma : \Gamma, u : \Sigma_{x : A(\gamma)}B(\gamma,x)$ and let us observe that $\pi_1u : A(\gamma)$ and that $\pi_2u : B(\gamma,\pi_1u)$. Therefore $f_1(\gamma,\pi_1u):A'(\f(\gamma))$ and $f_2(\gamma, \pi_1u,\pi_2u):B'(\f(\gamma),f_1(\gamma,\pi_1u))$ hence: $$\judge{\gamma,u}{f^\Sigma(\gamma,u)\equiv\langle f_1(\gamma,\pi_1u),f_2(\gamma, \pi_1u,\pi_2u)\rangle : [\Sigma_{\underline{x}:A'(\delta)}B'(\delta,\underline{x})](\f(\gamma))}.$$

\medskip

Vice versa, if $\gamma:\Gamma,\underline{u}':\Sigma_{\underline{x}':A'(\f(\gamma))}B'(\f(\gamma),\underline{x}')$ then $\pi_1\underline{u}' : A'(\f(\gamma))$ hence $g_1(\gamma,\pi_1\underline{u}'):A(\gamma)$. Moreover $\pi_2\underline{u}' : B'(\f(\gamma),\pi_1\underline{u}')$ hence $q_1(\gamma,\pi_1\underline{u}')^*\pi_2\underline{u}':B'(\f(\gamma),f_1(\gamma,g_1(\gamma,\pi_1\underline{u}')))$. Therefore: $$g_2(\gamma,g_1(\gamma,\pi_1\underline{u}'),q_1(\gamma,\pi_1\underline{u}')^*\pi_2\underline{u}') : B(\gamma,g_1(\gamma,\pi_1\underline{u}'))$$ hence: $$\judge{\gamma,\underline{u}'}{g^\Sigma(\gamma,\underline{u}')\equiv\langle g_1(\gamma,\pi_1\underline{u}'),g_2(\gamma,g_1(\gamma,\pi_1\underline{u}'),q_1(\gamma,\pi_1\underline{u}')^*\pi_2\underline{u}')\rangle: \Sigma_{x:A(\gamma)}B(\gamma,x)}.$$

We claim that $(f^\Sigma;g^\Sigma)$ is a homotopy equivalence between:\begin{center}$\Sigma_{x : A(\gamma)}B(\gamma,x)$ and $\Sigma_{\underline{x}:A'(\delta)}B'(\delta,\underline{x})$\end{center}relative to $(\f;\g)$. In order to verify that $g^\Sigma(\gamma,\underline{u}')=u$, where $\underline{u}'\equiv f^\Sigma(\gamma,u)$, it is enough to verify that: \begin{center}
     $p:g_1(\gamma,\pi_1\underline{u}')=\pi_1u$
     
     \smallskip
     
     $g_2(\gamma,g_1(\gamma,\pi_1\underline{u}'),q_1(\gamma,\pi_1\underline{u}')^*\pi_2\underline{u}')=p^*\pi_2u$
\end{center} as this implies that $g^\Sigma(\gamma,\underline{u}')=\langle \pi_1u,\pi_2u\rangle=u$. Since: $$g_1(\gamma, \beta_1(f_1(\gamma,\pi_1u),f_2(\gamma, \pi_1u,\pi_2u)))\bullet p_1(\gamma,\pi_1u):g_1(\gamma,\pi_1\underline{u}')=\pi_1u$$ we are left to verify that: $$g_2(\gamma,g_1(\gamma,\pi_1\underline{u}'),q_1(\gamma,\pi_1\underline{u}')^*\pi_2\underline{u}')=p^*\pi_2u$$ for $p\equiv g_1(\gamma, \beta_1(f_1(\gamma,\pi_1u),f_2(\gamma, \pi_1u,\pi_2u)))\bullet p_1(\gamma,\pi_1u)$. Let us observe that the following square:\[\begin{tikzcd}
f_1(\gamma,g_1(\gamma,\pi_1\underline{u}')) \arrow[dd, "q_1(\gamma\text{,}\pi_1\underline{u}')"', Rightarrow] \arrow[rrrrrrrr, "f_1(\gamma\text{,}g_1(\gamma\text{,}\beta_1(f_1(\gamma\text{,}\pi_1u)\text{,}f_2(\gamma\text{,}\pi_1u\text{,}\pi_2u))))", Rightarrow] &&&&&&&& f_1(\gamma,g_1(\gamma,f_1(\gamma,\pi_1u))) \arrow[dd, "q_1(\gamma\text{,}f_1(\gamma\text{,}\pi_1u))", Rightarrow] \\
\\
\pi_1\underline{u}' \arrow[rrrrrrrr, "\beta_1(f_1(\gamma\text{,}\pi_1u)\text{,}f_2(\gamma\text{,}\pi_1u\text{,}\pi_2u))", Rightarrow] &&&&&&&& f_1(\gamma,\pi_1u)
\end{tikzcd}\] commutes propositionally, because $\judge{\gamma:\Gamma}{q_1(\gamma,-):f_1(\gamma,g_1(\gamma,-))=-}$ is a homotopy. Therefore: $$\begin{aligned}q_1(\gamma,\pi_1\underline{u}')=\;&f_1(\gamma,g_1(\gamma,\beta_1(f_1(\gamma,\pi_1u),f_2(\gamma,\pi_1u,\pi_2u))))\bullet q_1(\gamma,f_1(\gamma,\pi_1u))\\&\bullet \beta_1(f_1(\gamma,\pi_1u),f_2(\gamma,\pi_1u,\pi_2u))^{-1}\\=\;&f_1(\gamma,g_1(\gamma,\beta_1(f_1(\gamma,\pi_1u),f_2(\gamma,\pi_1u,\pi_2u))))\bullet f_1(\gamma,p_1(\gamma,\pi_1u))\\&\bullet \beta_1(f_1(\gamma,\pi_1u),f_2(\gamma,\pi_1u,\pi_2u))^{-1}\\=\;&f_1(\gamma,g_1(\gamma,\beta_1(f_1(\gamma,\pi_1u),f_2(\gamma,\pi_1u,\pi_2u)))\bullet p_1(\gamma,\pi_1u))\\&\bullet \beta_1(f_1(\gamma,\pi_1u),f_2(\gamma,\pi_1u,\pi_2u))^{-1}\\\equiv\;& f_1(\gamma,p)\bullet \beta_1(f_1(\gamma,\pi_1u),f_2(\gamma,\pi_1u,\pi_2u))^{-1}\end{aligned}$$ where the second equality follows because $q_1(\gamma,f_1(\gamma,\pi_1u))=f_1(\gamma,p_1(\gamma,\pi_1u))$ without loss of generality and the third by propositional functoriality. By propositional functoriality and groupoidality, we deduce that: $$\begin{aligned}q_1(\gamma,\pi_1\underline{u}')^*\pi_2\underline{u}'&=f_1(\gamma,p)^* (\beta_1(f_1(\gamma,\pi_1u),f_2(\gamma,\pi_1u,\pi_2u))^{-1})^*\pi_2\underline{u}'\\&=f_1(\gamma,p)^*f_2(\gamma,\pi_1u,\pi_2u)\end{aligned}$$ hence: $$g_2(\gamma,g_1(\gamma,\pi_1\underline{u}'),q_1(\gamma,\pi_1\underline{u}')^*\pi_2\underline{u}')=g_2(\gamma,g_1(\gamma,\pi_1\underline{u}'),f_1(\gamma,p)^*f_2(\gamma,\pi_1u,\pi_2u)).$$ Therefore we are left to observe that: $$g_2(\gamma,g_1(\gamma,\pi_1\underline{u}'),f_1(\gamma,p)^*f_2(\gamma,\pi_1u,\pi_2u))=p^*\pi_2u$$ which follows by the judgement: $$\judge{p:x_1=x_2,y:B(\gamma,x_2)}{g_2(\gamma,x_1,f_1(\gamma,p)^*f_2(\gamma,x_2,y))=p^*y}$$ with $x_1\equiv  g_1(\gamma,\pi_1\underline{u}')$, $x_2\equiv \pi_1u$ and $y\equiv \pi_2u$. This last judgement is true by path elimination on $p$ and since $p_2(\gamma,x_1,y):g_2(\gamma,x_1,f_2(\gamma,x_1,y))=y$.

\medskip

Vice versa, in order to verify that $f^\Sigma(\gamma,u)=\underline{u}'$, where $u\equiv g^\Sigma(\gamma,\underline{u}')$, let us observe that: $$q\equiv f_1(\gamma,\beta_1(\overline{g}_1(\gamma,\underline{u}'),\overline{g}_2(\gamma,\underline{u}')))\bullet q_1(\gamma,\pi_1\underline{u}') :f_1(\gamma,\pi_1u)=\pi_1\underline{u}'$$ where: \[\begin{aligned}
    \overline{g}_1(\gamma,\underline{u}')&\equiv g_1(\gamma,\pi_1\underline{u}')\\
    \overline{g}_2(\gamma,\underline{u}')&\equiv g_2(\gamma,g_1(\gamma,\pi_1\underline{u}'),q_1(\gamma,\pi_1\underline{u}')^*\pi_2\underline{u}').
\end{aligned}\] Hence we are left to verify that: $$f_2(\gamma,\pi_1u,\pi_2u)=q^*\pi_2\underline{u}'$$ as this implies that $f^\Sigma(\gamma,u)=\langle \pi_1\underline{u}',\pi_2\underline{u}'\rangle=\underline{u}'$. This is the case, as: \[\begin{aligned} f_2(\gamma,\pi_1u,\pi_2u)&= f_1(\gamma,\beta_1(\overline{g}_1(\gamma,\underline{u}'),\overline{g}_2(\gamma,\underline{u}')))^*f_2(\gamma,\overline{g}_1(\gamma,\underline{u}'),\overline{g}_2(\gamma,\underline{u}'))\\ &= f_1(\gamma,\beta_1(\overline{g}_1(\gamma,\underline{u}'),\overline{g}_2(\gamma,\underline{u}')))^*q_1(\gamma,\pi_1\underline{u}')^*\pi_2\underline{u}' \\&= (\;f_1(\gamma,\beta_1(\overline{g}_1(\gamma,\underline{u}'),\overline{g}_2(\gamma,\underline{u}')))\bullet q_1(\gamma,\pi_1\underline{u}')\;)^*\pi_2\underline{u}' \\&\equiv q^*\pi_2\underline{u}' \end{aligned}\] where the second equality holds by propositional functoriality and since: $$q_2(\gamma,g_1(\gamma,\pi_1\underline{u}'),q_1(\gamma,\pi_1\underline{u}')^*\pi_2\underline{u}'):f_2(\gamma,\overline{g}_1(\gamma,\underline{u}'),\overline{g}_2(\gamma,\underline{u}'))=q_1(\gamma,\pi_1\underline{u}')^*\pi_2\underline{u}'$$ and the third by propositional functoriality: we are left to verify the first. Since: \begin{center}$\judge{x_1,x_2:A(\gamma),y_1:B(\gamma,x_1),y_2:B(\gamma,x_2),p:x_1=x_2,q:y_1=p^*y_2}$\\$f_2(\gamma,x_1,y_1)= f_1(\gamma,p)^*f_2(\gamma,x_2,y_2)$\end{center} by path elimination on $p$ and $q$, we are done if: $$p\equiv \beta_1(\overline{g}_1(\gamma,\underline{u}'),\overline{g}_2(\gamma,\underline{u}'))\text{ and }q\equiv \beta_2(\overline{g}_1(\gamma,\underline{u}'),\overline{g}_2(\gamma,\underline{u}')).$$

\medskip

Let us summarise the present subsection into the following:

\begin{lem}\label{sigma equiv}
Let us assume that $\gamma:\Gamma$ and $\delta:\Delta$ are contexts and that $(\f;\g)$ is a context homotopy equivalence between them. Moreover, let us assume that we are given judgements: \begin{itemize}
        \item $\judge{\gamma:\Gamma}{A(\gamma):\type}$ \, and \, $\judge{\gamma:\Gamma,x:A(\gamma)}{B(\gamma,x):\type}$
        \item $\judge{\delta:\Delta}{A'(\delta):\type}$ \, and \, $\judge{\delta:\Delta,\underline{x}:A'(\delta)}{B'(\delta,\underline{x}):\type}$
\end{itemize} together with a homotopy equivalence $(f_1;g_1)$: \begin{center}$\judge{\gamma:\Gamma,x:A(\gamma)}{f_1(\gamma,x):A'(\f(\gamma))}$\\ $\judge{\gamma:\Gamma,\underline{x}':A'(\f(\gamma))}{g_1(\gamma,\underline{x}'):A(\gamma)}$\end{center} between $A(\gamma)$ and $A'(\delta)$ relative to $(\f;\g)$, and a homotopy equivalence $(f_2;g_2)$: \begin{center}$\judge{\gamma:\Gamma,x:A(\gamma),y:B(\gamma,x)}{f_2(\gamma,x,y):B'(\f(\gamma),f_1(\gamma,x))}$ \\ $\judge{\gamma:\Gamma,x:A(\gamma),\underline{y}':B'(\f(\gamma),f_1(\gamma,x))}{g_2(\gamma,x,\underline{y}'):B(\gamma,x)}$\end{center} between $B(\gamma,x)$ and $B'(\delta,\underline{x})$ relative to the extension $(\f,\f_{m+1};\g,\g_{n+1})$ of $(\f;\g)$ via $(f_1;g_1)$ (\autoref{context homotopy equivalence expansion}).

In context: \begin{center} $\gamma:\Gamma$ and $u : \Sigma_{x:A(\gamma)}B(\gamma,x)$ and $\underline{u}':[\Sigma_{\underline{x}:A'(\delta)}B'(\delta,\underline{x})](\f(\gamma))$ \end{center} if we name: \begin{center}
    $\judge{\gamma,u}{f^\Sigma(\gamma,u)\equiv\langle f_1(\gamma,\pi_1u),f_2(\gamma, \pi_1u,\pi_2u)\rangle:[\Sigma_{\underline{x}:A'(\delta)}B'(\delta,\underline{x})](\f(\gamma))}$
    
    \medskip
    
    $\judge{\gamma,\underline{u}'}{g^\Sigma(\gamma,\underline{u}')\equiv\langle g_1(\gamma,\pi_1\underline{u}'),g_2(\gamma,g_1(\gamma,\pi_1\underline{u}'),q_1(\gamma,\pi_1\underline{u}')^*\pi_2\underline{u}')\rangle:\Sigma_{x:A(\gamma)}B(\gamma,x)}$
\end{center} then $(f^\Sigma;g^\Sigma)$ is a homotopy equivalence between $\Sigma_{x:A(\gamma)}B(\gamma,x)$ and $\Sigma_{\underline{x}:A'(\delta)}B'(\delta,\underline{x})$ relative to $(\f;\g)$.
\end{lem}

\subsection{Identity types over homotopy equivalent types}\label{section2.3.IV}

Again, let us describe some:

\medskip

\noindent \textbf{Data.} Let us assume that we are given a context homotopy equivalence $(\f;\g)$ between $\gamma:\Gamma$ and $\delta:\Delta$ and judgements: $$\judge{\gamma:\Gamma}{A(\gamma):\type} \textnormal{ and } \judge{\delta:\Delta}{A'(\delta):\type}$$ together with a homotopy equivalence $(f;g)$: \begin{center}$\judge{\gamma:\Gamma,x:A(\gamma)}{f(\gamma,x):A'(\f(\gamma))}$\\ $\judge{\gamma:\Gamma,\underline{x}':A'(\f(\gamma))}{g(\gamma,\underline{x}'):A(\gamma)}$\end{center} between $A(\gamma)$ and $A'(\delta)$ relative to $(\f;\g)$. Moreover let us assume that we are given judgements: \begin{center}$\judge{\gamma:\Gamma}{s_1(\gamma),s_2(\gamma):A(\gamma)}$\\ $\judge{\delta:\Delta}{t_1(\delta),t_2(\delta):A'(\delta)}$\end{center} together with: \begin{center}$\judge{\gamma:\Gamma}{\refl_1(\gamma):f(\gamma,s_1(\gamma))=t_1(\f(\gamma))}$\\ $\judge{\gamma:\Gamma}{\refl_2(\gamma):f(\gamma,s_2(\gamma))=t_2(\f(\gamma)).}$\end{center}

\noindent \textbf{Fact.} We can define a homotopy equivalence $(f^=;g^=)$ between $s_1(\gamma)=s_2(\gamma)$ and $t_1(\delta)=t_2(\delta)$ relative to $(\f;\g)$ as follows, leading to \autoref{id equiv} below.

\medskip

Let us consider homotopies: \begin{center}
    $\judge{\gamma:\Gamma,x:A(\gamma)}{p(\gamma,x):g(\gamma,f(\gamma,x))=x}$\\
    $\judge{\gamma:\Gamma,\underline{x}':A'(\f(\gamma))}{q(\gamma,\underline{x}'):f(\gamma,g(\gamma,\underline{x}'))=\underline{x}'}$
\end{center}

We define: $$\judge{\gamma:\Gamma,p : s_1(\gamma)=s_2(\gamma)}{f^=(\gamma,p)\equiv \refl_1(\gamma)^{-1}\bullet f(\gamma,p)\bullet \refl_2(\gamma):t_1(\f(\gamma))=t_2(\f(\gamma))}$$ and: \begin{center}$\judge{\gamma:\Gamma,\underline{p}':t_1(\f(\gamma))=t_2(\f(\gamma))}{}$\\${g^=(\gamma,\overline{p}')\equiv p(\gamma,s_1(\gamma))^{-1}\bullet g(\gamma,\refl_1(\gamma)\bullet\underline{p}'\bullet \refl_2(\gamma)^{-1})\bullet p(\gamma,s_2(\gamma)):s_1(\gamma)=s_2(\gamma)}.$\end{center} Let us observe that: $$g^=(\gamma,f^=(\gamma,p))=p(\gamma,s_1(\gamma))^{-1}\bullet g(\gamma,f(\gamma,p))\bullet p(\gamma,s_2(\gamma))=p$$ where the first identity type is inhabited by groupoidality and the second since: $$\judge{\gamma:\Gamma}{p(\gamma,-):g(\gamma,f(\gamma,-))=-}.$$ Moreover, by propositional functoriality: $$ f^=(\gamma,g^=(\gamma,\underline{p}'))=a(\gamma)\bullet f(\gamma,g(\gamma,\underline{p}')) \bullet b(\gamma) $$ where: \[\begin{aligned}
    a(\gamma)&\equiv \refl_1(\gamma)^{-1}\bullet f(\gamma,p(\gamma,s_1(\gamma)))^{-1}\bullet f(\gamma,g(\gamma,\refl_1(\gamma)))\\
    &=\refl_1(\gamma)^{-1}\bullet q(\gamma,f(\gamma,s_1(\gamma)))^{-1}\bullet f(\gamma,g(\gamma,\refl_1(\gamma)))\\
    b(\gamma)&\equiv f(\gamma,g(\gamma,\refl_2(\gamma)))^{-1}\bullet f(\gamma,p(\gamma,s_2(\gamma)))\bullet \refl_2(\gamma)\\
    &=f(\gamma,g(\gamma,\refl_2(\gamma)))^{-1}\bullet q(\gamma,f(\gamma,s_2(\gamma)))\bullet \refl_2(\gamma)
\end{aligned}\] where the identity types are inhabited because $\judge{\gamma,x}{f(\gamma,p(\gamma,x))=q(\gamma,f(\gamma,x))}$ without loss of generality. Observe that: \[\begin{aligned}
    a(\gamma)&=q(\gamma,t_1(\f(\gamma)))^{-1}
    \\
    b(\gamma)&=q(\gamma,t_2(\f(\gamma)))
\end{aligned}\] as the diagram: \[\begin{tikzcd}
f(\gamma,s_i(\gamma)) \arrow[rrr, "\refl_i(\gamma)", Rightarrow] &&& t_i(\f(\gamma))\\
f(\gamma,g(\gamma,f(\gamma,s_i(\gamma))))\arrow[u, "q(\gamma\text{,}f(\gamma\text{,}s_i(\gamma)))", Rightarrow] \arrow[rrr, "f(\gamma\text{,}g(\gamma\text{,}\refl_i(\gamma)))", Rightarrow] &&& f(\gamma,g(\gamma,t_i(\f(\gamma)))) \arrow[u, "q(\gamma\text{,}t_i(\f(\gamma)))"', Rightarrow]
\end{tikzcd}\] commutes propositionally for $i=1,2$, since $\judge{\gamma:\Gamma}{q(\gamma,-):f(\gamma,g(\gamma,-))=-}$. We conclude that: $$ f^=(\gamma,g^=(\gamma,\underline{p}'))=q(\gamma,t_1(\f(\gamma)))^{-1}\bullet f(\gamma,g(\gamma,\underline{p}')) \bullet q(\gamma,t_2(\f(\gamma)))=\underline{p}'  $$ by groupoidality and since $\judge{\gamma:\Gamma}{q(\gamma,-):f(\gamma,g(\gamma,-))=-}$. We summarise this fact into the following:

\begin{lem}\label{id equiv}
Let us assume that we are given a context homotopy equivalence $(\f;\g)$ between $\gamma:\Gamma$ and $\delta:\Delta$ and judgements: $$\judge{\gamma:\Gamma}{A(\gamma):\type} \textnormal{ and } \judge{\delta:\Delta}{A'(\delta):\type}$$ together with a homotopy equivalence $(f;g)$: \begin{center}$\judge{\gamma:\Gamma,x:A(\gamma)}{f(\gamma,x):A'(\f(\gamma))}$\\ $\judge{\gamma:\Gamma,\underline{x}':A'(\f(\gamma))}{g(\gamma,\underline{x}'):A(\gamma)}$\end{center} between $A(\gamma)$ and $A'(\delta)$ relative to $(\f;\g)$. Moreover let us assume that we are given judgements: \begin{center}$\judge{\gamma:\Gamma}{s_1(\gamma),s_2(\gamma):A(\gamma)}$\\ $\judge{\delta:\Delta}{t_1(\delta),t_2(\delta):A'(\delta)}$\end{center} together with: \begin{center}$\judge{\gamma:\Gamma}{\refl_1(\gamma):f(\gamma,s_1(\gamma))=t_1(\f(\gamma))}$\\ $\judge{\gamma:\Gamma}{\refl_2(\gamma):f(\gamma,s_2(\gamma))=t_2(\f(\gamma)).}$\end{center}

If we name: \begin{center} $\judge{\gamma:\Gamma,p : s_1(\gamma)=s_2(\gamma)}{f^=(\gamma,p)\equiv \refl_1(\gamma)^{-1}\bullet f(\gamma,p)\bullet \refl_2(\gamma):t_1(\f(\gamma))=t_2(\f(\gamma))}$

\smallskip

$\judge{\gamma:\Gamma,\underline{p}':t_1(\f(\gamma))=t_2(\f(\gamma))}{g^=(\gamma,\underline{p}')\equiv p(\gamma,s_1(\gamma))^{-1}\bullet g(\gamma,\refl_1(\gamma)\bullet\underline{p}'\bullet \refl_2(\gamma)^{-1})\bullet p(\gamma,s_2(\gamma)):s_1(\gamma)=s_2(\gamma)}$
\end{center} then $(f^=;g^=)$ is a homotopy equivalence between $s_1(\gamma)=s_2(\gamma)$ and $t_1(\delta)=t_2(\delta)$ relative to $(\f;\g)$.
\end{lem}

\subsection{Canonical homotopy equivalences between types}\label{canonical homotopy equivalences}\label{section2.3.V}

In \autoref{section2.3.VI} we are going to use \autoref{context homotopy equivalence expansion} to define a specific class of context homotopy equivalences. However, in order to do so, we first need to identify a particular notion of homotopy equivalence between types relative to a given context homotopy equivalence: we are going to apply \autoref{context homotopy equivalence expansion} to these ones only. In this subsection we inductively present this notion.

\medskip

Let $\gamma : \Gamma$ and $\delta :\Delta$ be contexts together with a context homotopy equivalence $(\f;\g)$ and let $\judge{\gamma:\Gamma}{S(\gamma):\type}$ and $\judge{\delta :\Delta}{T(\delta):\type}$ have h-propositional identities. We provide a list of inductive clauses determining the family  of the \textbf{canonical homotopy equivalences} $\equivpair{\phi}{\psi}$ between $S(\gamma)$ and $T(\delta)$ relative to $(\f;\g)$:
\begin{itemize}[align=left]
    \item[\texttt{(a)}] If
    $S(\gamma)\equiv S$ and $T(\delta)\equiv S(\gamma)\equiv S$, then the identity of $S$ (as judgement of the form $\judge{\gamma : \Gamma, s:S}{s:S}$) and itself
    constitute a canonical homotopy equivalence $\equivpair{\phi}{\psi}$ between $S(\gamma)$ and $T(\delta)$ relative to $(\f;\g)$. Observe that this homotopy equivalence between $S$ and itself is in fact \textit{relative to} $(\f;\g)$, as: \begin{itemize}
    \item $\judge{\gamma:\Gamma,s:S}{\phi(\gamma,s)\equiv s:S\equiv S(\f(\gamma))}$
    \item $\judge{\gamma:\Gamma,s:S\equiv S(\f(\gamma))}{\psi(\gamma,s)\equiv s:S}$.
\end{itemize}
    \item[\texttt{(b)}] If $S(\gamma)\equiv \Pi_{x:A(\gamma)}B(\gamma,x)$ and $T(\delta)\equiv \Pi_{\underline{x}:A'(\delta)}B'(\delta,\underline{x})$ for some judgements: \begin{itemize}
        \item $\judge{\gamma:\Gamma}{A(\gamma):\type}$ \, and \, $\judge{\gamma:\Gamma,x:A(\gamma)}{B(\gamma,x):\type}$
        \item $\judge{\delta:\Delta}{A'(\delta):\type}$ \, and \, $\judge{\delta:\Delta,\underline{x}:A'(\delta)}{B'(\delta,\underline{x}):\type}$
    \end{itemize} and if there are a canonical homotopy equivalence $\equivpair{f_1}{g_1}$: \begin{center}$\judge{\gamma:\Gamma,x:A(\gamma)}{f_1(\gamma,x):A'(\f(\gamma))}$\\ $\judge{\gamma:\Gamma,\underline{x}':A'(\f(\gamma))}{g_1(\gamma,\underline{x}'):A(\gamma)}$\end{center} between $A(\gamma)$ and $A'(\delta)$ relative to $(\f;\g)$ and a canonical homotopy equivalence $\equivpair{f_2}{g_2}$: \begin{center}$\judge{\gamma:\Gamma,x:A(\gamma),y:B(\gamma,x)}{f_2(\gamma,x,y):B'(\f(\gamma),f_1(\gamma,x))}$ \\ $\judge{\gamma:\Gamma,x:A(\gamma),\underline{y}':B'(\f(\gamma),f_1(\gamma,x))}{g_2(\gamma,x,\underline{y}'):B(\gamma,x)}$\end{center} between $B(\gamma,x)$ and $B'(\delta,\underline{x})$ relative to the extension $(\f,\f_{m+1};\g,\g_{n+1})$ of $(\f;\g)$ via $(f_1;g_1)$ (\autoref{context homotopy equivalence expansion}), then the homotopy equivalence $( f^\Pi;g^\Pi)$ of \autoref{pi equiv} is a canonical homotopy equivalence $\equivpair{\phi}{\psi}$ between $S(\gamma)$ and $T(\delta)$ relative to $(\f;\g)$.
    \item[\texttt{(c)}] If $S(\gamma)\equiv \Sigma_{x:A(\gamma)}B(\gamma,x)$ and $T(\delta)\equiv \Sigma_{\underline{x}:A'(\delta)}B'(\delta,\underline{x})$ for some judgements: \begin{itemize}
        \item $\judge{\gamma:\Gamma}{A(\gamma):\type}$ \, and \, $\judge{\gamma:\Gamma,x:A(\gamma)}{B(\gamma,x):\type}$
        \item $\judge{\delta:\Delta}{A'(\delta):\type}$ \, and \, $\judge{\delta:\Delta,\underline{x}:A'(\delta)}{B'(\delta,\underline{x}):\type}$
    \end{itemize} and if there are a canonical homotopy equivalence $\equivpair{f_1}{g_1}$: \begin{center}$\judge{\gamma:\Gamma,x:A(\gamma)}{f_1(\gamma,x):A'(\f(\gamma))}$\\ $\judge{\gamma:\Gamma,\underline{x}':A'(\f(\gamma))}{g_1(\gamma,\underline{x}'):A(\gamma)}$\end{center} between $A(\gamma)$ and $A'(\delta)$ relative to $(\f;\g)$ and a canonical homotopy equivalence $\equivpair{f_2}{g_2}$: \begin{center}$\judge{\gamma:\Gamma,x:A(\gamma),y:B(\gamma,x)}{f_2(\gamma,x,y):B'(\f(\gamma),f_1(\gamma,x))}$ \\ $\judge{\gamma:\Gamma,x:A(\gamma),\underline{y}':B'(\f(\gamma),f_1(\gamma,x))}{g_2(\gamma,x,\underline{y}'):B(\gamma,x)}$\end{center} between $B(\gamma,x)$ and $B'(\delta,\underline{x})$ relative to the extension $(\f,\f_{m+1};\g,\g_{n+1})$ of $(\f;\g)$ via $(f_1;g_1)$ (\autoref{context homotopy equivalence expansion}), then the homotopy equivalence $( f^\Sigma;g^\Sigma)$ of \autoref{sigma equiv} is a canonical homotopy equivalence $\equivpair{\phi}{\psi}$ between $S(\gamma)$ and $T(\delta)$ relative to $(\f;\g)$.
    \item[\texttt{(d)}] If $S(\gamma)\equiv s_1(\gamma)=s_2(\gamma)$ and $T(\delta)\equiv t_1(\delta)=t_2(\delta)$ for some judgements: \begin{center}
    $\judge{\gamma:\Gamma}{A(\gamma):\type}$ \, and \, $\judge{\delta:\Delta}{A'(\delta):\type}$
    \end{center} and some judgements:
    \begin{center}
        $\judge{\gamma:\Gamma}{s_1(\gamma),s_2(\gamma):A(\gamma)}$ \, and \, $\judge{\delta:\Delta}{t_1(\delta),t_2(\delta):A'(\delta)}$\end{center} and if there are a canonical homotopy equivalence $\equivpair{f}{g}$: \begin{center}$\judge{\gamma:\Gamma,x:A(\gamma)}{f(\gamma,x):A'(\f(\gamma))}$\\ $\judge{\gamma:\Gamma,\underline{x}':A'(\f(\gamma))}{g(\gamma,\underline{x}'):A(\gamma)}$\end{center} between $A(\gamma)$ and $A'(\delta)$ relative to $(\f;\g)$ and judgements: \begin{center}$\judge{\gamma:\Gamma}{\refl_1(\gamma):f(\gamma,s_1(\gamma))=t_1(\f(\gamma))}$\\ $\judge{\gamma:\Gamma}{\refl_2(\gamma):f(\gamma,s_2(\gamma))=t_2(\f(\gamma)).}$\end{center} then the homotopy equivalence $(f^=,g^=)$ of \autoref{id equiv} is a canonical homotopy equivalence $\equivpair{\phi}{\psi}$ between $S(\gamma)$ and $T(\delta)$ relative to $(\f;\g)$.
\end{itemize}

\subsection{Canonical homotopy equivalences between contexts}\label{canonical context homotopy equivalences}\label{section2.3.VI}

In this subsection we recursively define a family of context homotopy equivalences that we refer to as the \textit{canonical} ones. In fact, when we make the syntax of $\ptt$ into a model of $\ett$, we need to define equivalence classes of contexts modulo the context homotopy equivalences of this particular shape that we are going to identify.

As mentioned before, this proof strategy, based on defining a family of canonical equivalences between contexts---up to which contexts are identified---appears in a work by Hofmann \cite{phdhofmann,hofmannconservativity} in order to obtain an analogous conservativity result. Moreover, such a notion of canonical equivalence appears---in a similar formulation---in a work by Maietti \cite{MAIETTI2009319} in order to define an interpretation of the extensional level of Minimalist Foundation within the intensional one, as well as in one by Contente and Maietti \cite{contentemaietti} in order to interpret the former within Homotopy Type Theory, and in one by Maietti and Sabelli \cite{maiettisabelli} to define an interpretation within the former of itself extendend with an extensionality axiom for its propositions.

\medskip

Let $\gamma : \Gamma$ and $\delta : \Delta$ be two contexts. We give a list of inductive clauses determining the family of the \textbf{canonical context homotopy equivalences} $\equivpair{\ci}{\di}$ between them:
\begin{itemize}[align=left]
    \item[\texttt{(1)}] If $\gamma:\Gamma$ and $\delta:\Delta$ are the empty context then the empty list (as context morphism $\Gamma\to\Delta$) and itself (as context morphism $\Delta\to\Gamma$) constitute a canonical context homotopy equivalence $\equivpair{\ci}{\di}$.
    \item[\texttt{(2)}] If $\gamma  : \Gamma$ and $\delta:\Delta$ are of the form $\gamma'  : \Gamma',x : A(\gamma')$ and $\delta' :  \Delta',\underline{x}: A'(\delta')$ respectively and: \begin{center}$\judge{\gamma':\Gamma'}{\f(\gamma'):\Delta'}$\\ $\judge{\delta':\Delta'}{\g(\delta'):\Gamma'}$\end{center} is a canonical context homotopy equivalence $\equivpair{\f}{\g}$ and: \begin{center} $\judge{\gamma':\Gamma',x:A(\gamma')}{f(\gamma',x):A'(\f(\gamma'))}$\\ $\judge{\gamma':\Gamma',\underline{x}':A'(\f(\gamma'))}{g(\gamma',\underline{x}'):A(\gamma')}$\end{center} is a canonical homotopy equivalence $\equivpair{f}{g}$ between $A(\gamma')$ and $A'(\delta')$ relative to $\equivpair{\f}{\g}$ (in the sense of \autoref{section2.3.V}) then the extension: \begin{center} $\judge{\gamma' : \Gamma',  x:A(\gamma')}{(\f,\f_{m+1})(\gamma',x):\Delta',A'}$\\ $\judge{\delta' : \Delta', \underline{x}:A'(\delta')}{(\g,\g_{n+1})(\delta',\underline{x}):\Gamma',A}$\end{center} of $\equivpair{\f}{\g}$ via $\equivpair{f}{g}$ (\autoref{context homotopy equivalence expansion}) is a canonical context homotopy equivalence $\equivpair{\ci}{\di}$.
\end{itemize}

We end the current section with the following:

\begin{rem}\label{crucial example}
    Let $A:\type$, $a,b:A$, and $q:a=b$ be atomic judgements of PTT. By \texttt{(a)} of \autoref{section2.3.V}, the pair: $$\equivpair{\;\;\judge{x:A}{x:A}\;\;}{\;\;\judge{x:A}{x:A}\;\;}$$ constitutes a canonical homotopy equivalence between $A$ and itself relative to the canonical context homotopy equivalence of \texttt{(1)}. Therefore, by \texttt{(b)} of \autoref{section2.3.V}, the pair: $$\equivpair{\;\;\judge{p:a=a}{\refl(a)^{-1}\bullet p \bullet q:a=b}\;\;}{\;\;\judge{p':a=b}{\refl(a)\bullet p' \bullet q^{-1}:a=a}\;\;}$$ constitutes a canonical homotopy equivalence between the types $a=a$ and $a=b$ relative to the canonical context homotopy equivalence of \texttt{(1)}. By \texttt{(2)} we obtain a canonical context homotopy equivalence: $$\equivpair{\;\;\judge{p:a=a}{\refl(a)^{-1}\bullet p \bullet q:a=b}\;\;}{\;\;\judge{p':a=b}{\refl(a)\bullet p' \bullet q^{-1}:a=a}\;\;}$$ between the context $p:a=a$ and the context $p':a=b$.

    However, assuming that $a\equiv b$ (and w.l.o.g. that $p'\equiv p$) this homotopy equivalence is homotopic---i.e. propositionally equal---to the identity over the context $p:a=a$ if and only if the type $q=\refl(a)$ is inhabited---cf. \autoref{only one canonical context homotopy equivalence} on contexts with h-propositional identities.

    We will use this fundamental canonical context homotopy equivalence, which we have presented here in a particular case, in \autoref{section2.6.II} to show that the syntax of $\ptt$ with contexts identified modulo these canonical equivalences constitutes a model of extensional identity types.

\end{rem}

\section{Properties of the canonical homotopy equivalences}\label{section2.4}

In this section we prove those properties of the family of the canonical context homotopy equivalences (defined in \autoref{section2.3.VI}) that are needed in order to define a model of $\ett$ starting from the syntax of $\ptt$. Our approach is the following: we start by proving properties (e.g. the fact that an equivalence relation is induced between the contexts of $\ptt$) for the general family of canonical equivalences; secondly, we restrict to a smaller family of contexts, called \textit{contexts with h-propositional identities}, and we only consider canonical equivalences between them, so that additional properties (e.g. the uniqueness---up to homotopy---of a canonical equivalence between two given contexts) are satisfied.

\subsection{Properties of the family of canonical equivalences}\label{section2.4.I}

We remind that we use the notation $(\;\cdot\;;\;\cdot\;)$ to indicate any homotopy equivalence between contexts or types in context, reserving the symbol $\equivpair{\cdot}{\cdot}$ for the canonical ones. We start our list of results from the following observation:

\begin{lem}
Let $\gamma:\Gamma$ be a context $\gamma_1 : \Gamma_1,\gamma_2: \Gamma_2(\gamma_1),...,\gamma_n:\Gamma_n(\gamma_1,...,\gamma_{n-1})$ and let $\delta:\Delta$ be a context $\delta_1 : \Delta_1,\delta_2: \Delta_2(\delta_1),...,\delta_m:\Delta_m(\delta_1,...,\delta_{m-1})$. If there is a canonical context homotopy equivalence $\equivpair{\ci}{\di}$ between $\gamma:\Gamma$ and $\delta:\Delta$ then $n=m$.

\proof
By induction on the complexity of $\equivpair{\ci}{\di}$. If $\equivpair{\ci}{\di}$ is of the form \texttt{(1)} then $n=0=m$ and we are done. If $\equivpair{\ci}{\di}$ is of the form \texttt{(2)} for some canonical context homotopy equivalence $\equivpair{\f}{\g}$ and some canonical homotopy equivalence $\equivpair{f}{g}$ relative to $\equivpair{\f}{\g}$, then by inductive hypothesis $n-1=m-1$ hence we are done.
\endproof
\end{lem}

Secondly, we prove the properties of reflexivity, symmetry and transitivity of the relation between contexts induced by the family of the canonical equivalences between them.

\begin{prop}[Reflexivity]\label{there are context identities}
Let $\gamma : \Gamma$ be a context: $$\gamma_1 : \Gamma_1,\gamma_2: \Gamma_2(\gamma_1),...,\gamma_n:\Gamma_n(\gamma_1,...,\gamma_{n-1}).$$ Then there is a canonical context homotopy equivalence $\equivpair{\ci}{\di}$ from $\gamma$ to $\gamma$ such that $\judge{\gamma}{\ci(\gamma)=\gamma}$ (\,i.e. $\judge{\gamma}{\di(\gamma)=\gamma}$\,).

\proof
By induction on the length $n$ of $\gamma:\Gamma$. If $n=0$ then we are done by \texttt{(1)}. Otherwise, let $\gamma':\Gamma'$ be the context $\gamma_1,\gamma_2,...,\gamma_{n-1}$. By inductive hypothesis, there is a canonical context homotopy equivalence $\equivpair{\f}{\g}$ between $\gamma'$ and itself such that $\judge{\gamma'}{\alpha(\gamma'):\f(\gamma')=\gamma'}$ (\,i.e. $\judge{\gamma'}{\g(\gamma')=\gamma'}$\,).

Let us pretend that we know that there is a canonical homotopy equivalence $\equivpair{f}{g}$ between $\Gamma_n(\gamma')$ and itself relative to $\equivpair{\f}{\g}$ and such that $\judge{\gamma',\gamma_n:\Gamma_n(\gamma')}{f(\gamma',\gamma_n)=\alpha(\gamma')^*\gamma_n}$ (\,i.e. $\judge{\gamma',\gamma_n}{g(\gamma',\alpha(\gamma')^*\gamma_n)=\gamma_n}$\,). Then the extension: \begin{center} $\judge{\gamma' : \Gamma',  \gamma_n:\Gamma_n(\gamma')}{(\f,\f_{m+1})(\gamma',\gamma_n):\Gamma',\Gamma_n}$\\ $\judge{\gamma' : \Gamma',  \gamma_n:\Gamma_n(\gamma')}{(\g,\g_{n+1})(\gamma',\gamma_n):\Gamma',\Gamma_n}$\end{center} of $\equivpair{\f}{\g}$ via $\equivpair{f}{g}$ (\autoref{context homotopy equivalence expansion}) satisfies $(\f,\f_{m+1})(\gamma)=\gamma$ (\,hence $(\g,\g_{n+1})(\gamma)=\gamma$\,). Moreover, it is a canonical homotopy equivalence $\equivpair{\ci}{\di}$ by \texttt{(2)}.

Hence we are done if there is a canonical homotopy equivalence $\equivpair{f}{g}$ between $\Gamma_n(\gamma')$ and itself relative to $\equivpair{\f}{\g}$ and such that $\judge{\gamma',\gamma_n:\Gamma_n(\gamma')}{f(\gamma',\gamma_n)=\alpha(\gamma')^*\gamma_n}$. But this is true because of the following \autoref{there are identities}.
\endproof
\end{prop}

\begin{lem}\label{there are identities}
Let $\gamma:\Gamma$ be a context and let $(\f;\g)$ be a context homotopy equivalence between $\gamma$ and itself such that $\judge{\gamma}{\alpha(\gamma):\f(\gamma)=\gamma}$ ( and equivalently $\judge{\gamma}{\g(\gamma)=\gamma}$ ). Moreover, let us assume that $\judge{\gamma:\Gamma}{S(\gamma):\type}$. Then there exists a canonical homotopy equivalence $\equivpair{\phi}{\psi}$ between $S(\gamma)$ and itself relative to $(\f;\g)$ and such that: \begin{center} $\judge{\gamma,s:S(\gamma)}{\phi(\gamma,s)=\alpha(\gamma)^*s}$ (\,i.e. $\judge{\gamma,s}{\psi(\gamma,\alpha(\gamma)^*s)=s}$\,). \end{center}

\proof
By induction on the complexity of the type $S(\gamma)$.

If $S(\gamma)\equiv S$ then by \texttt{(a)} the pair: $$(\judge{\gamma,s:S}{\phi(\gamma,s)\equiv s:S};\;\judge{\gamma,s:S}{\psi(\gamma,s)\equiv s:S})$$ is a canonical homotopy equivalence between $S$ and itself relative to $(\f;\g)$. Moreover, let us observe that: $$\judge{\gamma,s}{\phi(\gamma,s)\equiv s=\alpha(\gamma)^*s}$$---by (generalised) path induction on $\alpha(\gamma)$---hence we are done.

If $S(\gamma)$ is of the form $\Pi_{x:A(\gamma)}B(\gamma,x)$ for some judgements $\judge{\gamma:\Gamma}{A(\gamma):\type}$ \, and \, $\judge{\gamma:\Gamma,x:A(\gamma)}{B(\gamma,x):\type}$, then by inductive hypothesis there are a canonical homotopy equivalence $\equivpair{f_1}{g_1}$: \begin{center}$\judge{\gamma:\Gamma,x:A(\gamma)}{f_1(\gamma,x):A(\f(\gamma))}$\\ $\judge{\gamma:\Gamma,x':A(\f(\gamma))}{g_1(\gamma,x'):A(\gamma)}$\end{center} between $A(\gamma)$ and itself relative to $(\f;\g)$ such that $\judge{\gamma,x}{\alpha_1(\gamma,x):f_1(\gamma,x)=\alpha(\gamma)^*x}$ and equivalently $\judge{\gamma,x}{g_1(\gamma,\alpha(\gamma)^*x)=x}$ and a canonical homotopy equivalence $\equivpair{f_2}{g_2}$: \begin{center}$\judge{\gamma:\Gamma,x:A(\gamma),y:B(\gamma,x)}{f_2(\gamma,x,y):B(\f(\gamma),f_1(\gamma,x))}$ \\ $\judge{\gamma:\Gamma,x:A(\gamma),y':B(\f(\gamma),f_1(\gamma,x))}{g_2(\gamma,x,\underline{y}'):B(\gamma,x)}$\end{center} between $B(\gamma,x)$ and itself relative to the extension $(\f,\f_{m+1};\g,\g_{n+1})$ of $(\f;\g)$ via $(f_1;g_1)$ (\autoref{context homotopy equivalence expansion}) such that $\judge{\gamma,x,y}{f_2(\gamma,x,y)=(\alpha(\gamma),\alpha_1(\gamma,x))^*y}$ and equivalently: $$\judge{\gamma,x,y}{g_2(\gamma,x,(\alpha(\gamma),\alpha_1(\gamma,x))^*y)=y}.$$ The homotopy equivalence $(f^\Pi;g^\Pi)$ of \autoref{pi equiv} is a canonical homotopy equivalence relative to $(\f;\g)$. We are left to verify that: $$\judge{\gamma,z:\Pi_{x:A(\gamma)}B(\gamma,x)}{g^\Pi(\gamma,\alpha(\gamma)^*z)=z}.$$ Let us observe that: $$\judge{\gamma,z,x}{\ev(\alpha(\gamma)^*z,f_1(\gamma,x))=\ev(\alpha(\gamma)^*z,\alpha(\gamma)^*x)=(\alpha(\gamma),\alpha_1(\gamma,x))^*\ev(z,x)}$$ where the second equality holds by multiple (generalised) path induction on $(\alpha(\gamma),\alpha_1(\gamma,x))$. Hence, by propositional functoriality: $$\judge{\gamma,z,x}{g_2(\gamma,x,\ev(\alpha(\gamma)^*z,f_1(\gamma,x)))=g_2(\gamma,x,(\alpha(\gamma),\alpha_1(\gamma,x))^*\ev(z,x))=\ev(z,x)}$$ where the second equality holds by substitution into the judgement: $$\judge{\gamma,x,y}{g_2(\gamma,x,(\alpha(\gamma),\alpha_1(\gamma,x))^*y)=y}.$$ Finally, by propositional functoriality and propositional $\eta$-expansion: $$\judge{\gamma,z}{g^\Pi(\gamma,\alpha(\gamma)^*z)\equiv \lambda x.g_2(\gamma,x,\ev(\alpha(\gamma)^*z,f_1(\gamma,x)))=\lambda x. \ev(z,x)=z}$$ and we are done.

If $S(\gamma)$ is of the form $\Sigma_{x:A(\gamma)}B(\gamma,x)$ for some judgements $\judge{\gamma:\Gamma}{A(\gamma):\type}$ \, and \, $\judge{\gamma:\Gamma,x:A(\gamma)}{B(\gamma,x):\type}$, then by inductive hypothesis there are a canonical homotopy equivalence $\equivpair{f_1}{g_1}$: \begin{center}$\judge{\gamma:\Gamma,x:A(\gamma)}{f_1(\gamma,x):A(\f(\gamma))}$\\ $\judge{\gamma:\Gamma,x':A(\f(\gamma))}{g_1(\gamma,x'):A(\gamma)}$\end{center} between $A(\gamma)$ and itself relative to $(\f;\g)$ such that $\judge{\gamma,x}{\alpha_1(\gamma,x):f_1(\gamma,x)=\alpha(\gamma)^*x}$ and equivalently $\judge{\gamma,x}{g_1(\gamma,\alpha(\gamma)^*x)=x}$ and a canonical homotopy equivalence $\equivpair{f_2}{g_2}$: \begin{center}$\judge{\gamma:\Gamma,x:A(\gamma),y:B(\gamma,x)}{f_2(\gamma,x,y):B(\f(\gamma),f_1(\gamma,x))}$ \\ $\judge{\gamma:\Gamma,x:A(\gamma),y':B(\f(\gamma),f_1(\gamma,x))}{g_2(\gamma,x,\underline{y}'):B(\gamma,x)}$\end{center} between $B(\gamma,x)$ and itself relative to the extension $(\f,\f_{m+1};\g,\g_{n+1})$ of $(\f;\g)$ via $(f_1;g_1)$ (\autoref{context homotopy equivalence expansion}) such that $\judge{\gamma,x,y}{f_2(\gamma,x,y)=(\alpha(\gamma),\alpha_1(\gamma,x))^*y}$ and equivalently: $$\judge{\gamma,x,y}{g_2(\gamma,x,(\alpha(\gamma),\alpha_1(\gamma,x))^*y)=y}.$$ The homotopy equivalence $(f^\Sigma;g^\Sigma)$ of \autoref{sigma equiv} is a canonical homotopy equivalence relative to $(\f;\g)$, hence we are left to verify that: $$\judge{\gamma,u:\Sigma_{x:A(\gamma)}B(\gamma,x)}{f^\Sigma(\gamma,u)\equiv\langle f_1(\gamma,\pi_1u),f_2(\gamma, \pi_1u,\pi_2u)\rangle=\alpha(\gamma)^*u}.$$ As $\alpha_1(\gamma,\pi_1u):f_1(\gamma,\pi_1u)=\alpha(\gamma)^*\pi_1u$ and: $$f_2(\gamma,\pi_1u,\pi_2u)=(\alpha(\gamma),\alpha_1(\gamma,\pi_1u))^*\pi_2u=\alpha_1(\gamma,\pi_1u)^*(\alpha(\gamma)^*\pi_2u)$$---where the last equality follows by (generalised) path induction on $\alpha_1(\gamma,\pi_1u)$---it is in fact the case that: $$\judge{\gamma,u}{f^\Sigma(\gamma,u)=\langle \alpha(\gamma)^*\pi_1u,\alpha(\gamma)^*\pi_2u\rangle=\alpha(\gamma)^*\langle \pi_1u,\pi_2u\rangle=\alpha(\gamma)^*u}$$ where the second equality follows by multiple (generalised) path induction on $\alpha(\gamma)$ and the third by propositional functoriality and propositional $\eta$-expansion.

If $S(\gamma)$ is of the form $s_1(\gamma)=s_2(\gamma)$ for some judgement $\judge{\gamma:\Gamma}{A(\gamma):\type}$, then by inductive hypothesis there is a canonical homotopy equivalence $\equivpair{f}{g}$: \begin{center}$\judge{\gamma:\Gamma,x:A(\gamma)}{f(\gamma,x):A(\f(\gamma))}$\\ $\judge{\gamma:\Gamma,x':A(\f(\gamma))}{g(\gamma,x'):A(\gamma)}$\end{center} between $A(\gamma)$ and itself relative to $(\f;\g)$ such that $\judge{\gamma,x}{\alpha_1(\gamma,x):f(\gamma,x)=\alpha(\gamma)^*x}$ and equivalently $\judge{\gamma,x}{g_1(\gamma,\alpha(\gamma)^*x)=x}$. For $i=1,2$ let us observe that: $$\judge{\gamma}{\alpha_1(\gamma,s_i(\gamma)):f(\gamma,s_i(\gamma))=\alpha(\gamma)^*s_i(\gamma)}\textnormal{ and }\judge{\gamma}{s_i(\alpha(\gamma)):s_i(\f(\gamma))=\alpha(\gamma)^*s_i(\gamma)}$$ hence: $$\judge{\gamma:\Gamma}{r_i(\gamma)\equiv \alpha_1(\gamma,s_i(\gamma))\bullet s_i(\alpha(\gamma))^{-1} :f(\gamma,s_i(\gamma))=s_i(\f(\gamma))}$$ and then the homotopy equivalence $(f^=,g^=)$ of \autoref{id equiv} is a canonical homotopy equivalence relative to $(\f;\g)$. We are left to show that $\judge{\gamma, p:s_1(\gamma)=s_2(\gamma)}{f^=(\gamma,p)=\alpha(\gamma)^*p}$. In fact: $$\begin{aligned}\judge{\gamma,p}{f^=(\gamma,p)&=s_1(\alpha(\gamma))\bullet \alpha_1(\gamma,s_1(\gamma))^{-1}\bullet f(\gamma,p)\bullet\alpha_1(\gamma,s_2(\gamma))\bullet s_2(\alpha(\gamma))^{-1}\\&=s_1(\alpha(\gamma))\bullet\alpha(\gamma)^*p\bullet s_2(\alpha(\gamma))^{-1}}\\&=\alpha(\gamma)^*p\end{aligned}$$ by propositional groupoidality and since $\judge{\gamma,x}{\alpha_1(\gamma,x):f(\gamma,x)=\alpha(\gamma)^*x}$. Here the last of these equalities follows by multiple (generalised) path induction on $\alpha(\gamma)$ and by propositional groupoidality.
\endproof
\end{lem}

\begin{prop}[Symmetry]\label{symmetry}
If $\equivpair{\ci}{\di}$ is a canonical context homotopy equivalence between $\gamma:\Gamma$ and $\delta:\Delta$, then there is a canonical context homotopy equivalence $\equivpair{\di'}{\ci'}$ between $\delta$ and $\gamma$ such that $\judge{\gamma:\Gamma}{\ci(\gamma)=\ci'(\gamma)}$ and equivalently $\judge{\delta:\Delta}{\di(\delta)=\di'(\delta)}$.

\proof
By induction on the complexity of $\equivpair{\ci}{\di}$. If $\equivpair{\ci}{\di}$ is of the form \texttt{(1)} then we are done. If $\equivpair{\ci}{\di}$ is of the form \texttt{(2)} for some canonical context homotopy equivalence $\equivpair{\f}{\g}$ between $\gamma':\Gamma'$ and $\delta':\Delta'$ and some canonical homotopy equivalence $\equivpair{f}{g}$ between $\judge{\gamma'}{A(\gamma')}$ and $\judge{\delta'}{A'(\delta')}$ relative to $\equivpair{\f}{\g}$ \begin{center}---hence $\gamma\equiv\gamma',x:A(\gamma') $ and $\delta\equiv \delta',\underline{x}:A'(\delta')$---\end{center} then by inductive hypothesis there is a canonical context homotopy equivalence $\equivpair{\g'}{\f'}$ such that $\judge{\gamma':\Gamma'}{\f(\gamma')=\f'(\gamma')}$ and equivalently $\judge{\delta':\Delta'}{\alpha(\delta'):\g(\delta')=\g'(\delta')}$. Moreover, by the following \autoref{symmetry of canonical equiv} there is a canonical homotopy equivalence $\equivpair{g'}{f'}$: \begin{center}$\judge{\delta':\Delta',\underline{x}:A'(\delta)}{g'(\delta',\underline{x}):A(\g'(\delta'))}$\\ $\judge{\delta':\Delta',x':A(\g'(\delta'))}{f'(\delta',x'):A'(\delta')}$\end{center} relative to $(\g';\f')$ such that: \begin{center}$\judge{\delta':\Delta',\underline{x}:A'(\delta')}{\alpha_1(\delta',\underline{x}):g(\g(\delta'),\q(\delta')^*\underline{x})=\alpha(\delta')^*g'(\delta',\underline{x})}$\\ (\,i.e. $\judge{\delta':\Delta',x':A(\g(\delta'))}{f(\g(\delta'),\alpha(\delta')^*x')=\q(\delta')^*f'(\delta',x')}$\,)\end{center} hence by \texttt{(2)} the extension $(\g',\g'_{n+1};\f',\f'_{m+1})$ of $\equivpair{\g'}{\f'}$ via $\equivpair{g'}{f'}$ (\autoref{context homotopy equivalence expansion}) is a canonical context homotopy equivalence $\equivpair{\di'}{\ci'}$ between $\delta',\underline{x}\equiv \delta$ and $\gamma',x\equiv \gamma$ and: $$\judge{\delta',\underline{x}}{(\alpha(\delta'),\alpha_1(\delta',\underline{x})):(\g,\g_{n+1})(\delta',\underline{x})=(\g',\g'_{n+1})(\delta',\underline{x})}$$ that is: $$\judge{\delta}{\di(\delta)\equiv(\g,\g_{n+1})(\delta)=(\g',\g'_{n+1})(\delta)\equiv \di'(\delta)}$$ hence we are done. \endproof
\end{prop}

\begin{lem}\label{symmetry of canonical equiv}
Suppose that we are given a canonical homotopy equivalence $(\f;\g)$ between $\gamma: \Gamma$ and $\delta: \Delta$ and suppose that there is a canonical homotopy equivalence $(\g';\f')$ between $\delta$ and $\gamma$ such that $\judge{\gamma}{\f(\gamma)=\f'(\gamma)}$ and equivalently $\judge{\delta}{\alpha(\delta):\g(\delta)=\g'(\delta)}$.

If we are given judgements: \begin{center}
        $\judge{\gamma:\Gamma}{S(\gamma):\type}$ \, and \, $\judge{\delta:\Delta}{T(\delta):\type}$
\end{center} and if there is a canonical homotopy equivalence $\equivpair{\phi}{\psi}$: \begin{center}$\judge{\gamma:\Gamma,x:S(\gamma)}{\phi(\gamma,x):T(\f(\gamma))}$\\ $\judge{\gamma:\Gamma,\underline{x}':T(\f(\gamma))}{\psi(\gamma,\underline{x}'):S(\gamma)}$\end{center} between $S(\gamma)$ and $T(\delta)$ relative to a context homotopy equivalence $(\f;\g)$, then there is a canonical homotopy equivalence $\equivpair{\psi'}{\phi'}$: \begin{center}$\judge{\delta:\Delta,\underline{x}:T(\delta)}{\psi'(\delta,\underline{x}):S(\g'(\delta))}$\\ $\judge{\delta:\Delta,x':S(\g'(\delta))}{\phi'(\delta,x'):T(\delta)}$\end{center} relative to $(\g';\f')$ such that: \begin{center}$\judge{\delta:\Delta,\underline{x}:T(\delta)}{\psi(\g(\delta),\q(\delta)^*\underline{x})=\alpha(\delta)^*\psi'(\delta,\underline{x})}$\\and equivalently $\judge{\delta:\Delta,x':S(\g'(\delta))}{\phi(\g(\delta),\alpha(\delta)^*x')=\q(\delta)^*\phi'(\delta,x')}$.\end{center}

\proof By induction on the complexity of $\equivpair{\phi}{\psi}$.

\medskip

If $\equivpair{\phi}{\psi}$ is of the form \texttt{(a)} then we are done.

\medskip

If $\equivpair{\phi}{\psi}$ is of the form \texttt{(b)} then $S(\gamma)\equiv \Pi_{x:A(\gamma)}B(\gamma,x)$ and $T(\delta)\equiv \Pi_{\underline{x}:A'(\delta)}B'(\delta,\underline{x})$ for some judgements: \begin{itemize}
        \item $\judge{\gamma:\Gamma}{A(\gamma):\type}$ \, and \, $\judge{\gamma:\Gamma,x:A(\gamma)}{B(\gamma,x):\type}$
        \item $\judge{\delta:\Delta}{A'(\delta):\type}$ \, and \, $\judge{\delta:\Delta,\underline{x}:A'(\delta)}{B'(\delta,\underline{x}):\type}$
\end{itemize} and there are a canonical homotopy equivalence $\equivpair{f_1}{g_1}$: \begin{center}$\judge{\gamma:\Gamma,x:A(\gamma)}{f_1(\gamma,x):A'(\f(\gamma))}$\\ $\judge{\gamma:\Gamma,\underline{x}':A'(\f(\gamma))}{g_1(\gamma,\underline{x}'):A(\gamma)}$\end{center} between $A(\gamma)$ and $A'(\delta)$ relative to $(\f;\g)$ and a canonical homotopy equivalence $\equivpair{f_2}{g_2}$: \begin{center}$\judge{\gamma:\Gamma,x:A(\gamma),y:B(\gamma,x)}{f_2(\gamma,x,y):B'(\f(\gamma),f_1(\gamma,x))}$ \\ $\judge{\gamma:\Gamma,x:A(\gamma),\underline{y}':B'(\f(\gamma),f_1(\gamma,x))}{g_2(\gamma,x,\underline{y}'):B(\gamma,x)}$\end{center} between $B(\gamma,x)$ and $B'(\delta,\underline{x})$ relative to the extension $(\f,\f_{m+1};\g,\g_{n+1})$ of $(\f;\g)$ via $(f_1;g_1)$ (\autoref{context homotopy equivalence expansion}), in such a way that the homotopy equivalence $( f^\Pi;g^\Pi)$ of \autoref{pi equiv} is the given $\equivpair{\phi}{\psi}$. By inductive hypothesis there are a canonical homotopy equivalence $\equivpair{g'_1}{f'_1}$: \begin{center}$\judge{\delta:\Delta,\underline{x}:A'(\delta)}{g'_1(\delta,\underline{x}):A(\g'(\delta))}$\\ $\judge{\delta:\Delta,x':A(\g'(\delta))}{f'_1(\delta,x'):A'(\delta)}$\end{center} between $A'(\delta)$ and $A(\gamma)$ relative to $(\g';\f')$ and a canonical homotopy equivalence $\equivpair{g'_2}{f'_2}$: \begin{center}$\judge{\delta:\Delta,\underline{x}:A'(\delta),\underline{y}:B'(\delta,x)}{g'_2(\delta,\underline{x},\underline{y}):B(\g'(\delta),g'_1(\delta,\underline{x}))}$ \\ $\judge{\delta:\Delta,\underline{x}:A'(\delta),y':B(\g'(\delta),g'_1(\delta,\underline{x}))}{f'_2(\delta,\underline{x},y'):B'(\delta,\underline{x})}$\end{center} between $B'(\delta,\underline{x})$ and $B(\gamma,x)$ relative to the extension $(\g',\g'_{n+1};\f',\f'_{m+1})$ of $(\g';\f')$ via $(g'_1;f'_1)$, in such a way that: \begin{center}$\judge{\delta:\Delta,\underline{x}:A'(\delta)}{\alpha_1(\delta,\underline{x}):g_1(\g(\delta),\q(\delta)^*\underline{x})=\alpha(\delta)^*g'_1(\delta,\underline{x})}$\\i.e. $\judge{\delta:\Delta,x':A(\g'(\delta))}{f_1(\g(\delta),\alpha(\delta)^*x')=\q(\delta)^*f'_1(\delta,x')}$\end{center} and that: \begin{center}$\judge{\delta:\Delta,\underline{x}:A'(\delta),\underline{y}:B'(\delta,\underline{x})}{g_2(\g(\delta),\g_{n+1}(\delta,\underline{x}),(\q(\delta),\q_{n+1}(\delta,\underline{x}))^*\underline{y})=(\alpha(\delta),\alpha_1(\delta,\underline{x}))^*g'_1(\delta,\underline{x},\underline{y})}$\\i.e. $\judge{\delta:\Delta,\underline{x}:A'(\delta),y':B(\g'(\delta),g_1'(\delta,\underline{x}))}{f_2(\g(\delta),\g_{n+1}(\delta,\underline{x}),(\alpha(\delta),\alpha_1(\delta,\underline{x}))^*y')=(\q(\delta),\q_{n+1}(\delta,\underline{x}))^*f'_2(\delta,\underline{x},y')}$\end{center} where we remind that: \[\begin{aligned}
    \g_{n+1}(\delta,\underline{x})&\equiv g_1(\g(\delta),\q(\delta)^*\underline{x})\\
    \q_{n+1}(\delta,\underline{x})&\equiv q_1(\g(\delta),\q(\delta)^*\underline{x})^{-1}
\end{aligned}\] in context $\delta,\underline{x}$, being $\judge{\gamma:\Gamma,\underline{x}':A'(\f(\gamma))}{q_1(\gamma,\underline{x}'):\underline{x}'=f_1(\gamma,g_1(\gamma,\underline{x}'))}$. Then $\equivpair{g'^\Pi}{f'^\Pi}$ is a canonical homotopy equivalence between $\Pi_{\underline{x}:A'(\delta)}B'(\delta,\underline{x})$ and $\Pi_{x:A(\gamma)}B(\gamma,x)$ and relative to $(\g';\f')$, by \texttt{(b)}. We are left to verify that: $$\judge{\delta,z':\Pi_{x':A(g'(\delta))}B(\g'(\delta),x')}{f^\Pi(\g(\delta),\alpha(\delta)^*z')=\q(\delta)^*f'^\Pi(\delta,z')}.$$ Let us observe that $f'^\Pi(\delta,z')\equiv$ \[\begin{aligned}
&\equiv \lambda \underline{x}:A'(\delta).f_2'(\delta,\underline{x},\ev(z',g_1'(\delta,\underline{x})))\\
&=\lambda \underline{x}.(\q(\delta)^{-1})^*q_1(\g(\delta),\q(\delta)^*\underline{x})^*f_2(\g(\delta),\g_{n+1}(\delta,\underline{x}),(\alpha(\delta),\alpha_1(\delta,\underline{x}))^*\ev(z',g_1'(\delta,\underline{x})))\\
&=\lambda \underline{x}.(\q(\delta)^{-1})^*q_1(\g(\delta),\q(\delta)^*\underline{x})^*f_2(\g(\delta),\g_{n+1}(\delta,\underline{x}),\ev(\alpha(\delta)^*z',g_1(\g(\delta),q(\delta)^*\underline{x})))
\end{aligned}\] in context $\delta,z'$, where the first propositional equality follows by: \begin{center}$\judge{\delta:\Delta,\underline{x}:A'(\delta),y':B(\g'(\delta),g_1'(\delta,\underline{x}))}$\\
${f_2(\g(\delta),\g_{n+1}(\delta,\underline{x}),(\alpha(\delta),\alpha_1(\delta,\underline{x}))^*y')=(\q(\delta),\q_{n+1}(\delta,\underline{x}))^*f'_2(\delta,\underline{x},y')}$\end{center} and the second by: $$\judge{\delta:\Delta,\underline{x}:A'(\delta)}{\alpha_1(\delta,\underline{x}):g_1(\g(\delta),\q(\delta)^*\underline{x})=\alpha(\delta)^*g'_1(\delta,\underline{x})}.$$ Therefore $q(\delta)^*f'^\Pi(\delta,z')=$ $$\begin{aligned}\text{ }&=q(\delta)^*\lambda \underline{x}.(\q(\delta)^{-1})^*q_1(\g(\delta),\q(\delta)^*\underline{x})^*f_2(\g(\delta),\g_{n+1}(\delta,\underline{x}),\ev(\alpha(\delta)^*z',\g_{n+1}(\delta,\underline{x})))\\&=\lambda \underline{x}'':A'(\f(\g(\delta))).q_1(\g(\delta),\underline{x}'')^*f_2(\g(\delta),g_1(\g(\delta),\underline{x}''),\ev(\alpha(\delta)^*z',g_1(\g(\delta),\underline{x}'')))\\&\equiv f^\Pi(\g(\delta),\alpha(\delta)^*z')\end{aligned}$$ in context ${\delta,z'}$, where the first propositional equality follows by propositional function extensionality and propositional functoriality and the second by multiple (generalised) path induction on $q(\delta)$. We are done.

\medskip

If $\equivpair{\phi}{\psi}$ is of the form \texttt{(c)} then as before $S(\gamma)\equiv \Sigma_{x:A(\gamma)}B(\gamma,x)$ and $T(\delta)\equiv \Sigma_{\underline{x}:A'(\delta)}B'(\delta,\underline{x})$ for some judgements: \begin{itemize}
        \item $\judge{\gamma:\Gamma}{A(\gamma):\type}$ \, and \, $\judge{\gamma:\Gamma,x:A(\gamma)}{B(\gamma,x):\type}$
        \item $\judge{\delta:\Delta}{A'(\delta):\type}$ \, and \, $\judge{\delta:\Delta,\underline{x}:A'(\delta)}{B'(\delta,\underline{x}):\type}$
\end{itemize} and there are a canonical homotopy equivalence $\equivpair{f_1}{g_1}$: \begin{center}$\judge{\gamma:\Gamma,x:A(\gamma)}{f_1(\gamma,x):A'(\f(\gamma))}$\\ $\judge{\gamma:\Gamma,\underline{x}':A'(\f(\gamma))}{g_1(\gamma,\underline{x}'):A(\gamma)}$\end{center} between $A(\gamma)$ and $A'(\delta)$ relative to $(\f;\g)$ and a canonical homotopy equivalence $\equivpair{f_2}{g_2}$: \begin{center}$\judge{\gamma:\Gamma,x:A(\gamma),y:B(\gamma,x)}{f_2(\gamma,x,y):B'(\f(\gamma),f_1(\gamma,x))}$ \\ $\judge{\gamma:\Gamma,x:A(\gamma),\underline{y}':B'(\f(\gamma),f_1(\gamma,x))}{g_2(\gamma,x,\underline{y}'):B(\gamma,x)}$\end{center} between $B(\gamma,x)$ and $B'(\delta,\underline{x})$ relative to the extension $(\f,\f_{m+1};\g,\g_{n+1})$ of $(\f;\g)$ via $(f_1;g_1)$ (\autoref{context homotopy equivalence expansion}), in such a way that the homotopy equivalence $( f^\Sigma;g^\Sigma)$ of \autoref{sigma equiv} is the given $\equivpair{\phi}{\psi}$. Again, by inductive hypothesis there are a canonical homotopy equivalence $\equivpair{g'_1}{f'_1}$: \begin{center}$\judge{\delta:\Delta,\underline{x}:A'(\delta)}{g'_1(\delta,\underline{x}):A(\g'(\delta))}$\\ $\judge{\delta:\Delta,x':A(\g'(\delta))}{f'_1(\delta,x'):A'(\delta)}$\end{center} between $A'(\delta)$ and $A(\gamma)$ relative to $(\g';\f')$ and a canonical homotopy equivalence $\equivpair{g'_2}{f'_2}$: \begin{center}$\judge{\delta:\Delta,\underline{x}:A'(\delta),\underline{y}:B'(\delta,x)}{g'_2(\delta,\underline{x},\underline{y}):B(\g'(\delta),g'_1(\delta,\underline{x}))}$ \\ $\judge{\delta:\Delta,\underline{x}:A'(\delta),y':B(\g'(\delta),g'_1(\delta,\underline{x}))}{f'_2(\delta,\underline{x},y'):B'(\delta,\underline{x})}$\end{center} between $B'(\delta,\underline{x})$ and $B(\gamma,x)$ relative to the extension $(\g',\g'_{n+1};\f',\f'_{m+1})$ of $(\g';\f')$ via $(g'_1;f'_1)$, in such a way that: \begin{center}$\judge{\delta:\Delta,\underline{x}:A'(\delta)}{\alpha_1(\delta,\underline{x}):g_1(\g(\delta),\q(\delta)^*\underline{x})=\alpha(\delta)^*g'_1(\delta,\underline{x})}$\\i.e. $\judge{\delta:\Delta,x':A(\g'(\delta))}{f_1(\g(\delta),\alpha(\delta)^*x')=\q(\delta)^*f'_1(\delta,x')}$\end{center} and that: \begin{center}$\judge{\delta:\Delta,\underline{x}:A'(\delta),\underline{y}:B'(\delta,\underline{x})}{g_2(\g(\delta),\g_{n+1}(\delta,\underline{x}),(\q(\delta),\q_{n+1}(\delta,\underline{x}))^*\underline{y})=(\alpha(\delta),\alpha_1(\delta,\underline{x}))^*g'_1(\delta,\underline{x},\underline{y})}$\\i.e. $\judge{\delta:\Delta,\underline{x}:A'(\delta),y':B(\g'(\delta),g_1'(\delta,\underline{x}))}{f_2(\g(\delta),\g_{n+1}(\delta,\underline{x}),(\alpha(\delta),\alpha_1(\delta,\underline{x}))^*y')=(\q(\delta),\q_{n+1}(\delta,\underline{x}))^*f'_2(\delta,\underline{x},y')}$\end{center} where, as before, we remind that: \[\begin{aligned}
    \g_{n+1}(\delta,\underline{x})&\equiv g_1(\g(\delta),\q(\delta)^*\underline{x})\\
    \q_{n+1}(\delta,\underline{x})&\equiv q_1(\g(\delta),\q(\delta)^*\underline{x})^{-1}
\end{aligned}\] in context $\delta,\underline{x}$, being $\judge{\gamma:\Gamma,\underline{x}':A'(\f(\gamma))}{q_1(\gamma,\underline{x}'):\underline{x}'=f_1(\gamma,g_1(\gamma,\underline{x}'))}$. Then $\equivpair{g'^\Sigma}{f'^\Sigma}$ is a canonical homotopy equivalence between $\Sigma_{\underline{x}:A'(\delta)}B'(\delta,\underline{x})$ and $\Sigma_{x:A(\gamma)}B(\gamma,x)$ and relative to $(\g';\f')$, by \texttt{(c)}. We are left to verify that: $$\judge{\delta,\underline{u}:\Sigma_{x':A'(\delta)}B'(\delta,x')}{g^\Sigma(\g(\delta),\q(\delta)^*\underline{u})=\alpha(\delta)^*g'^\Sigma(\delta,\underline{u})}.$$ Let us observe that: \[\begin{aligned}
    \judge{\delta,\underline{u}}{\alpha_1(\delta,\pi_1\underline{u}):\g_{n+1}(\delta,\pi_1\underline{u})&=\alpha(\delta)^*g_1'(\delta,\pi_1\underline{u})}\\
    \judge{\delta,\underline{u}}{g_2(\g(\delta),\g_{n+1}(\delta,\pi_1\underline{u}),(\q(\delta),\q_{n+1}(\delta,\pi_1\underline{u}))^*\pi_2\underline{u})&=(\alpha(\delta),\alpha_1(\delta,\pi_1\underline{u}))^*g_2'(\delta,\pi_1\underline{u},\pi_2\underline{u})}\\
    &=\alpha_1(\delta,\pi_1\underline{u})^*\alpha(\delta)^*g_2'(\delta,\pi_1\underline{u},\pi_2\underline{u})
\end{aligned}\] where the first equality in the second judgement follows by: \begin{center}
    $\judgectx{\delta:\Delta,\underline{x}:A'(\delta),\underline{y}:B'(\delta,\underline{x})}$\\
    ${g_2(\g(\delta),\g_{n+1}(\delta,\underline{x}),(\q(\delta),\q_{n+1}(\delta,\underline{x}))^*\underline{y})=(\alpha(\delta),\alpha_1(\delta,\underline{x}))^*g'_1(\delta,\underline{x},\underline{y})}$
\end{center} and the second by (generalised) path induction on $\alpha_1(\delta,\pi_1\underline{u})$. Therefore: $$\begin{aligned} \langle \g_{n+1}(\delta,\pi_1\underline{u}),& g_2(\g(\delta),\g_{n+1}(\delta,\pi_1\underline{u}),(\q(\delta),\q_{n+1}(\delta,\pi_1\underline{u}))^*\pi_2\underline{u}) \rangle\\=&\langle \alpha(\delta)^*g_1'(\delta,\pi_1\underline{u}), \alpha(\delta)^*g_2'(\delta,\pi_1\underline{u},\pi_2\underline{u})\rangle\\=&\alpha(\delta)^*\langle g_1'(\delta,\pi_1\underline{u}), g_2'(\delta,\pi_1\underline{u},\pi_2\underline{u})\rangle\\\equiv& \alpha(\delta)^*g'^\Sigma(\delta,\underline{u}) \end{aligned}$$ in context ${\delta,\underline{u}}$, where the second propositional equality follows by multiple (generalised) path induction on $\alpha(\delta)$. Finally $g^\Sigma(\g(\delta),\q(\delta)^*\underline{u})\equiv$ $$\begin{aligned}\text{ }&\equiv\langle g_1(\g(\delta),\pi_1\q(\delta)^*\underline{u}), g_2(\g(\delta),g_1(\g(\delta),\pi_1\q(\delta)^*\underline{u}),q_1(\g(\delta),\pi_1\q(\delta)^*\underline{u})^*\pi_2\q(\delta)^*\underline{u}) \rangle \\&=\langle \g_{n+1}(\delta,\pi_1\underline{u}), g_2(\g(\delta),\g_{n+1}(\delta,\pi_1\underline{u}),(\q(\delta),\q_{n+1}(\delta,\pi_1\underline{u}))^*\pi_2\underline{u}) \rangle \end{aligned}$$ in context ${\delta,\underline{u}}$, where the propositional equality follows by multiple (generalised) path induction on $\q(\delta)$ followed by (generalised) path induction on $q_1(\g(\delta),\pi_1\q(\delta)^*\underline{u})$. We are done.

\medskip

\textit{Notation for the last section of the proof.} Whenever: \begin{center}
$\judge{\omega:\Omega}{x(\omega),y(\omega):O(\omega)}$\\ $\judge{\omega}{\alpha:x(\omega)=y(\omega)}$\\ $\judge{\omega,\omega':\Omega}{p:\omega'=\omega}$ \end{center} we denote as $p^*\alpha:x(\omega')=y(\omega')$ the usual transport of the term $\alpha(\omega)$ from the type $x(\omega)=y(\omega)$ to the type $x(\omega')=y(\omega')$. However, the operation $p^*$ is defined on the terms of $O(\omega)$ as well, hence the terms $p^*x(\omega),p^*y(\omega):O(\omega')$ are defined. Therefore the operation $p^*$ might be extended as usual to a propositional functor from $O(\omega)$ to $O(\omega')$ and a term $p^*\alpha: p^*x(\omega)=p^*y(\omega)$ is defined. In the last section, in order to distinguish between the terms: \begin{center} $p^*\alpha:x(\omega')=y(\omega')$ and $p^*\alpha:p^*x(\omega)=p^*y(\omega)$ \end{center} that we usually indicate by the same notation, we indicate the \textit{latter} as $\underline{p^*\alpha}$. Observe that: \begin{center}
    $\judge{p:\omega'=\omega}{}$ $\begin{tikzcd}x(\omega') \arrow[d, Rightarrow, "x(p)"'] \arrow[r, Rightarrow, "p^*\alpha"] & y(\omega') \arrow[d, Rightarrow, "y(p)"] \\ p^*x(\omega) \arrow[r, Rightarrow, "\underline{p^*\alpha}"'] & p^*y(\omega). \arrow[shorten <=14pt, shorten >=14pt, Rightarrow, no head, from=2-1, to=1-2] \end{tikzcd}$
\end{center} In order to verify this, since $\refl(\omega')^*\alpha=\alpha$ and since $x(\refl(\omega'))$ and $y(\refl(\omega'))$ are propositionally equal to the canonical proofs $\xi(x(\omega'))$ and $\xi(y(\omega'))$ that $x(\omega')=\refl(\omega')^*x(\omega')$ and $y(\omega')=\refl(\omega')^*y(\omega')$ and by path induction, it is enough to verify that: \begin{center}
    $\judge{\omega}{}$ $\begin{tikzcd}x(\omega') \arrow[d, Rightarrow, "\xi(x(\omega'))"'] \arrow[r, Rightarrow, "\alpha"] & y(\omega') \arrow[d, Rightarrow, "\xi(y(\omega'))"] \\ \refl(\omega')^*x(\omega') \arrow[r, Rightarrow, "\underline{\refl(\omega')^*\alpha}"'] & \refl(\omega')^*y(\omega'). \arrow[shorten <=21pt, shorten >=21pt, Rightarrow, no head, from=2-1, to=1-2] \end{tikzcd}$
\end{center} which is true because $\xi(-)$ is a homotopy.

\medskip

If $\equivpair{\phi}{\psi}$ is of the form \texttt{(d)} then $S(\gamma)\equiv s_1(\gamma)=s_2(\gamma)$ and $T(\delta)\equiv t_1(\delta)=t_2(\delta)$ for some judgements: \begin{center}
$\judge{\gamma:\Gamma}{A(\gamma):\type}$ \, and \, $\judge{\delta:\Delta}{A'(\delta):\type}$
\end{center} and some judgements:
\begin{center} $\judge{\gamma:\Gamma}{s_1(\gamma),s_2(\gamma):A(\gamma)}$ \, and \, $\judge{\delta:\Delta}{t_1(\delta),t_2(\delta):A'(\delta)}$\end{center} and there are a canonical homotopy equivalence $\equivpair{f}{g}$: \begin{center}$\judge{\gamma:\Gamma,x:A(\gamma)}{f(\gamma,x):A'(\f(\gamma))}$\\ $\judge{\gamma:\Gamma,\underline{x}':A'(\f(\gamma))}{g(\gamma,\underline{x}'):A(\gamma)}$\end{center} between $A(\gamma)$ and $A'(\delta)$ relative to $(\f;\g)$ and judgements: \begin{center}$\judge{\gamma:\Gamma}{r_1(\gamma):f(\gamma,s_1(\gamma))=t_1(\f(\gamma))}$\\ $\judge{\gamma:\Gamma}{r_2(\gamma):f(\gamma,s_2(\gamma))=t_2(\f(\gamma)).}$\end{center} in such a way that the homotopy equivalence $(f^=,g^=)$ of \autoref{id equiv} is the given $\equivpair{\phi}{\psi}$. By inductive hypothesis there is a canonical homotopy equivalence $\equivpair{g'}{f'}$: \begin{center}$\judge{\delta:\Delta,\underline{x}:A'(\delta)}{g'(\delta,\underline{x}):A(\g'(\delta))}$\\ $\judge{\delta:\Delta,x':A(\g'(\delta))}{f'(\delta,x'):A'(\delta)}$\end{center} between $A'(\delta)$ and $A(\gamma)$ relative to $(\g';\f')$ such that: \begin{center}$\judge{\delta:\Delta,\underline{x}:A'(\delta)}{\alpha_1(\delta,\underline{x}):g(\g(\delta),\q(\delta)^*\underline{x})=\alpha(\delta)^*g'(\delta,\underline{x})}$\\i.e. $\judge{\delta:\Delta,x':A(\g'(\delta))}{f(\g(\delta),\alpha(\delta)^*x')=\q(\delta)^*f'(\delta,x')}$.\end{center} For $i \in \{1,2\}$, let us observe that in context $\delta$: \begin{center}
    $\alpha_1(\delta,t_i(\delta))^{-1}:\alpha(\delta)^*g'(\delta,t_i(\delta))=g(\g(\delta),\q(\delta)^*t_i(\delta))$\\
    $g(\g(\delta),t_i(\q(\delta))^{-1}):g(\g(\delta),\q(\delta)^*t_i(\delta))=g(\g(\delta),t_i(\f(\g(\delta))))$\\
    $g(\g(\delta),r_i(\g(\delta))^{-1}):g(\g(\delta),t_i(\f(\g(\delta))))=g(\g(\delta),f(\g(\delta),s_i(\g(\delta))))$\\
    $p(\g(\delta),s_i(\g(\delta))):g(\g(\delta),f(\g(\delta),s_i(\g(\delta))))=s_i(\g(\delta))$\\
    $s_i(\alpha(\gamma)):s_i(\g(\delta))=\alpha(\delta)^*s_i(\g'(\delta))$
\end{center} and let us call $a_i(\delta):\alpha(\delta)^*g'(\delta,t_i(\delta))=\alpha(\delta)^*s_i(\g'(\delta))$. Let $r_i'(\delta):g'(\delta,t_i(\delta))=s_i(\g'(\delta))$ be such that: $$\underline{\alpha(\delta)^*r_i'(\delta)}=a_i(\delta).$$ Hence $\equivpair{g'^=}{f'^=}$ is a canonical homotopy equivalence between $t_1(\delta)=t_2(\delta)$ and $s_1(\gamma)=s_2(\gamma)$ relative to $(\g';\f')$, by \texttt{(d)}, and we are left to verify that: $$\judge{\delta:\Delta,\underline{p}:t_1(\delta)=t_2(\delta)}{g^=(\g(\delta),\q(\delta)^*\underline{p})=\alpha(\delta)^*g'^=(\delta,\underline{p})}$$ that is in fact true, as $\judge{\delta,\underline{p}}{\alpha(\delta)^*g'^=(\delta,\underline{p})\equiv}$ $$\begin{aligned}\;&\equiv\alpha(\delta)^*(r_1'(\delta)^{-1}\bullet  g'(\delta,\underline{p})  \bullet r_2'(\delta))\\
&=s_1(\alpha(\delta))\bullet\underline{\alpha(\delta)^*(r_1'(\delta)^{-1}\bullet  g'(\delta,\underline{p})  \bullet r_2'(\delta))}\bullet s_2(\alpha(\delta))^{-1}\\
&=s_1(\alpha(\delta))\bullet\underline{\alpha(\delta)^*r_1'(\delta)^{-1}}\bullet  \underline{\alpha(\delta)^*g'(\delta,\underline{p})}  \bullet \underline{\alpha(\delta)^*r_2'(\delta)}\bullet s_2(\alpha(\delta))^{-1}\\
&=s_1(\alpha(\delta))\bullet a_1(\gamma)^{-1} \bullet  \underline{\alpha(\delta)^*g'(\delta,\underline{p})}  \bullet a_2(\gamma) \bullet s_2(\alpha(\delta))^{-1}\\
&=p(\g(\delta),s_1(\g(\delta)))^{-1}\bullet g(\g(\delta),r_1(\g(\delta)))\bullet g(\g(\delta),t_1(\q(\delta)))\bullet\alpha_1(\delta,t_1(\delta))\bullet \underline{\alpha(\delta)^*g'(\delta,\underline{p})}\;\bullet\\
&\;\;\;\;\;\alpha_1(\delta,t_2(\delta))^{-1}\bullet g(\g(\delta),t_2(\q(\delta))^{-1}) \bullet g(\g(\delta),r_2(\g(\delta))^{-1})\bullet p(\g(\delta),s_2(\g(\delta)))\\
&=p(\g(\delta),s_1(\g(\delta)))^{-1}\bullet g(\g(\delta),r_1(\g(\delta)))\bullet g(\g(\delta),t_1(\q(\delta)))\bullet g(\g(\delta),\underline{\q(\delta)^*\underline{p}})\;\bullet\\
&\;\;\;\;\; g(\g(\delta),t_2(\q(\delta))^{-1}) \bullet g(\g(\delta),r_2(\g(\delta))^{-1})\bullet p(\g(\delta),s_2(\g(\delta)))\\
&=p(\g(\delta),s_1(\g(\delta)))^{-1}\bullet g(\g(\delta),r_1(\g(\delta)))\bullet g(\g(\delta),t_1(\q(\delta))\bullet \underline{\q(\delta)^*\underline{p}}\bullet t_2(\q(\delta))^{-1})\;\bullet\\
&\;\;\;\;\; g(\g(\delta),r_2(\g(\delta))^{-1})\bullet p(\g(\delta),s_2(\g(\delta)))\\
&=p(\g(\delta),s_1(\g(\delta)))^{-1}\bullet g(\g(\delta),r_1(\g(\delta)))\bullet g(\g(\delta),\q(\delta)^*\underline{p})\;\bullet\\
&\;\;\;\;\; g(\g(\delta),r_2(\g(\delta))^{-1})\bullet p(\g(\delta),s_2(\g(\delta)))\\
&= p(\g(\delta),s_1(\g(\delta)))^{-1}\bullet g(\g(\delta),r_1(\g(\delta))\bullet\q(\delta)^*\underline{p}\bullet r_2(\g(\delta))^{-1})\bullet p(\g(\delta),s_2(\g(\delta)))\\
&\equiv g^=(\g(\delta),\q(\delta)^*\underline{p})
\end{aligned}$$ hence we are done.\endproof
\end{lem}

\begin{lem}\label{lemma x}
Let $\gamma:\Gamma$, $\delta:\Delta$ and $\omega:\Omega$ be contexts, let $(\f;\g)$ be a canonical context homotopy equivalence between $\gamma$ and $\delta$ and let $(\f',\g')$ be a context homotopy equivalence between $\delta$ and $\omega$.

If we are given judgements $\judge{\gamma}{S(\gamma):\type}$, $\judge{\delta}{T(\delta):\type}$ and $\judge{\omega}{U(\omega):\type}$ and canonical homotopy equivalences $\equivpair{\phi}{\psi}$ between $S(\gamma)$ and $T(\delta)$ relative to $(\f;\g)$ and $\equivpair{\phi'}{\psi'}$ between $T(\delta)$ and $U(\omega)$ relative to $(\f';\g')$, then there exists a canonical homotopy equivalence $\equivpair{{\phi}''}{{\psi}''}$ between $S(\gamma)$ and $U(\omega)$ relative to $(\f'\f;\g\g')$ and such that: \begin{center} $\judge{\gamma,s:S(\gamma)}{{\phi}''(\gamma,s)=\phi'(\f(\gamma),\phi(\gamma,s))}$\\ and equivalently $\judge{\gamma,u'':U(\f'(\f(\gamma)))}{{\psi}''(\gamma,u'')=\psi(\gamma,\psi'(\f(\gamma),u''))}$. \end{center}

\proof
By induction on the complexity of $\equivpair{\phi}{\psi}$.
\endproof
\end{lem}

\begin{lem}\label{lemma y}
Let $\gamma:\Gamma$ and $\delta:\Delta$ be contexts and let $(\f;\g)$ and $(\f';\g')$ be context homotopy equivalences between $\gamma$ and $\delta$. Let us assume that $\judge{\gamma}{\alpha(\gamma):\f(\gamma)=\f'(\gamma)}$ and equivalently $\judge{\delta}{\g(\delta)=\g'(\delta)}$.

If we are given judgements $\judge{\gamma}{S(\gamma):\type}$ and $\judge{\delta}{T(\delta):\type}$ and a canonical homotopy equivalence $\equivpair{\phi}{\psi}$ between $S(\gamma)$ and $T(\delta)$ relative to $(\f;\g)$, then there is a canonical homotopy equivalence $\equivpair{\phi'}{\psi'}$ between $S(\gamma)$ and $T(\delta)$ relative to $(\f';\g')$ such that: \begin{center} $\judge{\gamma,s:S(\gamma)}{\phi(\gamma,s)=\alpha(\gamma)^*\phi'(\gamma,s)}$\\ and equivalently $\judge{\gamma,t':T(\f'(\gamma))}{{\psi}(\gamma,\alpha(\gamma)^*t')=\psi'(\gamma,t')}$. \end{center}

\proof
By induction on the complexity of $\equivpair{\phi}{\psi}$.
\endproof
\end{lem}

\begin{prop}[Transitivity]\label{transitivity}
Let $\gamma:\Gamma$, $\delta:\Delta$ and $\omega:\Omega$ be contexts. If $\equivpair{\ci}{\di}$ is a canonical context homotopy equivalence between $\gamma$ and $\delta$ and $\equivpair{\ci'}{\di'}$ is a canonical context homotopy equivalence between $\delta$ and $\omega$ then there is a canonical context homotopy equivalence $\equivpair{\ci''}{\di''}$ between $\gamma$ and $\omega$ such that $\judge{\gamma}{\ci''(\gamma)=\ci'(\ci(\gamma))}$ and equivalently $\judge{\delta}{\di''(\delta)=\di(\di'(\delta))}$.

\proof
By induction on the complexity of $\equivpair{\ci}{\di}$. If $\equivpair{\ci}{\di}$ is of the form \texttt{(a)} then $\equivpair{\ci'}{\di'}$ is of the form \texttt{(a)} as well and we are done. Let $\equivpair{\ci}{\di}$ be of the form \texttt{(b)}, hence $\equivpair{\ci'}{\di'}$ is of the form \texttt{(b)} as well. Then $\gamma  : \Gamma$, $\delta:\Delta$ and $\omega:\Omega$ are of the form $\gamma'  : \Gamma',x : A(\gamma')$ and $\delta' :  \Delta',\underline{x}: A'(\delta')$ and $\omega':\Omega',\underline{\underline{x}}:A''(\Omega')$ respectively. Moreover there are a canonical context homotopy equivalence: \begin{center}$\judge{\gamma':\Gamma'}{\f(\gamma'):\Delta'}$\\ $\judge{\delta':\Delta'}{\g(\delta'):\Gamma'}$\end{center} and a canonical homotopy equivalence: \begin{center} $\judge{\gamma':\Gamma',x:A(\gamma')}{f(\gamma',x):A'(\f(\gamma'))}$\\ $\judge{\gamma':\Gamma',\underline{x}':A'(\f(\gamma'))}{g(\gamma',\underline{x}'):A(\gamma')}$\end{center} between $A(\gamma')$ and $A'(\delta')$ relative to $\equivpair{\f}{\g}$ in such a way that the extension: \begin{center} $\judge{\gamma' : \Gamma',  x:A(\gamma')}{(\f,\f_{m+1})(\gamma',x):\Delta',A'}$\\ $\judge{\delta' : \Delta', \underline{x}:A'(\delta')}{(\g,\g_{n+1})(\delta',\underline{x}):\Gamma',A}$\end{center} of $\equivpair{\f}{\g}$ via $\equivpair{f}{g}$ (\autoref{context homotopy equivalence expansion}) is the given $\equivpair{\ci}{\di}$. Analogously, there are a canonical context homotopy equivalence: \begin{center}$\judge{\delta':\Delta'}{\f'(\delta'):\Omega'}$\\ $\judge{\omega':\Omega'}{\g(\omega'):\Delta'}$\end{center} and a canonical homotopy equivalence: \begin{center} $\judge{\delta':\Delta',\underline{x}:A'(\delta')}{f'(\delta',\underline{x}):A''(\f'(\delta'))}$\\ $\judge{\delta':\Delta',\underline{\underline{x}}:A''(\f'(\delta'))}{g'(\delta',\underline{\underline{x}}'):A'(\delta')}$\end{center} between $A'(\delta')$ and $A''(\omega')$ relative to $\equivpair{\f'}{\g'}$ in such a way that the extension: \begin{center} $\judge{\delta':\Delta',\underline{x}:A'(\delta')}{(\f',\f'_{m+1})(\delta',\underline{x}):\Omega',A''}$\\ $\judge{\omega' : \Omega', \underline{\underline{x}}:A''(\omega')}{(\g',\g'_{n+1})(\omega',\underline{\underline{x}}):\Delta',A'}$\end{center} of $\equivpair{\f'}{\g'}$ via $\equivpair{f'}{g'}$ (\autoref{context homotopy equivalence expansion}) is the given $\equivpair{\ci'}{\di'}$. By inductive hypothesis there is a canonical context homotopy equivalence $\equivpair{\f''}{\g''}$ between $\gamma'$ and $\omega'$ such that: \begin{center} $\judge{\gamma'}{\f''(\gamma')=\f'(\f(\gamma'))}$ and $\judge{\omega'}{\g''(\omega')=\g(\g'(\omega'))}$ \end{center} Let us consider the canonical homotopy equivalence: \begin{center} $\judge{\gamma':\Gamma',x:A(\gamma')}{\tilde{f}(\gamma',x)\equiv f'(\f(\gamma'),f(\gamma',x)):A''(\f'(\f(\gamma')))}$\\ $\judge{\gamma':\Gamma',\underline{\underline{x}}'':A''(\f'(\f(\gamma')))}{\tilde{g}(\gamma',\underline{\underline{x}}'')\equiv g(\gamma',g'(\f(\gamma'),\underline{\underline{x}}'')):A(\gamma')}$\end{center} relative to $(\f'\f;\g\g')$ and let us observe that: $$\begin{aligned}
    (\f'(\f(\gamma)),f'(\f(\gamma'),f(\gamma',x)))&\equiv (\f'(\f(\gamma)),\f_{m+1}'(\f(\gamma'),\f_{m+1}(\gamma',x)))\\&\equiv (\f',\f'_{m+1})(\f,\f_{m+1})(\gamma',x)\\&\equiv \ci'(\ci(\gamma))
\end{aligned}$$ hence the extension of $(\f'\f;\g\g')$ via $(\tilde{f};\tilde{g})$ is pairwise propositionally equal to $(\ci'\ci;\di\di')$. By \autoref{lemma x}, there is canonical homotopy equivalence $\equivpair{\overline{f}}{\overline{g}}$ between $A(\gamma')$ and $A''(\omega')$ relative to $(\f'\f;\g\g')$ such that $\judge{\gamma'}{\overline{f}(\gamma')=\tilde{f}(\gamma')}$ and equivalently $\judge{\omega'}{\overline{g}(\omega')=\tilde{g}(\omega')}$. Then, again, the extension of $(\f'\f;\g\g')$ via $\equivpair{\overline{f}}{\overline{g}}$ is pairwise propositionally equal to $(\ci'\ci;\di\di')$. Since $\judge{\gamma'}{\alpha(\gamma'):\f'(\f(\gamma'))=\f''(\gamma')}$ and equivalently $\judge{\omega'}{\g(\g'(\omega'))=\g''(\omega')}$, by \autoref{lemma y} there is a canonical homotopy equivalence $\equivpair{f''}{g''}$ between $A(\gamma')$ and $A''(\omega')$ relative to $\equivpair{\f''}{\g''}$ and such that: \begin{center} $\judge{\gamma',x:A(\gamma')}{\overline{f}(\gamma',x)=\alpha(\gamma')^*f''(\gamma',x)}$\\ and equivalently $\judge{\gamma',\underline{\underline{x}}':A''(\f''(\gamma'))}{\overline{g}(\gamma',\alpha(\gamma')^*\underline{\underline{x}}')=g''(\gamma',\underline{\underline{x}}')}.$ \end{center} Therefore the extension $\equivpair{\f'',\f''_{m+1}}{\g'',\g''_{n+1}}$ of $\equivpair{\f''}{\g''}$ via $\equivpair{f''}{g''}$ is a canonical context homotopy equivalence by \texttt{(b)} and $\equivpair{\f'',\f''_{m+1}}{\g'',\g''_{n+1}}$ is pairwise propositionally equal to the extension of $(\f'\f;\g\g')$ via $\equivpair{\overline{f}}{\overline{g}}$. In particular $\equivpair{\f'',\f''_{m+1}}{\g'',\g''_{n+1}}$ is a canonical context homotopy equivalence between $\gamma$ and $\delta$ that is pairwise propositionally equal to $(\ci'\ci;\di\di')$. \endproof
\end{prop}

We conclude the current subsection with a brief list of technical results that we are using in \autoref{section2.5} and \autoref{section2.6}.

\begin{lem}
\label{lemma z}
Let $\gamma:\Gamma$, $\gamma':\Gamma'$ and $\delta:\Delta$ be contexts, let $( \f;\g)$ be a context homotopy equivalence between $\gamma$ and $\gamma'$ and let $a(\gamma),a'(\gamma'):\Delta$ be context morphisms $\gamma\to\delta$ and $\gamma'\to\delta$ respectively such that: \[\begin{tikzcd}
	\gamma \\
	& \delta
	\\
	{\gamma'}
	\arrow[""{name=0, anchor=center, inner sep=0}, "{a(\gamma)}", bend left, from=1-1, to=2-2]
	\arrow["{a'(\gamma')}"', bend right, from=3-1, to=2-2]
	\arrow["{\f(\gamma)}"', from=1-1, to=3-1]
	\arrow["{\alpha(\gamma)}"{description}, shorten <=15pt, shorten >=15pt, Rightarrow, from=0, to=3-1]
\end{tikzcd}\] Then, whenever $\judge{\delta}{S(\delta):\type}$, there is a canonical homotopy equivalence $\equivpair{\phi}{\psi}$ between $S(a(\gamma))$ and $S(a'(\gamma'))$ relative to $(\f;\g)$ such that:\begin{center} $\judge{\gamma,x:S(a(\gamma))}{\alpha(\gamma)^*\phi(\gamma,x)=x}$ \;\;\; and equivalently \;\;\; $\judge{\gamma,\underline{x}':S(a'(\f(\gamma)))}{\psi(\gamma,\underline{x}')=\alpha(\gamma)^*\underline{x}'}$. \end{center}

\proof
By induction on the complexity of $S(\delta)$.
\endproof
\end{lem}

\begin{cor}\label{semi destiny}
Let $\gamma:\Gamma$ and $\gamma':\Gamma'$ be contexts and let $( \f;\g)$ be a context homotopy equivalence between $\gamma$ and $\gamma'$. Then, whenever $\judge{\gamma'}{S(\gamma'):\type}$, there is a canonical homotopy equivalence $\equivpair{\phi}{\psi}$ between $S(\f(\gamma))$ and $S(\gamma')$ relative to $(\f;\g)$ and such that: $$\judge{\gamma,x:S(\f(\gamma))}{\phi(\gamma,x)=\psi(\gamma,x)=x}.$$

\proof
Follows by \autoref{lemma z}, where $\delta \equiv \gamma'$, $a(\gamma)\equiv \f(\gamma)$, $a'(\gamma')\equiv \gamma'$ and $\alpha(\gamma)\equiv \refl(\f(\gamma))$.
\endproof
\end{cor}

\begin{lem}\label{lemma w}
Let $\gamma:\Gamma$ and $\gamma':\Gamma'$ be contexts and let $( \f;\g)$ be a context homotopy equivalence between $\gamma$ and $\gamma'$. Suppose that $\judge{\gamma}{S(\gamma):\type}$ and $\judge{\gamma'}{T(\gamma'):\type}$ are canonically homotopy equivalent via some $\equivpair{\phi}{\psi}$ relative to $(\f;\g)$. Let $\equivpair{\ci}{\di}$ be some canonical context homotopy equivalence between $\gamma$ and itself such that $\equivpair{\ci}{\di}$ is pairwise propositionally equal to $( \judge{\gamma}{\gamma};\judge{\gamma}{\gamma})$ (whose existence is ensured by \autoref{there are context identities}) and let $\judge{\gamma}{\alpha(\gamma):\ci(\gamma)=\gamma}$.

Then $S(\gamma)$ and $T(\f(\gamma))$ are canonically homotopy equivalent via some $\equivpair{\phi'}{\psi'}$ relative to $\equivpair{\ci}{\di}$ in such a way that: $$\judge{\gamma,x:S(\gamma)}{\phi'(\gamma,x)=\alpha(\gamma)^*\phi(\gamma,x)}$$ and equivalently $\judge{\gamma,\underline{x}':T(\f(\gamma))}{\psi'(\gamma,\alpha(\gamma)^*\underline{x}')=\psi(\gamma,\underline{x}')}$.

\proof
By induction on the complexity of $\equivpair{\phi}{\psi}$.
\endproof
\end{lem}

\begin{lem}\label{small reindexing}
Let $\gamma:\Gamma$ and $\delta:\Delta$ be contexts, let $(\f;\g)$ be a canonical homotopy equivalence between $\delta$ and itself such that $\judge{\delta}{\alpha(\delta):\f(\delta)=\delta}$ and let $(\f';\g')$ be a canonical homotopy equivalence between $\gamma$ and itself such that $\judge{\gamma}{\alpha'(\gamma):\gamma=\f'(\gamma)}$. Moreover, let $\judge{\gamma}{a(\gamma):\Delta}$ be a context morphism and let $\judge{\delta}{S(\delta):\type}$ and $\judge{\delta}{T(\delta):\type}$ be such that there is a canonical   homotopy equivalence $\equivpair{\phi}{\psi}$ between $S(\delta)$ and $T(\delta)$ relative to $(\f;\g)$. Then there is a canonical homotopy equivalence $\equivpair{\phi'}{\psi'}$ between $S(a(\gamma))$ and $T(a(\gamma))$ relative to $(\f';\g')$ and such that: $$\judge{\gamma,x:S(a(\gamma))}{\phi(a(\gamma),x)=\alpha''(\gamma)^*\phi'(\gamma,x)}$$ i.e. $\judge{\gamma,\underline{x}:T(a(\f'(\gamma)))}{\psi(a(\gamma),\alpha''(\gamma)^*\underline{x})=\psi'(\gamma,\underline{x})}$ where $\alpha''(\gamma)\equiv \alpha(a(\gamma))\bullet a(\alpha'(\gamma))$.

\proof
By induction on the complexity of $\equivpair{\phi}{\psi}$.
\endproof
\end{lem}

\subsection{Properties of a restriction of the family}\label{section2.4.II}

In this subsection we restrict the family of the canonical context homotopy equivalences and prove additional properties satisfied by this restriction. Again, we remind that we use the notation $(\;\cdot\;;\;\cdot\;)$ to indicate any homotopy equivalence between contexts or types in context, reserving the symbol $\equivpair{\cdot}{\cdot}$ for the canonical ones. Moreover, we remind that a type judgement $\judge{\delta}{T(\delta):\type}$ is said to be an \textit{h-proposition} if $\judge{\delta,x_1,x_2:T(\delta)}{\alpha(\delta,x_1,x_2):x_1=x_2}$. We start by giving the following inductive definition:

\begin{defi}\label{concrete type}
A type judgement $\judge{\delta}{T(\delta):\type}$ in some context $\delta:\Delta$ \textbf{has h-propositional identities} (or is \textbf{with h-propositional identities}) if it belongs to the smallest family $\mathcal{F}$ of type judgements that satisfies the following clauses: \begin{itemize}
    \item a judgement $\judge{\gamma}{S:\type}$ (where $S$ is a atomic type) belongs to $\mathcal{F}$;
    \item a judgement $\judge{\gamma}{\Pi_{x : A(\gamma)}B(\gamma,x):\type}$, for some $\judge{\gamma}{A(\gamma)}:\type$ of $\mathcal{F}$ and some $\judge{\gamma,x:A(\gamma)}{B(\gamma,x):\type}$ of $\mathcal{F}$, belongs to $\mathcal{F}$;
    \item a judgement $\judge{\gamma}{\Sigma_{x : A(\gamma)}B(\gamma,x):\type}$, for some $\judge{\gamma}{A(\gamma)}:\type$ of $\mathcal{F}$ and some $\judge{\gamma,x:A(\gamma)}{B(\gamma,x):\type}$ of $\mathcal{F}$, belongs to $\mathcal{F}$;
    \item an h-proposition $\judge{\gamma}{s_1(\gamma)=s_2(\gamma)}:\type$, for some $\judge{\gamma}{A(\gamma)}:\type$ of $\mathcal{F}$ and some $\judge{\gamma}{s_1(\gamma),s_2(\gamma):A(\gamma)}$, belongs to $\mathcal{F}$.
\end{itemize}
\end{defi}

Despite our choice of terminology, we would like to clarify that the notion of a type with h-propositional identities does not refer to the types $T$ whose identity types $\judge{x,y:T}{x=y:\type}$ are h-propositions, i.e. the \textit{h-sets}. Rather, it characterises types whose equalities \textit{involved in their construction} are h-propositions.

\begin{defi}\label{concrete context}
Let $\gamma : \Gamma$ be a context $\gamma_1 : \Gamma_1,\gamma_2: \Gamma_2(\gamma_1),...,\gamma_n:\Gamma_n(\gamma_1,...,\gamma_{n-1})$, where $n$ might be $0$. For any $i\in\{1,...,n\}$, let $\gamma^i$ be the context $\gamma_1,...,\gamma_i$. We say that $\gamma$ \textbf{has h-propositional identities} (or is \textbf{with h-propositional identities}) if, for every $i \in \{1,...,n\}$, the judgement $\judge{\gamma^i}{\Gamma_i(\gamma^i)}$ has h-propositional identities.
\end{defi}

\begin{lem}
Let $\gamma:\Gamma$ be a context $\gamma_1 : \Gamma_1,\gamma_2: \Gamma_2(\gamma_1),...,\gamma_n:\Gamma_n(\gamma_1,...,\gamma_{n-1})$ and let $\delta:\Delta$ be a context $\delta_1 : \Delta_1,\delta_2: \Delta_2(\delta_1),...,\delta_m:\Delta_m(\delta_1,...,\delta_{m-1})$. If there is a canonical context homotopy equivalence $\equivpair{\ci}{\di}$ between $\gamma:\Gamma$ and $\delta:\Delta$ then $\gamma$ has h-propositional identities if and only if $\delta$ has h-propositional identities.

\proof
By \autoref{symmetry} we can assume w.l.o.g. that $\gamma$ has h-propositional identities and we are left to verify that $\delta$ has h-propositional identities as well. By induction on the complexity of $\equivpair{\ci}{\di}$. If $\equivpair{\ci}{\di}$ is of the form \texttt{(1)} then we are done. If $\equivpair{\ci}{\di}$ is of the form \texttt{(2)} then $\gamma  : \Gamma$ and $\delta:\Delta$ are of the form $\gamma'  : \Gamma',x : A(\gamma')$ and $\delta' :  \Delta',\underline{x}: A'(\delta')$ respectively and there are a canonical context homotopy equivalence $\equivpair{\f}{\g}$: \begin{center}$\judge{\gamma':\Gamma'}{\f(\gamma'):\Delta'}$\\ $\judge{\delta':\Delta'}{\g(\delta'):\Gamma'}$\end{center} and a canonical homotopy equivalence $\equivpair{f}{g}$: \begin{center} $\judge{\gamma':\Gamma',x:A(\gamma')}{f(\gamma',x):A'(\f(\gamma'))}$\\ $\judge{\gamma':\Gamma',\underline{x}':A'(\f(\gamma'))}{g(\gamma',\underline{x}'):A(\gamma')}$\end{center} between $A(\gamma')$ and $A'(\delta')$ relative to $\equivpair{\f}{\g}$, in such a way that the extension: \begin{center} $\judge{\gamma' : \Gamma',  x:A(\gamma')}{(\f,\f_{m+1})(\gamma',x):\Delta',A'}$\\ $\judge{\delta' : \Delta', \underline{x}:A'(\delta')}{(\g,\g_{n+1})(\delta',\underline{x}):\Gamma',A}$\end{center} of $\equivpair{\f}{\g}$ via $\equivpair{f}{g}$ (\autoref{context homotopy equivalence expansion}) is the canonical context homotopy equivalence $\equivpair{\ci}{\di}$. By inductive hypothesis and being $\gamma'$ a context with h-propositional identities, the context $\delta'$ has h-propositional identities. Moreover, by the following \autoref{p} and having $\judge{\gamma'}{A(\gamma'):\type}$ h-propositional identities, we deduce that $\judge{\delta'}{A'(\delta'):\type}$ has h-propositional identities. Therefore $\delta$ has h-propositional identities and we are done.
\endproof
\end{lem}

\begin{lem}\label{p}
Let $\gamma:\Gamma$ and $\delta:\Delta$ be contexts and let $(\f;\g)$ be a canonical context homotopy equivalence between $\gamma$ and $\delta$. Let $\judge{\gamma:\Gamma}{S(\gamma):\type}$ and $\judge{\gamma:\Delta}{T(\delta):\type}$ be such that there is a canonical homotopy equivalence $\equivpair{\phi}{\psi}$ between $S(\gamma)$ and $T(\delta)$ relative to $\equivpair{\f}{\g}$. Then $\judge{\gamma:\Gamma}{S(\gamma):\type}$ has h-propositional identities if and only if $\judge{\gamma:\Delta}{T(\delta):\type}$ has h-propositional identities.

\proof
By induction on the complexity of $\equivpair{\phi}{\psi}$. If $\equivpair{\phi}{\psi}$ is of the form \texttt{(a)} then we are done. If $\equivpair{\phi}{\psi}$ is of the form \texttt{(b)} then $S(\gamma)$ is of the form $\Pi_{x:A(\gamma)}B(\gamma,x)$ and $T(\delta)$ is of the form $\Pi_{x':A'(\delta)}B(\delta,x')$, where $\judge{\gamma}{A(\gamma):\type}$ and $\judge{\gamma,x}{B(\gamma,x):\type}$ have h-propositional identities, $A(\gamma)$ is canonically homotopy equivalent to $A'(\delta)$ relative to $(\f;\g)$ and $B(\gamma,x)$ is canonically homotopy equivalent to $B'(\delta,x')$ relative to the extension of $(\f;\g)$ via the given canonical homotopy equivalence between $A(\gamma)$ and $A'(\delta)$ relative to $(\f;\g)$. By inductive hypothesis, the judgements $\judge{\delta}{A'(\delta):\type}$ and $\judge{\delta,x'}{B'(\delta,x'):\type}$ have h-propositional identities, hence $T(\delta)$ has h-propositional identities and we are done. If $\equivpair{\phi}{\psi}$ is of the form \texttt{(c)} then we infer that $T(\delta)$ has h-propositional identities as for the case \texttt{(b)}. Finally, if $\equivpair{\phi}{\psi}$ is of the form \texttt{(d)} then $S(\gamma)$ is an h-proposition of the form $s_1(\gamma)=s_2(\gamma)$ for some $\judge{\gamma}{A(\gamma):\type}$ with h-propositional identities and some $\judge{\gamma}{s_1(\gamma),s_2(\gamma):A(\gamma)}$, while $T(\delta)$ is of the form $t_1(\delta)=t_2(\delta)$ for some $\judge{\delta}{A'(\delta):\type}$ and some $\judge{\delta}{t_1(\delta),t_2(\delta):A'(\delta)}$, being $A(\gamma)$ and $A'(\delta)$ canonically homotopy equivalent relative to $(\f;\g)$ in such a way that the canonical homotopy equivalence propositionally identifies $s_i(\gamma)$ and $t_i(\gamma)$, for $i \in  \{1,2\}$. By inductive hypothesis, the judgement $\judge{\delta}{A'(\delta):\type}$ has h-propositional identities. Moreover $t_1(\delta)=t_2(\delta)$ happens to be an h-proposition, since it is homotopy equivalent (via $\equivpair{\phi}{\psi}$) to an h-set. Hence $T(\delta)$ has h-propositional identities and we are done.
\endproof
\end{lem}

\begin{prop}\label{only one canonical context homotopy equivalence}
Let $\gamma:\Gamma$ and $\delta:\Delta$ be contexts with h-propositional identities. If $\equivpair{\ci}{\di}$ and $\equivpair{\ci'}{\di'}$ are canonical context homotopy equivalences between $\gamma:\Gamma$ and $\delta:\Delta$, then $\judge{\gamma:\Gamma}{\ci(\gamma)=\ci'(\gamma)}$ (i.e. $\judge{\delta:\Delta}{\di(\delta)=\di'(\delta)}$).

\proof
By induction on the complexity of $\equivpair{\ci}{\di}$. If $\equivpair{\ci}{\di}$ is of the form \texttt{(1)} then $\gamma : \Gamma$ and $\delta : \Delta$ are the empty contexts hence $\ci\equiv \ci'$ and $\di\equiv\di'$. If $\equivpair{\ci}{\di}$ is of the form \texttt{(2)} then $\gamma  : \Gamma$ and $\delta:\Delta$ are of the form $\gamma'  : \Gamma',x : A(\gamma')$ and $\delta' :  \Delta',\underline{x}: A'(\delta')$ respectively and there are a canonical context homotopy equivalence $\equivpair{\f}{\g}$: \begin{center}$\judge{\gamma':\Gamma'}{\f(\gamma'):\Delta'}$\\ $\judge{\delta':\Delta'}{\g(\delta'):\Gamma'}$\end{center} and a canonical homotopy equivalence $\equivpair{f}{g}$ : \begin{center} $\judge{\gamma':\Gamma',x:A(\gamma')}{f(\gamma',x):A'(\f(\gamma'))}$\\ $\judge{\gamma':\Gamma',\underline{x}':A'(\f(\gamma'))}{g(\gamma',\underline{x}'):A(\gamma')}$\end{center} between $A(\gamma')$ and $A'(\delta')$ relative to $\equivpair{\f}{\g}$, in such a way that the extension: \begin{center} $\judge{\gamma' : \Gamma',  x:A(\gamma')}{(\f,\f_{n+1})(\gamma',x):\Delta',A'}$\\ $\judge{\delta' : \Delta', \underline{x}:A'(\delta')}{(\g,\g_{n+1})(\delta',\underline{x}):\Gamma',A}$\end{center} of $\equivpair{\f}{\g}$ via $\equivpair{f}{g}$ (\autoref{context homotopy equivalence expansion}) is the given $\equivpair{\ci}{\di}$. Therefore $\equivpair{\ci'}{\di'}$ cannot be of the form \texttt{(1)} as (e.g.) $\delta:\Delta$ is not the empty context, which means that $\equivpair{\ci'}{\di'}$, being canonical, is of the form \texttt{(2)}. Hence, as before, there are a canonical context homotopy equivalence $\equivpair{\f'}{\g'}$: \begin{center}$\judge{\gamma':\Gamma'}{\f'(\gamma'):\Delta'}$\\ $\judge{\delta':\Delta'}{\g'(\delta'):\Gamma'}$\end{center} and a canonical homotopy equivalence $\equivpair{f'}{g'}$ : \begin{center} $\judge{\gamma':\Gamma',x:A(\gamma')}{f'(\gamma',x):A'(\f'(\gamma'))}$\\ $\judge{\gamma':\Gamma',\underline{\underline{x}}':A'(\f'(\gamma'))}{g'(\gamma',\underline{\underline{x}}'):A(\gamma')}$\end{center} between $A(\gamma')$ and $A'(\delta')$ relative to $\equivpair{\f'}{\g'}$, in such a way that the extension: \begin{center} $\judge{\gamma' : \Gamma',  x:A(\gamma')}{(\f',\f'_{n+1})(\gamma',x):\Delta',A'}$\\ $\judge{\delta' : \Delta', \underline{x}:A'(\delta')}{(\g',\g'_{n+1})(\delta',\underline{x}):\Gamma',A}$\end{center} of $\equivpair{\f'}{\g'}$ via $\equivpair{f'}{g'}$ (\autoref{context homotopy equivalence expansion}) is the given $\equivpair{\ci'}{\di'}$. Now $\judge{\gamma':\Gamma'}{\alpha(\gamma'):\f(\gamma')=\f'(\gamma')}$ (i.e. $\judge{\delta':\Delta'}{\tilde{\alpha}(\delta'):\g(\delta')=\g'(\delta')}$) by inductive hypothesis.

Let us pretend that we know that: $$\judge{\gamma',x}{f'(\gamma',x)=\alpha(\gamma')^*f'(\gamma',x)}$$ (\,i.e. $\judge{\gamma',\underline{\underline{x}}'}{g(\gamma',\alpha(\gamma')^*\underline{\underline{x}}')=g'(\gamma',\underline{\underline{x}}')}$\,), that is: \begin{center} $\judge{\gamma',x}{\f_{n+1}(\gamma',x)=\alpha(\gamma')^*\f'_{n+1}(\gamma',x)}$ (\,i.e. $\judge{\delta',\underline{x}}{\g_{n+1}(\delta',\underline{x})=\tilde{\alpha}(\delta')^*\g'_{n+1}(\delta',\underline{x})}$\,).
\end{center} Then: $$\judge{\gamma',x}{(\f,\f_{n+1})(\gamma',x)=(\f',\f'_{n+1})(\gamma',x)}$$ (\,i.e. $\judge{\delta',\underline{x}}{(\g,\g_{n+1})(\delta',\underline{x})=(\g',\g'_{n+1})(\delta',\underline{x})}$\,) that is: $$\judge{\gamma:\Gamma}{\ci(\gamma)=\ci'(\gamma)}$$ (\,i.e. $\judge{\delta:\Delta}{\di(\delta)=\di'(\delta)}$\,).

Hence we are done if $\judge{\gamma',x}{f'(\gamma',x)=\alpha(\gamma')^*f'(\gamma',x)}$ (\,i.e. $\judge{\gamma',\underline{\underline{x}}'}{g(\gamma',\alpha(\gamma')^*\underline{\underline{x}}')=g'(\gamma',\underline{\underline{x}}')}$\,), but this is actually the case by the following \autoref{only one canonical homotopy equivalence}.
\endproof
\end{prop}

\begin{lem}\label{only one canonical homotopy equivalence}
Let $\gamma:\Gamma$ and $\delta:\Delta$ be contexts together with context homotopy equivalences $(\f;\g)$ and $(\f';\g')$ such that: \begin{center}  $\judge{\gamma:\Gamma}{\alpha(\gamma):\f(\gamma)=\f'(\gamma)}$ (i.e. $\judge{\delta:\Delta}{\tilde{\alpha}(\delta):\g(\delta)=\g'(\delta)}$) \end{center} and let $\judge{\gamma:\Gamma}{S(\gamma):\type}$ and $\judge{\delta :\Delta}{T(\delta):\type}$ have h-propositional identities. Moreover, let $\equivpair{\phi}{\psi}$ be a canonical homotopy equivalence between $S(\gamma)$ and $T(\delta)$ relative to $(\f;\g)$ and let $\equivpair{\phi'}{\psi'}$ be a canonical homotopy equivalence between $S(\gamma)$ and $T(\delta)$ relative to $(\f';\g')$. Then $\judge{\gamma,s:S(\gamma)}{\phi(\gamma,s)=\alpha(\gamma)^*\phi'(\gamma,s)}$ i.e. $\judge{\gamma,\underline{t}':T(\f'(\gamma))}{\psi(\gamma,\alpha(\gamma)^*\underline{t}')=\psi'(\gamma,\underline{t}')}$.

\proof
By induction on the complexity of $S(\gamma)$.

\medskip

If $S(\gamma)$ is a base type $S$ then, being $\equivpair{\phi}{\psi}$ a canonical homotopy equivalence, the context $\delta$ needs to be $\gamma$ and $T(\delta)$ needs to be $S$, hence both $\equivpair{\phi}{\psi}$ and $\equivpair{\phi'}{\psi'}$ need to be of the form \texttt{(a)} and we are done.

\medskip

If $S(\gamma)$ is of the form $\Pi_{x:A(\gamma)}B(\gamma,x)$ for some $\judge{\gamma}{A(\gamma):\type}$ with h-propositional identities and some $\judge{\gamma,x}{B(\gamma,x):\type}$ with h-propositional identities then, being $\equivpair{\phi}{\psi}$ a canonical homotopy equivalence, the type $T(\delta)$ needs to be of the form $\Pi_{x':A'(\delta)}B'(\delta,x')$ for some $\judge{\delta}{A'(\delta):\type}$ with h-propositional identities and some $\judge{\delta,x'}{B(\delta,x'):\type}$ with h-propositional identities. Moreover, there are a canonical homotopy equivalence $\equivpair{f_1}{g_1}$: \begin{center}$\judge{\gamma:\Gamma,x:A(\gamma)}{f_1(\gamma,x):A'(\f(\gamma))}$\\ $\judge{\gamma:\Gamma,\underline{x}':A'(\f(\gamma))}{g_1(\gamma,\underline{x}'):A(\gamma)}$\end{center} between $A(\gamma)$ and $A'(\delta)$ relative to $(\f;\g)$ and a canonical homotopy equivalence $\equivpair{f_2}{g_2}$: \begin{center}$\judge{\gamma:\Gamma,x:A(\gamma),y:B(\gamma,x)}{f_2(\gamma,x,y):B'(\f(\gamma),f_1(\gamma,x))}$ \\ $\judge{\gamma:\Gamma,x:A(\gamma),\underline{y}':B'(\f(\gamma),f_1(\gamma,x))}{g_2(\gamma,x,\underline{y}'):B(\gamma,x)}$\end{center} between $B(\gamma,x)$ and $B'(\delta,\underline{x})$ relative to the extension $(\f,\f_{m+1};\g,\g_{n+1})$ of $(\f;\g)$ via $(f_1;g_1)$ (see \autoref{context homotopy equivalence expansion}), such that the homotopy equivalence $( f^\Pi;g^\Pi)$ of \autoref{pi equiv} is the given $\equivpair{\phi}{\psi}$. Analogously, there are a canonical homotopy equivalence $\equivpair{f'_1}{g'_1}$: \begin{center}$\judge{\gamma:\Gamma,x:A(\gamma)}{f'_1(\gamma,x):A'(\f'(\gamma))}$\\ $\judge{\gamma:\Gamma,\underline{x}':A'(\f'(\gamma))}{g'_1(\gamma,\underline{x}'):A(\gamma)}$\end{center} between $A(\gamma)$ and $A'(\delta)$ relative to $(\f';\g')$ and a canonical homotopy equivalence $\equivpair{f'_2}{g'_2}$: \begin{center}$\judge{\gamma:\Gamma,x:A(\gamma),y:B(\gamma,x)}{f'_2(\gamma,x,y):B'(\f'(\gamma),f'_1(\gamma,x))}$ \\ $\judge{\gamma:\Gamma,x:A(\gamma),\underline{y}':B'(\f'(\gamma),f'_1(\gamma,x))}{g'_2(\gamma,x,\underline{y}'):B(\gamma,x)}$\end{center} between $B(\gamma,x)$ and $B'(\delta,\underline{x})$ relative to the extension $(\f',\f'_{m+1};\g',\g'_{n+1})$ of $(\f';\g')$ via $(f'_1;g'_1)$ (see \autoref{context homotopy equivalence expansion}), such that the homotopy equivalence $( f'^\Pi;g'^\Pi)$ of \autoref{pi equiv} is the given $\equivpair{\phi'}{\psi'}$. By inductive hypothesis: \begin{center} $\judge{\gamma,x:A(\gamma)}{\alpha_1(\gamma,x):f_1(\gamma,x)=\alpha(\gamma)^*f'_1(\gamma,x)}$ i.e. $\judge{\gamma,\underline{x}':A'(\f'(\gamma))}{g_1(\gamma,\alpha(\gamma)^*\underline{x}')=g'_1(\gamma,\underline{x}')}$ \end{center} and: \begin{center} $\judge{\gamma,x:A(\gamma),y:B(\gamma,x)}{f_2(\gamma,x,y)=(\alpha(\gamma),\alpha_1(\gamma,x))^*f'_2(\gamma,x,y)}$\\ i.e. $\judge{\gamma,x:A(\gamma),\underline{y}':B(\f'(\gamma),f'_1(\gamma,x))}{g_2(\gamma,x,(\alpha(\gamma),\alpha_1(\gamma,x))^*\underline{y}')=g'_2(\gamma,x,\underline{y}')}$. \end{center} Then: \[\begin{aligned}
g_2(\gamma,x,\ev(\alpha(\gamma)^*\underline{z}',f_1(\gamma,x)))&=g_2(\gamma,x,\alpha_1(\gamma,x)^*\ev(\alpha(\gamma)^*\underline{z}',\alpha(\gamma)^*f'_1(\gamma,x)))\\
&=g_2(\gamma,x,\alpha_1(\gamma,x)^*(\alpha(\gamma),\refl(\alpha(\gamma)^*f'_1(\gamma,x)))^*\ev(\underline{z}',f'_1(\gamma,x)))\\
&=g_2(\gamma,x,\alpha_1(\gamma,x)^*\alpha(\gamma)^*\ev(\underline{z}',f'_1(\gamma,x)))\\
&=g_2(\gamma,x,(\alpha(\gamma),\alpha_1(\gamma,x))^*\ev(\underline{z}',f'_1(\gamma,x)))\\
&=g'_2(\gamma,x,\ev(\underline{z}',f'_1(\gamma,x)))
\end{aligned}\] in context ${\gamma,\underline{z}':[\Pi_{\underline{x}:A'(\delta)}B'(\delta,\underline{x})](\f'(\gamma))}$, where the first and the fourth equalities follows by (generalised) path induction on $\alpha_1(\gamma,x)$ and the second and the third by multiple (generalised) path induction on $\alpha(\gamma)$. By propositional function extensionality, we deduce that: $$\judge{\gamma,\underline{z}':[\Pi_{\underline{x}:A'(\delta)}B'(\delta,\underline{x})](\f'(\gamma))}{g^\Pi(\gamma,\alpha(\gamma)^*\underline{z}')=g'^\Pi(\gamma,\underline{z}')}$$ and we are done.

\medskip

If $S(\gamma)$ is of the form $\Sigma_{x:A(\gamma)}B(\gamma,x)$ for some $\judge{\gamma}{A(\gamma):\type}$ with h-propositional identities and some $\judge{\gamma,x}{B(\gamma,x):\type}$ with h-propositional identities then, being $\equivpair{\phi}{\psi}$ a canonical homotopy equivalence, the type $T(\delta)$ needs to be of the form $\Sigma_{x':A'(\delta)}B'(\delta,x')$ for some $\judge{\delta}{A'(\delta):\type}$ with h-propositional identities and some $\judge{\delta,x'}{B(\delta,x'):\type}$ with h-propositional identities. Moreover, there are a canonical homotopy equivalence $\equivpair{f_1}{g_1}$: \begin{center}$\judge{\gamma:\Gamma,x:A(\gamma)}{f_1(\gamma,x):A'(\f(\gamma))}$\\ $\judge{\gamma:\Gamma,\underline{x}':A'(\f(\gamma))}{g_1(\gamma,\underline{x}'):A(\gamma)}$\end{center} between $A(\gamma)$ and $A'(\delta)$ relative to $(\f;\g)$ and a canonical homotopy equivalence $\equivpair{f_2}{g_2}$: \begin{center}$\judge{\gamma:\Gamma,x:A(\gamma),y:B(\gamma,x)}{f_2(\gamma,x,y):B'(\f(\gamma),f_1(\gamma,x))}$ \\ $\judge{\gamma:\Gamma,x:A(\gamma),\underline{y}':B'(\f(\gamma),f_1(\gamma,x))}{g_2(\gamma,x,\underline{y}'):B(\gamma,x)}$\end{center} between $B(\gamma,x)$ and $B'(\delta,\underline{x})$ relative to the extension $(\f,\f_{m+1};\g,\g_{n+1})$ of $(\f;\g)$ via $(f_1;g_1)$ (\autoref{context homotopy equivalence expansion}), such that the homotopy equivalence $( f^\Sigma;g^\Sigma)$ of \autoref{sigma equiv} is the given $\equivpair{\phi}{\psi}$. Analogously, there are a canonical homotopy equivalence $\equivpair{f'_1}{g'_1}$: \begin{center}$\judge{\gamma:\Gamma,x:A(\gamma)}{f'_1(\gamma,x):A'(\f'(\gamma))}$\\ $\judge{\gamma:\Gamma,\underline{x}':A'(\f'(\gamma))}{g'_1(\gamma,\underline{x}'):A(\gamma)}$\end{center} between $A(\gamma)$ and $A'(\delta)$ relative to $(\f';\g')$ and a canonical homotopy equivalence $\equivpair{f'_2}{g'_2}$: \begin{center}$\judge{\gamma:\Gamma,x:A(\gamma),y:B(\gamma,x)}{f'_2(\gamma,x,y):B'(\f'(\gamma),f'_1(\gamma,x))}$ \\ $\judge{\gamma:\Gamma,x:A(\gamma),\underline{y}':B'(\f'(\gamma),f'_1(\gamma,x))}{g'_2(\gamma,x,\underline{y}'):B(\gamma,x)}$\end{center} between $B(\gamma,x)$ and $B'(\delta,\underline{x})$ relative to the extension $(\f',\f'_{m+1};\g',\g'_{n+1})$ of $(\f';\g')$ via $(f'_1;g'_1)$ (\autoref{context homotopy equivalence expansion}), such that the homotopy equivalence $( f'^\Sigma;g'^\Sigma)$ of \autoref{sigma equiv} is the given $\equivpair{\phi'}{\psi'}$. By inductive hypothesis: \begin{center} $\judge{\gamma,x:A(\gamma)}{\alpha_1(\gamma,x):f_1(\gamma,x)=\alpha(\gamma)^*f'_1(\gamma,x)}$ \end{center} and: \begin{center} $\judge{\gamma,x:A(\gamma),y:B(\gamma,x)}{f_2(\gamma,x,y)=(\alpha(\gamma),\alpha_1(\gamma,x))^*f'_2(\gamma,x,y)}$
\end{center} hence: \[\begin{aligned}
\judge{\gamma,u:\Sigma_{x:A(\gamma)}B(\gamma,x)}{f^\Sigma(\gamma,u)&\equiv\langle f_1(\gamma,\pi_1u),f_2(\gamma,\pi_1u,\pi_2u) \rangle}\\
&=\langle \alpha(\gamma)^*f'_1(\gamma,\pi_1u),\alpha(\gamma)^*f'_2(\gamma,\pi_1u,\pi_2u) \rangle\\
&= \alpha(\gamma)^*\langle f'_1(\gamma,\pi_1u),f'_2(\gamma,\pi_1u,\pi_2u) \rangle\\
&\equiv \alpha(\gamma)^*f'^\Sigma(\gamma,u)
\end{aligned}\] where the second propositional equality follows by multiple (generalised) path induction on $\alpha(\gamma)$. We are done.

\medskip

If $S(\gamma)$ is an h-proposition $s_1(\gamma)=s_2(\gamma)$ for some $\judge{\gamma}{A(\gamma)}:\type$ with h-propositional identities and some $\judge{\gamma}{s_1(\gamma),s_2(\gamma):A(\gamma)}$ then, being $\equivpair{\phi}{\psi}$ a canonical homotopy equivalence, the type $T(\delta)$ needs to be an h-proposition $t_1(\delta)=t_2(\delta)$ for some $\judge{\delta}{A'(\delta)}:\type$ with h-propositional identities and some $\judge{\delta}{t_1(\delta),t_2(\delta):A'(\delta)}$. We observe that $\phi(\gamma,p):t_1(\f(\gamma))=t_2(\f(\gamma))$ and $\phi'(\gamma,p):t_1(\f'(\gamma))=t_2(\f'(\gamma))$ in context $\gamma,p:s_1(\gamma)=s_2(\gamma)$, hence $\alpha(\gamma)^*\phi'(\gamma,p):t_1(\f(\gamma))=t_2(\f(\gamma))$. Since $t_1(\f(\gamma))=t_2(\f(\gamma))$ an h-proposition, then $\judge{\gamma,p}{\phi(\gamma,p)=\alpha(\gamma)^*\phi'(\gamma,p)}$ and we are done. \endproof
\end{lem}

\begin{rem}
We observe that the soundness of the last paragraph of the previous proof is the key reason why we need to work in the restricted family of type judgements of \autoref{concrete type}.
\end{rem}

In \autoref{section2.5} a further restriction of the contexts of $\ptt$ is adopted.

\section{Making the syntax of Propositional Type Theory\texorpdfstring{\\}{} into a category with attributes}\label{section2.5}

In this section we use the syntax of $\ptt$ in order to define a model of $\ett$. We remind that a type judgement $\judge{\delta}{T(\delta):\type}$ is said to be an \textit{h-set} if every type judgement $\judge{\delta,x_1,x_2:T(\delta)}{x_1=x_2:\type}$ is an h-proposition. We start by defining a further restriction on the type family that we allow to build contexts:

\begin{defi}\label{simple type}
A type judgement $\judge{\delta}{T(\delta):\type}$ in some context $\delta:\Delta$ is \textbf{h-elementary} if it belongs to the smallest family $\mathcal{F}$ of type judgements that satisfies the following clauses: \begin{itemize}
    \item a judgement $\judge{\gamma}{S:\type}$ (where $S$ is an atomic type) belongs to $\mathcal{F}$, whenever $S$ is an h-set;
    \item a judgement $\judge{\gamma}{\Pi_{x : A(\gamma)}B(\gamma,x):\type}$, for some $\judge{\gamma}{A(\gamma)}:\type$ of $\mathcal{F}$ and some $\judge{\gamma,x:A(\gamma)}{B(\gamma,x):\type}$ of $\mathcal{F}$, belongs to $\mathcal{F}$;
    \item a judgement $\judge{\gamma}{\Sigma_{x : A(\gamma)}B(\gamma,x):\type}$, for some $\judge{\gamma}{A(\gamma)}:\type$ of $\mathcal{F}$ and some $\judge{\gamma,x:A(\gamma)}{B(\gamma,x):\type}$ of $\mathcal{F}$, belongs to $\mathcal{F}$;
    \item a judgement $\judge{\gamma}{s_1(\gamma)=s_2(\gamma)}:\type$, for some $\judge{\gamma}{A(\gamma)}:\type$ of $\mathcal{F}$ and some $\judge{\gamma}{s_1(\gamma),s_2(\gamma):A(\gamma)}$, belongs to $\mathcal{F}$.
\end{itemize}
\end{defi}

The h-elementary types and the types with h-propositional identities (see \autoref{concrete type}) differ in their restrictions: while the latter only restrict the identity type formation to when it produces h-propositions, h-elementary types further restrict the atomic types to the h-sets.

\begin{defi}\label{simple context}
Let $\gamma : \Gamma$ be a context $\gamma_1 : \Gamma_1,\gamma_2: \Gamma_2(\gamma_1),...,\gamma_n:\Gamma_n(\gamma_1,...,\gamma_{n-1})$, where $n$ might be $0$. For any $i\in\{1,...,n\}$, let $\gamma^i$ be the context $\gamma_1,...,\gamma_i$. We say that $\gamma$ is \textbf{h-elementary} if, for every $i \in \{1,...,n\}$, the judgement $\judge{\gamma^i}{\Gamma_i(\gamma^i)}$ is h-elementary.
\end{defi}

We recall that in an intensional type theory the family of the h-sets is closed under $\Pi$, $\Sigma$, and $=$ operations. The same holds in $\ptt$, i.e. whenever $\judge{\gamma}{A(\gamma):\type}$ and $\judge{\gamma,x:A(\gamma)}{B(\gamma,x):\type}$ are h-sets, the judgements: $$\judge{\gamma}{\Sigma_{x:A(\gamma)}B(x,\gamma)}\textnormal{, }\judge{\gamma}{\Pi_{x:A(\gamma)}B(x,\gamma)}\textnormal{, and }\judge{\gamma,x,y:A(\gamma)}{x=y}$$ are h-sets as well. The proof works formally as in the intensional case. We refer the reader to \cite[Section 2.3]{rijke2015sets} and to \cite[Chapter 3]{hottbook}. Therefore, by induction on the complexity of an h-elementary type, we infer that:

\begin{cor}\label{simpleisset}
Every h-elementary type is an h-set.
\end{cor} and therefore, again by induction on the complexity of an h-elementary type, one proves that:

\begin{rem}
Every h-elementary type has h-propositional identities, hence every h-elementary context is a context with h-propositional identities. In particular, the results of \autoref{section2.4.II} continue being true for the family of h-elementary types and contexts.
\end{rem}

By induction on the complexity of a canonical homotopy equivalence, one proves that:

\begin{lem}\label{[]=>||types}
    Let $\gamma:\Gamma$ and $\delta:\Delta$ be h-elementary contexts together with a context homotopy equivalence $(\f;\g)$ such that $|\gamma:\Gamma|\equiv|\delta:\Delta|$ and: \begin{center} $\judgext{|\gamma|}{|\f(\gamma)|\equiv|\gamma|}$ \; (\;i.e. $\judgext{|\gamma|}{|\g(\delta)|\equiv|\gamma|}$\;)
    \end{center} and let $\judge{\gamma:\Gamma}{S(\gamma):\type}$ and $\judge{\delta :\Delta}{T(\delta):\type}$ be h-elementary. Moreover, let $\equivpair{\phi}{\psi}$ be a canonical homotopy equivalence between $S(\gamma)$ and $T(\delta)$ relative to $(\f;\g)$. Then $\judgext{|\gamma|}{|S(\gamma)|\equiv|T(\delta)|}$ \;(\:$\equiv|T(\f(\gamma))|$\;) and: \begin{center} $\judgext{|\gamma|,|s|}{|\phi(\gamma,s)|\equiv|s|}$ \; (\;i.e. $\judgext{|\gamma|,|t'|\equiv|t|\equiv|s|}{|\psi(\delta,t')|\equiv|s|}$\;).
    \end{center} 
\end{lem} and therefore, by induction on the complexity of a canonical context homotopy equivalence, one infers that:

\begin{prop}\label{[]=>||contexts}
    Let $\gamma:\Gamma$ and $\delta:\Delta$ be h-elementary contexts together with a canonical context homotopy equivalence $\equivpair{\ci}{\di}$ from $\gamma$ to $\delta$. Then $|\gamma|\equiv|\delta|$ and $\judgext{|\gamma|}{|\ci(\gamma)|=|\gamma|}$ i.e. $\judgext{|\gamma|}{|\di(\delta)|=|\gamma|}$.
\end{prop}

\subsection{The semantic context category}\label{section2.5.I}

We say that two h-elementary contexts of $\ptt$ are equivalent if there is a canonical context homotopy equivalence between them. This request defines an equivalence relation: its reflexivity follows by \autoref{there are context identities}, its symmetry follows by \autoref{symmetry} and its transitivity from \autoref{transitivity}.

From now on, whenever we speak about contexts and types of $\ptt$, we will actually refer to \textit{h-elementary} contexts and types of $\ptt$.

Let us assume that we are given two (h-elementary) context morphisms $\gamma \xrightarrow{f(\gamma)}\delta$ and $\gamma' \xrightarrow{f'(\gamma')}\delta'$ (of $\ptt$). We say that $f(\gamma)$ and $f'(\gamma')$ are equivalent if and only if: \begin{enumerate}
    \item the contexts $\gamma$ and $\gamma'$ are equivalent;
    \item the contexts $\delta$ and $\delta'$ are equivalent;
    \item if $\equivpair{\ci}{\di}$ is a canonical context homotopy equivalence $\gamma \to \gamma'$ and $\equivpair{\ci'}{\di'}$ is one $\delta\to \delta'$, then:\[\begin{tikzcd}
	\gamma && \delta \\
	\\
	{\gamma'} && {\delta'}
	\arrow["{\ci'(\delta)}", from=1-3, to=3-3]
	\arrow["{\ci(\delta)}"', from=1-1, to=3-1]
	\arrow["{f(\gamma)}", from=1-1, to=1-3]
	\arrow["{f'(\gamma')}"', from=3-1, to=3-3]
	\arrow[shift left=0.5, shorten <=27pt, shorten >=27pt, no head, from=3-1, to=1-3]
	\arrow[shift right=0.5, shorten <=27pt, shorten >=27pt, no head, from=3-1, to=1-3]
\end{tikzcd}\] i.e. $\judge{\gamma}{\ci'(f(\gamma))=f'(\ci(\delta))}$, where we observe that this condition does not depend on the choice of $\equivpair{\ci}{\di}$ and $\equivpair{\ci'}{\di'}$ by \autoref{only one canonical context homotopy equivalence} (we remind that h-elementary contexts have h-propositional identities).
\end{enumerate}

Then we observe that:

\begin{prop}\label{semantic context cat}
    There is a category $\mathcal{C}$ such that: \begin{itemize}
        \item the objects of $\mathcal{C}$ are the equivalence classes $[\gamma:\Gamma]$ of h-elementary contexts;
        \item the arrows $[\gamma:\Gamma]\to[\delta:\Delta]$ of $\mathcal{C}$ are the equivalence classes $[\judge{\gamma':\Gamma'}{f(\gamma'):\Delta'}]$ of morphisms of contexts $\judge{\gamma':\Gamma'}{f(\gamma'):\Delta'}$ where $[\gamma':\Gamma']=[\gamma]$ and $[\delta':\Delta']=[\delta]$;
        \item if we are given arrows:\begin{center} $[\judge{\gamma':\Gamma'}{f(\gamma'):\Delta'}]:[\gamma:\Gamma]\to[\delta:\Delta]$ and $[\judge{\delta'':\Delta''}{g(\delta''):\Omega'}]:[\delta:\Delta]\to[\omega:\Omega]$ \end{center} then, whenever $\equivpair{\ci}{\di}$ is a canonical homotopy equivalence $\delta'\to\delta''$, the composition arrow: \begin{center} $[\judge{\delta'':\Delta''}{g(\delta''):\Omega'}]\circ[\judge{\gamma':\Gamma'}{f(\gamma'):\Delta'}]:[\gamma:\Gamma]\to[\omega:\Omega]$ \end{center} is the arrow $[\judge{\gamma':\Gamma'}{g(\ci(f(\gamma'))):\Omega'}].$
    \end{itemize} Moreover, the category $\mathcal{C}$ has a terminal object.
\proof Observe that the composition operation is in fact well-defined: if we are given: \begin{center}$[\judge{\gamma'':\Gamma''}{f'(\gamma''):\Delta'''}]=[\judge{\gamma':\Gamma'}{f(\gamma'):\Delta'}]$ and $[\judge{\delta^{iv}:\Delta^{iv}}{g'(\delta^{iv}):\Omega''}]=[\judge{\delta'':\Delta''}{g(\delta''):\Omega'}]$\end{center} then: \[\begin{tikzcd}
	{\gamma'} && {\delta'} && {\delta''} && {\omega'} \\
	\\
	{\gamma''} && {\delta'''} && {\delta^{iv}} && {\omega''}
	\arrow["{\ci''(\delta')}"{description}, from=1-3, to=3-3]
	\arrow["{\ci^{iv}(\gamma')}"', from=1-1, to=3-1]
	\arrow["{f'(\gamma'')}"', from=3-1, to=3-3]
	\arrow["{\ci(\delta')}", from=1-3, to=1-5]
	\arrow["{\ci'(\delta''')}"', from=3-3, to=3-5]
	\arrow["{f(\gamma')}", from=1-1, to=1-3]
	\arrow["{g(\delta'')}", from=1-5, to=1-7]
	\arrow["{g'(\delta^{iv})}"', from=3-5, to=3-7]
	\arrow["{\ci'''(\delta'')}"{description}, from=1-5, to=3-5]
	\arrow["{\ci^{v}(\omega')}", from=1-7, to=3-7]
	\arrow[shift left=0.5, no head, from=3-1, to=1-3, shorten <=27pt, shorten >=27pt]
	\arrow[shift right=0.5, no head, from=3-1, to=1-3, shorten <=27pt, shorten >=27pt]
	\arrow[shift left=0.5, no head, from=3-3, to=1-5, shorten <=27pt, shorten >=27pt]
	\arrow[shift right=0.5, no head, from=3-3, to=1-5, shorten <=27pt, shorten >=27pt]
	\arrow[shift left=0.5, no head, from=3-5, to=1-7, shorten <=27pt, shorten >=27pt]
	\arrow[shift right=0.5, no head, from=3-5, to=1-7, shorten <=27pt, shorten >=27pt]
\end{tikzcd}\] for any choice of canonical context homotopy equivalences: \[\begin{aligned}\equivpair{\ci'}{\di'}&:\delta'''\to\delta^{iv}\\
\equivpair{\ci''}{\di''}&:\delta'\to\delta'''\\
\equivpair{\ci'''}{\di'''}&:\delta''\to\delta^{iv}\\
\equivpair{\ci^{iv}}{\di^{iv}}&:\gamma' \to\gamma''\\
\equivpair{\ci^{v}}{\di^{v}}&:\omega' \to\omega''
\end{aligned}\] where the propositional commutativity of the inner square is a consequence of \autoref{transitivity} and \autoref{only one canonical context homotopy equivalence}. Hence: $$[\judge{\gamma'':\Gamma''}{g'(\ci'(f'(\gamma''))):\Omega''}]=[\judge{\gamma':\Gamma'}{g(\ci(f(\gamma'))):\Omega'}]$$ and we are done. This operation is associative and the identity of $[\gamma:\Gamma]$ is the class $[\judge{\gamma}{\gamma}]$, hence $\mathcal{C}$ is actually a category.

We observe that a terminal object of $\mathcal{C}$ is the class $[\;\_:\_\;]$ represented by the empty context. In fact, whenever we are given an object $[\gamma:\Gamma]$ of $\mathcal{C}$, an arrow $[\gamma:\Gamma]\to[\;\_:\_\;]$ is the one represented by the judgement $\judge{\gamma:\Gamma}{\_:\_}$. Moreover, it is unique: if we are given an arrow $[\gamma:\Gamma]\to[\;\_:\_\;]$ represented by some judgement $\judge{\gamma':\Gamma'}{\_:\_}$, then the diagram:  \[\begin{tikzcd}
	\gamma \\
	& {\_} \\
	{\gamma'}
	\arrow["\ci"', from=1-1, to=3-1]
	\arrow["{\_}"', from=3-1, to=2-2, bend right]
	\arrow["{\_}", from=1-1, to=2-2, bend left]
\end{tikzcd}\] commutes (for every $\equivpair{\ci}{\di}$ canonical $\gamma\to\gamma'$) even strictly by the terminality of $\_$ in the category of contexts of $\ptt$, hence $\judge{\gamma'}{\_}$ represents $[\judge{\gamma}{\_}]$. \endproof
\end{prop}

We conclude the current subsection with the following:

\begin{rem}\label{giocare}\text{ }

\begin{enumerate}
    \item By \autoref{there are context identities}, if we are given two parallel morphisms of contexts $\judge{\gamma}{f(\gamma),f'(\gamma):\Delta}$, then $[\judge{\gamma}{f(\gamma)}]=[\judge{\gamma}{f'(\gamma)}]$ precisely when $\judge{\gamma}{f(\gamma)=f'(\gamma)}$.
    \item By \autoref{symmetry} and by \autoref{transitivity}, a judgement $\judge{\gamma':\Gamma'}{f(\gamma'):\Gamma''}$ represents $[\judge{\gamma}{\gamma}]$ (where $[\gamma']=[\gamma]=[\gamma'']$) if and only if $\judge{\gamma'}{f(\gamma')=\ci(\gamma')}$ for some $\equivpair{\ci}{\di}$ between $\gamma'$ and $\gamma''$. In particular, if $\gamma'\equiv\gamma''$ then $\judge{\gamma'}{f(\gamma')}$ represents $[\judge{\gamma}{\gamma}]$ if and only if $\judge{\gamma'}{f(\gamma')=\gamma'}$, by \autoref{there are context identities}. Therefore $\judge{\gamma'}{\gamma'}$ represents $[\judge{\gamma}{\gamma}]$ whenever $[\gamma']=[\gamma]$.
    \item If we are given an arrow $[\judge{\gamma':\Gamma'}{f(\gamma'):\Delta'}]:[\gamma:\Gamma]\to[\delta:\Delta]$ then, if $\equivpair{\ci}{\di}$ and $\equivpair{\ci'}{\di'}$ are canonical between $\gamma$ and $\gamma'$ and between $\delta'$ and $\delta$ respectively, we infer by 1. (or by \autoref{symmetry} and by \autoref{only one canonical context homotopy equivalence}) that: $$[\judge{\gamma':\Gamma'}{f(\gamma'):\Delta'}]=[\judge{\gamma:\Gamma}{\ci'(f(\ci(\gamma))):\Delta}].$$ Hence without loss of generality we can always assume that we are given a representative of $[\judge{\gamma':\Gamma'}{f(\gamma'):\Delta'}]$ of the form $\judge{\gamma:\Gamma}{f'(\gamma):\Delta}$.
    \item If we are given an arrow $[\gamma:\Gamma]\to [\delta:\Delta]$ represented by $\judge{\gamma:\Gamma}{f(\gamma):\Delta}$ and an arrow $[\delta:\Delta]\to[\omega:\Omega]$ represented by $\judge{\delta:\Delta}{g(\delta):\Omega}$, then their composition is represented by $g(\ci(f(\gamma)))$ for any canonical homotopy equivalence $\equivpair{\ci}{\di}$ between $\delta$ and itself. Hence $g(f(\gamma))$ represents their composition as well, by \autoref{there are context identities}, \autoref{only one canonical context homotopy equivalence} and 1. of \autoref{giocare}.
\end{enumerate}
\end{rem}

\subsection{The presheaf of semantic types}\label{section2.5.II}

In this subsection we define a presheaf of semantic types $\tp$ associated to the category $\mathcal{C}$ of \autoref{semantic context cat}. If $\gamma:\Gamma$ and $\gamma':\Gamma'$ are equivalent h-elementary contexts and $\judge{\gamma}{A(\gamma):\type}$ and $\judge{\gamma'}{A'(\gamma'):\type}$ are h-elementary then we say that the judgements $\judge{\gamma}{A(\gamma):\type}$ and $\judge{\gamma'}{A'(\gamma'):\type}$ are equivalent if there is a canonical homotopy equivalence between $A(\gamma)$ and $A'(\gamma')$ relative to some canonical context homotopy equivalence between $\gamma$ and $\gamma'$. Equivalently (see \autoref{context homotopy equivalence expansion} and the notion of canonical context homotopy equivalence), if the contexts $\gamma,x:A(\gamma)$ and $\gamma',x':A(\gamma')$ are equivalent.

\begin{itemize}
    \item If we are given an h-elementary context $\gamma:\Gamma$, we define $\tp[\gamma:\Gamma]$ as the family of the classes $[\judge{\gamma'}{A(\gamma'):\type}]$ where $[\gamma']=[\gamma]$ and $\judge{\gamma'}{A(\gamma'):\type}$ is h-elementary.
    \item If we are given an arrow $[\judge{\gamma':\Gamma'}{f(\gamma'):\Delta'}]:[\gamma:\Gamma]\to[\delta:\Delta]$, we define the map $\tp[\judge{\gamma':\Gamma'}{f(\gamma'):\Delta'}]:\tp[\delta:\Delta]\to\tp[\gamma:\Gamma]$ as the one such that: $$[\judge{\delta''}{A(\delta''):\type}]\mapsto[\judge{\gamma}{A(\ci'(f(\ci(\gamma)))):\type}]$$ where $\equivpair{\ci}{\di}$ and $\equivpair{\ci'}{\di'}$ are canonical between $\gamma$ and $\gamma'$ and between $\delta'$ and $\delta''$ respectively. This relation is in fact a mapping:
    \begin{itemize}
        \item Assuming that $[\judge{\gamma'':\Gamma''}{f'(\gamma''):\Delta'''}]=[\judge{\gamma':\Gamma'}{f(\gamma'):\Delta'}]$ and letting $\equivpair{\ci''}{\di''}$ and $\equivpair{\ci'''}{\di'''}$ be canonical between $\gamma$ and $\gamma''$ and between $\delta'''$ and $\delta''$ respectively, we are left to verify that: $$[\judge{\gamma}{A(\ci'''(f'(\ci''(\gamma)))):\type}]=[\judge{\gamma}{A(\ci'(f(\ci(\gamma)))):\type}].$$ Let $\equivpair{\ci^{iv}}{\di^{iv}}$ and $\equivpair{\ci^{v}}{\di^{v}}$ be canonical between $\gamma'$ and $\gamma''$ and between $\delta'$ and $\delta'''$. Then: \[\begin{tikzcd}
	& {\gamma'} && {\delta'} \\
	\gamma &&&& {\delta''} \\
	& {\gamma''} && {\delta'''}
	\arrow["{f(\gamma')}", from=1-2, to=1-4]
	\arrow["{f'(\gamma'')}", from=3-2, to=3-4]
	\arrow["\ci", from=2-1, to=1-2]
	\arrow["{\ci''}"', from=2-1, to=3-2]
	\arrow["{\ci'}", from=1-4, to=2-5]
	\arrow["{\ci'''}"', from=3-4, to=2-5]
	\arrow[""{name=0, anchor=center, inner sep=0}, "{\ci^{iv}}"{description}, from=1-2, to=3-2]
	\arrow[""{name=1, anchor=center, inner sep=0}, "{\ci^v}"{description}, from=1-4, to=3-4]
	\arrow[shorten <=33pt, shorten >=33pt, Rightarrow, no head, from=3-2, to=1-4]
	\arrow[shorten <=13pt, shorten >=13pt, Rightarrow, no head, from=2-1, to=0]
	\arrow[shorten <=13pt, shorten >=13pt, Rightarrow, no head, from=1, to=2-5]
\end{tikzcd}\] by \autoref{transitivity}. By \autoref{there are context identities}, there is a canonical context homotopy equivalence $\equivpair{\ci^{vi}}{\di^{vi}}$ between $\gamma$ and itself such that $\equivpair{\ci^{vi}}{\di^{vi}}$ is pairwise homotopic to $(\judge{\gamma}{\gamma};\judge{\gamma}{\gamma})$. Hence: \[\begin{tikzcd}
	\gamma \\
	& {\delta''} \\
	\gamma
	\arrow["{\ci^{vi}}"', from=1-1, to=3-1]
	\arrow[""{name=0, anchor=center, inner sep=0}, "{\ci'(f(\ci(\gamma)))}", from=1-1, to=2-2, bend left]
	\arrow[""{name=1, anchor=center, inner sep=0}, "{\ci'''(f'(\ci''(\gamma)))}"', from=3-1, to=2-2, bend right]
	\arrow[shift right=3, shorten <=17pt, shorten >=17pt, Rightarrow, no head, from=0, to=1]
\end{tikzcd}\] and therefore we are done by \autoref{lemma z}.

\item Assuming that $[\judge{\delta'''}{A'(\delta'''):\type}]=[\judge{\delta''}{A(\delta''):\type}]$, and letting $\equivpair{\ci''}{\di''}$ be canonical between $\delta'$ and $\delta'''$, we are left to verify that: $$[\judge{\gamma}{A'(\ci''(f(\ci(\gamma)))):\type}]=[\judge{\gamma}{A(\ci'(f(\ci(\gamma)))):\type}].$$ Let $\equivpair{\ci'''}{\di'''}$ be a canonical homotopy equivalence between $\delta''$ and $\delta'''$. Then: $$\begin{aligned}\relax[\judge{\delta'}{A(\ci'(\delta)):\type}]&=[\judge{\delta''}{A(\delta''):\type}]\\&=[\judge{\delta''}{A'(\ci'''(\delta'')):\type}]\\&=[\judge{\delta'}{A'(\ci''(\delta'))}]\end{aligned}$$ by \autoref{semi destiny}, by \autoref{lemma w} and by \autoref{lemma z} respectively. Then we are done by \autoref{small reindexing} with $\delta\equiv \delta'$ and $a(\gamma)\equiv f(\ci(\gamma))$.
\end{itemize} \end{itemize} In the remainder of the current subsection we verify that: \begin{prop}\label{presheaf of semantic types}
    The mapping $\tp$ is a functor. \proof\text{ }\begin{itemize}
    \item Let us consider the diagram: $$[\gamma]\xrightarrow{[\judge{\gamma'}{f(\gamma'):\Delta'}]}[\delta]\xrightarrow{[\judge{\delta''}{g(\delta''):\Omega'}]}[\omega]$$ and let $[\judge{\omega''}{A(\omega''):\type}]$ be in $\tp[\omega]$. Let us consider canonical context homotopy equivalences:\[\begin{aligned}\equivpair{\ci}{\di}&:\gamma\to\gamma'\\
\equivpair{\ci'}{\di'}&:\delta'\to\delta\\
\equivpair{\ci''}{\di''}&:\delta\to\delta''\\
\equivpair{\ci'''}{\di'''}&:\omega' \to\omega'''\\
\equivpair{\ci^{iv}}{\di^{iv}}&:\delta' \to\delta''.
\end{aligned}\] Then, to verify that $\tp$ preserves the composition, it is enough to verify that: $$[\judge{\gamma}{A(\ci'''(g(\ci''(\ci'(f(\ci(\gamma)))))))}]=[\judge{\gamma}{A(\ci'''(g(\ci^{iv}(f(\ci(\gamma))))))}]$$ which is in fact true by \autoref{transitivity} and \autoref{lemma z}.

\item Let $[\gamma]$ be a context and let $[\judge{\gamma'}{A(\gamma'):\type}]$ be in $\tp[\gamma]$. As $[\gamma']=[\gamma]$, the morphism $\judge{\gamma'}{\gamma'}$ represents the identity of $[\gamma]$ (see 2. of \autoref{giocare}). Let us consider canonical context homotopy equivalences: \[\begin{aligned}\equivpair{\ci}{\di}&:\gamma'\to\gamma'\\
\equivpair{\ci'}{\di'}&:\gamma'\to\gamma'.
\end{aligned}\] Then $[\judge{\gamma'}{A(\ci'(\ci(\gamma')))}]=[\judge{\gamma'}{A(\gamma')}]$ by \autoref{semi destiny}, hence $\tp$ preserves the identities and we are done.\qedhere
\end{itemize}
\endproof
\end{prop}

We end the current subsection with the following:

\begin{rem}\label{simpatico}\text{ }

\begin{enumerate}
    \item By \autoref{semi destiny}, if we are given a semantic type $[\judge{\gamma'}{A(\gamma'):\type}]$ in semantic context $[\gamma:\Gamma]$, then: $$[\judge{\gamma}{A(\ci(\gamma)):\type}]=[\judge{\gamma'}{A(\gamma'):\type}]$$ for every $\equivpair{\ci}{\di}$ canonical $\gamma\to\gamma'$. Hence we can always assume that we are given a representative of a semantic type in semantic context $[\gamma:\Gamma]$ whose context coincides with $\gamma$ itself.
    \item If we are given an arrow $[\judge{\gamma':\Gamma'}{f(\gamma'):\Delta'}]:[\gamma:\Gamma]\to[\delta:\Delta]$, we observe that the image of $[\judge{\delta''}{A(\delta''):\type}]$ via $\tp[\judge{\gamma':\Gamma'}{f(\gamma'):\Delta'}]:\tp[\delta:\Delta]\to\tp[\gamma:\Gamma]$, which is $[\judge{\gamma}{A(\ci'(f(\ci(\gamma)))):\type}]$ by definition, also coincides with: $$[\judge{\gamma'}{A(\ci'(f(\gamma'))):\type}]$$ by \autoref{semi destiny}, being $\equivpair{\ci}{\di}$ and $\equivpair{\ci'}{\di'}$ canonical between $\gamma$ and $\gamma'$ and between $\delta'$ and $\delta''$. We use this particular presentation in \autoref{section2.5.III}.
\end{enumerate}
\end{rem}

\subsection{The semantic context extension}\label{section2.5.III}

The semantic context extension $-.-$ associated to $\mathcal{C}$ and $\tp$ (see \autoref{section2.5.I} and \autoref{section2.5.II} respectively) is defined as follows:
\begin{itemize}
    \item An object $([\gamma:\Gamma],[\judge{\gamma'}{A(\gamma'):\type}])$ of $\gr\tp$ is sent to: $$[\gamma:\Gamma].[\judge{\gamma'}{A(\gamma'):\type}]:=[\gamma',x:A(\gamma')].$$ This mapping is well-defined because $[\judge{\gamma'}{A(\gamma'):\type}]=[\judge{\gamma''}{A'(\gamma''):\type}]$ precisely when $[\gamma',x:A(\gamma')]=[\gamma'',\underline{x}:A'(\gamma'')]$, by the notion of canonical context homotopy equivalence.
    
    \item An arrow $([\delta:\Delta],[\judge{\delta'}{B(\delta'):\type}])\to ([\gamma:\Gamma],[\judge{\gamma'}{A(\gamma'):\type}])$ in $\gr\tp$ is an arrow $[\judge{\delta'':\Delta''}{f(\delta''):\Gamma''}]$ such that: $$\tp[\judge{\delta'':\Delta''}{f(\delta''):\Gamma''}]:[\judge{\gamma'}{A(\gamma'):\type}]\mapsto[\judge{\delta'}{B(\delta'):\type}].$$ Hence, if $\equivpair{\ci}{\di}$ is a canonical context homotopy equivalence $\gamma'' \to \gamma'$, then: $$[\judge{\delta'}{B(\delta'):\type}]=[\judge{\delta''}{A(\ci(f(\delta''))):\type}]$$ by 2. or \autoref{simpatico}. Moreover: $$[\judge{\gamma'}{A(\gamma'):\type}]=[\judge{\gamma''}{A(\ci(\gamma'')):\type}]$$ by \autoref{semi destiny}. We stipulate that $-.-$ maps such an arrow of $\gr\tp$ to the arrow: $$[\delta'',\underline{x}:A(\ci(f(\delta'')))]\xrightarrow{[\;\judge{\delta'',\;\underline{x}\;}{f(\delta''),\;\underline{x}}\;]}[\gamma'',x:A(\ci(\gamma''))].$$ Let us prove that this relation is actually a map. Let us consider a representative: $$[\judge{\delta''':\Delta'''}{f'(\delta'''):\Gamma'''}]=[\judge{\delta'':\Delta''}{f(\delta''):\Gamma''}]$$ and let $\equivpair{\ci'}{\di'}$ be canonical $\gamma'''\to\gamma'$. We are left to verify that the corresponding: $$[\delta''',\underline{x}':A(\ci'(f'(\delta''')))]\xrightarrow{[\;\judge{\delta''',\;\underline{x}'\;}{f'(\delta'''),\;\underline{x}'}\;]}[\gamma''',x':A(\ci'(\gamma'''))]$$ satisfies: $$[\;\judge{\delta'',\;\underline{x}\;}{f(\delta''),\;\underline{x}}\;]=[\;\judge{\delta''',\;\underline{x}'\;}{f'(\delta'''),\;\underline{x}'}\;].$$ Let $\equivpair{\ci''}{\di''}$ and $\equivpair{\ci'''}{\di'''}$ be canonical $\delta''\to \delta'''$ and $\gamma''\to\gamma'''$ and let us consider the diagram: \[\begin{tikzcd}
	{\delta''} && {\gamma''} \\
	&&& {\gamma'} \\
	{\delta'''} && {\gamma'''}
	\arrow["{f(\delta'')}", from=1-1, to=1-3]
	\arrow["{f'(\delta''')}"', from=3-1, to=3-3]
	\arrow[""{name=0, anchor=center, inner sep=0}, "\ci", from=1-3, to=2-4, bend left]
	\arrow["{\ci'}"', from=3-3, to=2-4, bend right]
	\arrow["{\ci''}", from=1-1, to=3-1]
	\arrow["{\ci'''}"', from=1-3, to=3-3]
	\arrow["{\beta(\delta'')}"{description}, shorten <=22pt, shorten >=22pt, Rightarrow, from=1-3, to=3-1]
	\arrow["{\alpha(\gamma'')}"{description}, shorten <=13pt, shorten >=14pt, Rightarrow, from=0, to=3-3]
\end{tikzcd}\] where $\judge{\gamma''}{\alpha(\gamma''):\ci(\gamma'')=\ci'(\ci'''(\gamma''))}$ exists by \autoref{transitivity} and \autoref{only one canonical context homotopy equivalence} and $\judge{\delta''}{\beta(\delta''):\ci'''(f(\delta''))=f'(\ci''(\delta''))}$ exists since $f(\delta'')$ and $f'(\delta''')$ represent the same arrow of $\mathcal{C}$. If: $$\omega(\delta'')\equiv \alpha(f(\delta''))\bullet \ci'(\beta(\delta'')):\ci(f(\delta''))=\ci'(f'(\ci''(\delta'')))$$ then by \autoref{lemma z} there exists a canonical homotopy equivalence: $$\equivpair{\phi}{\psi}$$ between $A(\ci(f(\delta'')))$ and $A(\ci'(f'(\delta''')))$ relative to $\equivpair{\ci''}{\di''}$ and such that: $$\judge{\delta'',\underline{x}:A(\ci(f(\delta'')))}{\omega(\delta'')^*\phi(\delta'',\underline{x})=\underline{x}}.$$ Anologously, again by \autoref{lemma z} there is a canonical homotopy equivalence $\equivpair{\phi'}{\psi'}$ between $A(\ci(\gamma''))$ and $A(\ci'(\gamma'''))$ relative to $\equivpair{\ci''}{\di''}$ and such that: $$\judge{\gamma'',x:A(\ci(\gamma''))}{\alpha(\gamma'')^*\phi'(\gamma'',x)=x}.$$ By the notion of canonical context homotopy equivalence and by \autoref{context homotopy equivalence expansion}, we obtain canonical context homotopy equivalences: \[\begin{aligned}
    \judge{\delta'',\underline{x}:A(\ci(f(\delta'')))&}{\ci''(\delta''):\Delta''',\;\phi(\delta'',\underline{x}):A(\ci'(f'(\ci''(\delta''))))}\\
    \judge{\gamma'',x:A(\ci(\gamma''))&}{\ci'''(\gamma''):\Gamma''',\;\phi'(\gamma'',x):A(\ci'(\ci'''(\gamma'')))}
\end{aligned}\] between $\delta'',\underline{x}:A(\ci(f(\delta'')))$ and $\delta''',\underline{x}':A(\ci'(f'(\delta''')))$ and between $\gamma'',x:A(\ci(\gamma''))$ and $\gamma''',x':A(\ci'(\gamma'''))$ respectively. Therefore, we are left to verify that: \[\begin{tikzcd}
	{\delta'',\underline{x}} && {\gamma'',x} \\
	\\
	{\delta''',\underline{x}'} && {\gamma''',x'}
	\arrow["{\ci''(\delta''),\;\phi(\delta'',\underline{x})}"', from=1-1, to=3-1]
	\arrow["{\ci'''(\gamma''),\;\phi'(\gamma'',x)}", from=1-3, to=3-3]
	\arrow["{f(\delta''),\;\underline{x}}", from=1-1, to=1-3]
	\arrow["{f'(\delta'''),\;\underline{x}'}"', from=3-1, to=3-3]
	\arrow[shorten <=27pt, shorten >=27pt, Rightarrow, no head, from=3-1, to=1-3]
\end{tikzcd}\] i.e. that: $$\judge{\delta'',\underline{x}}{\ci'''(f(\delta'')),\;\phi'(f(\delta''),\underline{x})=f'(\ci''(\delta'')),\;\phi(\delta'',\underline{x})}.$$ Since $\judge{\delta''}{\beta(\delta''):\ci'''(f(\delta''))=f'(\ci''(\delta''))}$, we are left to verify that $\judge{\delta'',\underline{x}}{\phi'(f(\delta''),\underline{x})=\beta(\delta'')^*\phi(\delta'',\underline{x})}$. But since $\judge{\delta'',\underline{x}}{\beta(\delta'')^*\phi(\delta'',\underline{x})=\ci'(\beta(\delta''))^*\phi(\delta'',\underline{x})}$ by multiple (generalised) path induction on $\beta(\delta'')$, we  are left to verify that: $$\judge{\delta'',\underline{x}}{\phi'(f(\delta''),\underline{x})=\ci'(\beta(\delta''))^*\phi(\delta'',\underline{x})}.$$ This is equivalent to verifying that: $$\judge{\delta'',\underline{x}}{\alpha(f(\delta''))^*\phi'(f(\delta''),\underline{x})=\alpha(f(\delta''))^*\ci'(\beta(\delta''))^*\phi(\delta'',\underline{x})}$$ and this is true as $\alpha(f(\delta''))^*\phi'(f(\delta''),\underline{x})=\underline{x}$ and $$\alpha(f(\delta''))^*\ci'(\beta(\delta''))^*\phi(\delta'',\underline{x})=\omega(\delta'')^*\phi(\delta'',\underline{x})=\underline{x}.$$
\end{itemize} We are left to verify that:

\begin{prop}\label{semantic context extension}
    The mapping $-.-$ defines a functor.
\proof \text{ }
\begin{itemize}
    \item Suppose that we are given two arrows: \[\begin{aligned}\;&([\omega:\Omega],[\judge{\omega'}{C(\omega'):\type}])\xrightarrow{[{g(\omega''):\Delta'''}]}([\delta:\Delta],[\judge{\delta'}{B(\delta'):\type}])\\
    &([\delta:\Delta],[\judge{\delta'}{B(\delta'):\type}])\xrightarrow{[{f(\delta''):\Gamma''}]} ([\gamma:\Gamma],[\judge{\gamma'}{A(\gamma'):\type}])\end{aligned}\] in $\gr\tp$. We can rewrite the objects as follows:
    \[\begin{aligned}
    ([\gamma:\Gamma],[\judge{\gamma'}{A(\gamma'):\type}])&=([\gamma:\Gamma],[\judge{\gamma''}{A(\ci(\gamma'')):\type}])\\
    ([\delta:\Delta],[\judge{\delta'}{B(\delta'):\type}])&=([\delta:\Delta],[\judge{\delta''}{A(\ci(f(\delta''))):\type}])\\
    ([\omega:\Omega],[\judge{\omega'}{C(\omega'):\type}])&=([\omega:\Omega],[\judge{\omega''}{A(\ci(f(\ci'(g(\omega''))))):\type}])
    \end{aligned}\] for some $\equivpair{\ci}{\di}$ and $\equivpair{\ci'}{\di'}$ canonical $\gamma''\to\gamma'$ and $\delta'''\to\delta''$ respectively. The composition of $\gr\tp$ sends this diagram to the arrow $[f(\ci'(g(\omega''))):\Gamma'']$ of the form: $$([\omega:\Omega],[\judge{\omega''}{A(\ci(f(\ci'(g(\omega''))))):\type}])\to([\gamma:\Gamma],[\judge{\gamma''}{A(\ci(\gamma'')):\type}]).$$ Now, the arrows $[f(\delta'')]$, $[g(\omega'')]$ and $[f(\ci'(g(\omega'')))]$ are sent by $-.-$ to the arrows:
    $$\begin{aligned}
    \relax [\delta'', \underline{x} : A(\ci(f(\delta'')))] &\xrightarrow{[f(\delta''),\underline{x}]}[\gamma'',x:A(\ci(\gamma''))]\\
    [\omega'',\underline{\underline{x}}:A(\ci(f(\ci'(g(\omega'')))))]&\xrightarrow{[g(\omega''),\underline{\underline{x}}]}[\delta''',\underline{x}':A(\ci(f(\ci'(\delta'''))))]\\
    [\omega'',\underline{\underline{x}}:A(\ci(f(\ci'(g(\omega'')))))]&\xrightarrow{[f(\ci'(g(\omega''))),\underline{\underline{x}}]}[\gamma'',x:A(\ci(\gamma''))]
    \end{aligned}$$
    of $\mathcal{C}$ respectively. In order to conclude that $-.-$ actually preserves the composition, we are left to verify that the composition of the first and the second of these yields the third. By \autoref{semi destiny}, there is $\equivpair{\phi}{\psi}$ canonical between $A(\ci(f(\ci'(\delta'''))))$ and $A(\ci(f(\delta'')))$ relative to $\equivpair{\ci'}{\di'}$ and such that $\judge{\delta''',\underline{x}'}{\phi(\delta''',\underline{x}')=\underline{x}'}$. Hence, by \autoref{context homotopy equivalence expansion} and by the notion of canonical context homotopy equivalence, there is a canonical context homotopy equivalence $\delta''',\underline{x}'\to\delta'',\underline{x}$ whose first component is $\ci'(\delta'''),\phi(\delta''',\underline{x}')$ and therefore the composition of $[f(\delta''),\underline{x}]\circ [g(\omega''),\underline{\underline{x}}]$ is represented by: $$\judge{\omega'',\underline{\underline{x}}}{f(\ci'(g(\omega''))),\phi(g(\omega''),\underline{\underline{x}})}$$ which is in fact propositionally equal to $\judge{\omega'',\underline{\underline{x}}}{f(\ci'(g(\omega''))),\underline{\underline{x}}}$ since $\judge{\omega'',\underline{\underline{x}}}{\phi(g(\omega''),\underline{\underline{x}})=\underline{\underline{x}}}$. We are done by 1. of \autoref{giocare}.
    \item By 2. of \autoref{giocare}, the identity over $([\gamma],[\judge{\gamma'}{A(\gamma'):\type}])$ is represented by the judgement $\judge{\gamma'}{\gamma'}$. By definition, it is sent by $-.-$ to the arrow: $$[\gamma',x':A(\ci(\gamma'))]\xrightarrow{[\judge{\gamma',x'}{\gamma',x'}]}[\gamma',x':A(\ci(\gamma'))]$$ for some $\equivpair{\ci}{\di}$ canonical $\gamma'\to\gamma'$. This is in fact the identity of $[\gamma',x':A(\ci(\gamma'))]=[\gamma',x:A(\gamma')]$. Hence identities are preserved and we are done.\qedhere
\end{itemize}
\endproof
\end{prop}

\subsection{The display map family}\label{section2.5.IV}

We define a cartesian natural transformation $P$ from the semantic context extension and the projection, as functors $\gr\tp\to\mathcal{C}$. Whenever $([\gamma],[\judge{\gamma'}{A(\gamma'):\type}])$ is an object of $\gr\tp$ then let the component of $P$ in $([\gamma],[\judge{\gamma'}{A(\gamma'):\type}])$ be the arrow: $$[\gamma',x:A(\gamma')]\xrightarrow{[\judge{\gamma',x\;}{\;\gamma'}]}[\gamma'].$$ The mapping $P$ happens to be well-defined: as long as $[\judge{\gamma''}{A'(\gamma''):\type}]=[\judge{\gamma'}{A(\gamma'):\type}]$ then the diagram: \[\begin{tikzcd}
	{\gamma',x} && {\gamma'} \\
	\\
	{\gamma'',\underline{x}} && {\gamma''}
	\arrow["{\ci(\gamma'),\phi(\gamma',x)}"', from=1-1, to=3-1]
	\arrow["{\ci(\gamma')}", from=1-3, to=3-3]
	\arrow["{\gamma'}", from=1-1, to=1-3]
	\arrow["{\gamma''}"', from=3-1, to=3-3]
\end{tikzcd}\] commutes even strictly whenever $\equivpair{\phi}{\psi}$ is canonical $A(\gamma')\to A'(\gamma'')$ relative to some canonical $\equivpair{\ci}{\di}$ between $\gamma'$ and $\gamma''$. Let us verify that:

\begin{prop}\label{display map family}
    The family $P$ is natural and cartesian.
    \proof\text{ }
\begin{itemize}
    \item \textit{Naturality.} Let us consider an arrow: $$[\judge{\delta'':\Delta''}{f(\delta''):\Gamma''}]:([\delta:\Delta],[\judge{\delta'}{B(\delta'):\type}])\to ([\gamma:\Gamma],[\judge{\gamma'}{A(\gamma'):\type}])$$ in $\gr\tp$. By 3. of \autoref{giocare}, we can assume without loss of generality that $\delta''\equiv\delta'$ and $\gamma''\equiv\gamma'$. If $\equivpair{\ci}{\di}$ is canonical between $\gamma'$ and itself, then: \[\begin{aligned}
    \relax [\judge{\gamma'}{A(\gamma'):\type}]&=[\judge{\gamma'}{A(\ci(\gamma')):\type}]\\
    [\judge{\delta'}{B(\delta'):\type}]&=[\judge{\delta'}{A(\ci(f(\delta'))):\type}]
\end{aligned}\] Then the images of $[\judge{\delta'}{f(\delta)':\Gamma'}]$ via the semantic context extension and via the projection admit the presentations: \[\begin{aligned}
\relax [\delta',x':A(\ci(f(\delta')))]&\xrightarrow{[\judge{\delta',x'}{f(\delta'),x'}]}[\gamma',x:A(\ci(\gamma'))]\\
[\delta']&\xrightarrow{[\judge{\delta'}{f(\delta')}]}[\gamma']
\end{aligned}\] and we are left to verify that they commute with $P$, in order to conclude the naturality of $P$ itself. By considering the representatives  $\judge{\delta',x'}{\delta'}$ and $\judge{\gamma',x}{\gamma'}$ of the $[\delta',x']$-component and the $[\gamma',x]$-component of $P$ respectively, we are done by 1. and 4. of \autoref{giocare}.
\item \textit{Cartesianity.} Again, let us consider an arrow: $$[\judge{\delta':\Delta'}{f(\delta'):\Gamma'}]:([\delta:\Delta],[\judge{\delta'}{B(\delta'):\type}])\to ([\gamma:\Gamma],[\judge{\gamma'}{A(\gamma'):\type}])$$ in $\gr\tp$ (the choice of its representative is justified by 3. of \autoref{giocare}) and let us verify that the commutative square: \[\begin{tikzcd}
	{[\delta',x']} && {[\gamma',x]} \\
	\\
	{[\delta']} && {[\gamma']}
	\arrow["{[f(\delta'),x']}", from=1-1, to=1-3]
	\arrow["{[f(\delta')]}", from=3-1, to=3-3]
	\arrow["{[\delta']}"', from=1-1, to=3-1]
	\arrow["{[\gamma']}", from=1-3, to=3-3]
\end{tikzcd}\] of $\mathcal{C}$ is a pullback of $\mathcal{C}$, being $x':A(\ci(f(\delta')))$ and $x:A(\ci(\gamma'))$ and being $\equivpair{\ci}{\di}$ canonical $\gamma\to\gamma'$. Hence, let us assume that we are given two arrows $\chi:[\omega]\to[\delta']$ and $\chi':[\omega]\to[\gamma',x]$ of $\mathcal{C}$ such that:
\[\begin{tikzcd}
	{[\omega]} && {[\gamma',x]} \\
	\\
	{[\delta']} && {[\gamma']}
	\arrow[from=1-1, to=1-3]
	\arrow["{[f(\delta')]}", from=3-1, to=3-3]
	\arrow[from=1-1, to=3-1]
	\arrow["{[\gamma']}", from=1-3, to=3-3]
\end{tikzcd}\] commutes. By 3. of \autoref{giocare} there are representatives of the form: \begin{center} $\judge{\omega}{\alpha_1(\omega):\Delta'}$ and $\judge{\omega}{\alpha_2(\omega):\Gamma',\alpha_3(\omega):A(\ci(\alpha_2(\omega)))}$ \end{center} of $\chi$ and $\chi'$ respectively and $\judge{\omega}{p(\omega):f(\alpha_1(\omega))=\alpha_2(\omega)}$, by 1. and 4. of \autoref{giocare}. Therefore $p(\omega)^*\alpha_3(\omega):A(\ci(f(\alpha_1(\omega))))$ and we obtain a morphism of contexts: $$\judge{\omega}{\alpha_1(\omega):\Delta',\;p(\omega)^*\alpha_3(\omega):A(\ci(f(\alpha_1(\omega))))}$$ which represents a morphism of semantic contexts $\underline{\chi}:[\omega]\to[\delta',x']$ in $\mathcal{C}$. Post-composing $\underline{\chi}$ via $[\delta',x']\xrightarrow{[\delta']}[\delta']$ and via $[\delta',x']\xrightarrow{[f(\delta'),x']}[\gamma',x]$ we get morphisms $[\omega]\to[\delta']$ and $[\omega]\to[\gamma',x]$ respectively represented by: \begin{center} $\judge{\omega}{\alpha_1(\omega):\Delta'}$ and $\judge{\omega}{f(\alpha_1(\omega)):\Gamma',\;p(\omega)^*\alpha_3(\omega):A(\ci(f(\alpha_1(\omega))))}$ \end{center} respectively (by 4. of \autoref{giocare}). The former is clearly $\chi$, while the latter is $\chi'$ because: \[\begin{aligned}
    \judge{\omega&}{p(\omega):f(\alpha_1(\omega))=\alpha_2(\omega)}\\
    \judge{\omega&}{p(\omega)^*\alpha_3(\omega)\equiv p(\omega)^*\alpha_3(\omega)}
\end{aligned}\] which means (1. of \autoref{giocare}) that the morphisms: $$f(\alpha_1(\omega)),\;p(\omega)^*\alpha_3(\omega)\text{ and }\alpha_2(\omega),\alpha_3(\omega)$$ represent the same arrow of $\mathcal{C}$.

Now, let $\underline{\chi}'$ be a morphism $[\omega]\to[\delta',x']$ such that $[\delta']\underline{\chi}'=\chi$ and $[f(\delta'),x']\underline{\chi}'=\chi'$. By 3. of \autoref{giocare}, the morphism $\underline{\chi}'$ is represented by a morphism of contexts of the form: $$\judge{\omega}{\beta_1(\omega):\Delta',\;\beta_2(\omega):A(\ci(f(\beta_1(\omega))))}.$$ Moreover: \[\begin{aligned} \judge{\omega&}{p_1(\omega):\alpha_1(\omega)=\beta_1(\omega)}\\ \judge{\omega&}{p_2(\omega):\alpha_2(\omega)=f(\beta_1(\omega))}\\ \judge{\omega&}{p_3(\omega):\alpha_3(\omega)=p_2(\omega)^*\beta_2(\omega)} \end{aligned}\] by 1. and 4. or \autoref{giocare}, hence $p_1(\omega):\alpha_1(\omega)=\beta_1(\omega)$ and: \[\begin{aligned}
p(\omega)^*\alpha_3(\omega)&=p(\omega)^*p_2(\omega)^*\beta_2(\omega)\\
&=(p(\omega)\bullet p_2(\omega))^*\beta_2(\omega)\\
&=f(p_1(\omega))^*\beta_2(\omega)\\
&=p_1(\omega)^*\beta_2(\omega)
\end{aligned}\] where the first and the second equalities are instances of propositional functorialities, the fourth holds by multiple (generalised) path induction on $p_1(\omega)$ and the third because $\gamma':\Gamma'$ is an h-elementary context and by \autoref{simpleisset}. Therefore, the judgements $\alpha_1(\omega),\;p(\omega)^*\alpha_3(\omega)$ and $\beta_1(\omega),\beta_2(\omega)$ represent the same arrow $[\omega]\to[\delta',x']$ i.e. $\underline{\chi}=\underline{\chi}'$.\qedhere
\end{itemize}
\endproof
\end{prop}

\begin{rem}
The last part of the proof of the cartesianity of the natural transformation $P$ is the crucial point where we needed to work with h-elementary contexts. Any other result that we obtained so far is true for the (generally larger) family of contexts with h-propositional identities.
\end{rem}

We summarise the content of the current section into the following:

\begin{thm}\label{category with attributes}
The category $\mathcal{C}$ of \autoref{semantic context cat}, together with the presheaf $\tp$ of \autoref{presheaf of semantic types}, the semantic context extension $-.-$ of \autoref{semantic context extension} and the natural transformation $P$ of \autoref{display map family}, forms a category with attributes $\C:=\cattypes$.
\end{thm}

We end the current section with the following:

\begin{rem}\label{[]=>||morphismofsemantictypes} A terminal object preserving functor from the base category of $\C$ to the base category of $\ettm$ mapping: $$(\;[\delta:\Delta]\xrightarrow{[\judge{\delta}{f(\delta):\Gamma}]}[\gamma:\Gamma]\;)\mapsto(\;|\delta|\xrightarrow{|\judge{\delta}{f(\delta):\Gamma}|}|\gamma|\;)=(\;|\delta|\xrightarrow{\judgext{|\delta|}{|f(\delta)|:|\Gamma|}}|\gamma|\;)$$ is well-defined by \autoref{[]=>||contexts} and extends to a morphism of semantic types $\C\to\ettm$ if we stipulate that: $$[\judge{\gamma}{A(\gamma):\type}]\;\mapsto\;|\judge{\gamma}{A(\gamma):\type}|\;=\;\judgext{|\gamma|}{|A(\gamma)|:\type}$$ (see \autoref{canonical interpretation} for more details). Observe in fact that this operation is well-defined by \autoref{[]=>||types} and is natural since $\judgext{|\gamma|}{|A(f(\delta))|\equiv|A|(|f(\delta)|)}$ by \autoref{canonical interpretation} (here we are implicitly using that: $$(\tp[\judge{\delta}{f(\delta):\Gamma}])[\judge{\gamma}{A(\gamma):\type}]=[\judge{\delta}{A(f(\delta)):\type}]$$ and this presentation of the substitution in $\C$ is justified by \autoref{there are context identities}, \autoref{semi destiny} and \autoref{small reindexing}). One can verify that the semantic context extension and the display map family are preserved.
\end{rem}

\autoref{section2.6} is devoted to proving that $\C$, together with the structure specified in \autoref{category with attributes}, is a model of $\ett$.

\section{A model of Extensional Type Theory}\label{section2.6}

In this section, we show that $\C$ verifies the requirements of \autoref{modelofsomeT} where $\T$ is $\ett$. We start by observing that, if $T$ is an atomic type of $\ett$ and $t:T$ is an atomic term of $\ett$, then $T$ is an atomic h-set of $\ptt$, hence the class $[\;\judge{}{T:\type}\;]$ is a semantic type of $\C$ in semantic context $[\;\_:\_\;]$ and $[\;\judge{}{t:T}\;]$ is a section of the corresponding display map, hence a semantic term of $[\;\judge{}{T:\type}\;]$ in $\C$. Therefore the mappings: \[\begin{aligned}
    T&\mapsto[\;\judge{}{T:\type}\;]\\
    t&\mapsto[\;\judge{}{t:T}\;]
\end{aligned}\] define a choice function as the one required in \autoref{modelofsomeT}.

After studying an opportune presentation of the semantic terms of $\C$ in \autoref{section2.6.I}, we show that $\C$ has semantic extensional identity types, semantic dependent product types, and semantic dependent sum types. We present the proof for the extensional identities in full form, while we only leave a sketch of the corresponding ones for dependent products and sums. Finally, we deduce the conservativity result in \autoref{section2.6.IV}.

\subsection{A presentation of the sections in the quotient syntax}\label{section2.6.I}

Let $[\judge{\gamma}{A(\gamma):\type}]$ be some semantic type in some semantic context $[\gamma:\Gamma]$ (by 1. of \autoref{simpatico} we do not lose generality if we assume this presentation for $[\judge{\gamma}{A(\gamma):\type}]$). Then, its component of $P$ admits the presentation $[\judge{\gamma,x:A(\gamma)}{\gamma}]$. Let us consider a section: $$[\gamma]\to [\gamma].[\judge{\gamma}{A(\gamma):\type}]=[\gamma,x]$$ of $[\judge{\gamma,x}{\gamma}]$, that, by 3. of \autoref{giocare}, admits the presentation $[\judge{\gamma}{\tilde{a}(\gamma):\Gamma,a(\gamma):A(\tilde{a}(\gamma))}]$. Then $\judge{\gamma}{\alpha(\gamma):\gamma=\tilde{a}(\gamma)}$, by 1. of \autoref{giocare} and being $[\judge{\gamma}{\tilde{a}(\gamma),a(\gamma)}]$ a section of $[\judge{\gamma,x}{\gamma}]$. Hence the morphism  of contexts: $$\judge{\gamma}{\gamma,\alpha(\gamma)^*a(\gamma):A(\gamma)}$$ is context propositionally equal to $\judge{\gamma}{\tilde{a}(\gamma),a(\gamma)}$ and therefore it continues representing the given section $[\gamma]\to[\gamma,x]$.

We conclude that:

\begin{rem}\label{nuovo gioco}
Without loss of generality, every section of a display map $[\judge{\gamma,x}{\gamma}]$ is of  the form: $$[\judge{\gamma}{\gamma,a(\gamma)}]$$ for some term $\judge{\gamma}{a(\gamma):A(\gamma)}$.
\end{rem}

\subsection{Semantic extensional identity types}\label{section2.6.II}

Let $[\judge{\gamma}{A(\gamma):\type}]$ be a semantic type in semantic context $[\gamma:\Gamma]$ (see 1. or \autoref{simpatico}). A presentation of $[\gamma:\Gamma].[\judge{\gamma}{A(\gamma):\type}]$ is $[\gamma,x:A(\gamma)]$, hence: $$(\tp[\judge{\gamma,x}{\gamma}])([\judge{\gamma}{A(\gamma):\type}])=[\judge{\gamma,x}{A(\ci(\gamma)):\type}]$$ for some canonical homotopy equivalence $\equivpair{\ci}{\di}$ between $\gamma$ and itself. In particular: $$(\tp[\judge{\gamma,x}{\gamma}])([\judge{\gamma}{A(\gamma):\type}])=[\judge{\gamma,x}{A(\gamma):\type}]$$ by Propostion \autoref{there are context identities} and by \autoref{semi destiny}. Therefore we obtain the presentation: $$[\gamma,x:A(\gamma),y:A(\gamma)]$$ of $([\gamma:\Gamma].[\judge{\gamma}{A(\gamma):\type}]).((\tp[\judge{\gamma,x}{\gamma}])([\judge{\gamma}{A(\gamma):\type}]))$. We define the semantic type $\id_{[\judge{\gamma}{A(\gamma):\type}]}$ in context $[\gamma,x,y]$ as the one represented by the type judgement: $$\judge{\gamma,x,y}{x=y:\type}$$ and the morphism of contexts $\refl_{[\judge{\gamma}{A(\gamma):\type}]}$ between $[\gamma,x]$ and $[\gamma,x,y].[\judge{\gamma,x,y}{x=y:\type}]=[\gamma,x,y,p:x=y]$ as the one represented by the context morphism: $$\judge{\gamma,x}{\gamma,x,x,\refl(x)}.$$ We verified that \textit{formation} and \textit{introduction} of \autoref{semantic extensional identity types} are satisfied. We are left to verify that:

\begin{prop} The remaining conditions of \autoref{semantic extensional identity types}---i.e. \textnormal{extensionality} and \textnormal{compatibility with the substitution}---are satisfied by the above choice of: $$\id_{[\judge{\gamma}{A(\gamma):\type}]}\text{ and }\refl_{[\judge{\gamma}{A(\gamma):\type}]}$$ hence $\C$ is equipped with semantic extensional identity types.
\proof\text{ }
\begin{itemize}
    \item \textit{Extensionality}. Let us consider two semantic terms: $$[\judge{\gamma}{\gamma:\Gamma,a(\gamma):A(\gamma)}]\textnormal{ and }[\judge{\gamma}{\gamma:\Gamma,b(\gamma):A(\gamma)}]$$ of the semantic type $[\judge{\gamma}{A(\gamma):\type}]$ in the semantic context $[\gamma:\Gamma]$. They admit such a presentation because of 3. of \autoref{giocare} and \autoref{nuovo gioco}. Now, if we consider $[\judge{\gamma}{\gamma:\Gamma,a(\gamma):A(\gamma)}]$ as an arrow of $\gr\tp$ of source $([\gamma],[\judge{\gamma}{A(\gamma):\type}])$ and target $([\gamma,x],[\judge{\gamma,x}{A(\gamma):\type}])$, its semantic context extension: $$[\judge{\gamma}{\gamma:\Gamma,a(\gamma):A(\gamma)}].[\judge{\gamma}{A(\gamma):\type}]$$ is: $$[\gamma,x]\xrightarrow{[\judge{\gamma}{\gamma,a(\gamma),x}]}[\gamma,x,y]$$ hence the arrow: $$[\gamma]\xrightarrow{[\judge{\gamma}{\gamma,b(\gamma):A(\gamma)}]}[\gamma,x]\xrightarrow{[\judge{\gamma}{\gamma:\Gamma,a(\gamma):A(\gamma)}].[\judge{\gamma}{A(\gamma):\type}]}[\gamma,x,y]$$ of \autoref{encode two terms in one morphism} happens to be represented by $\judge{\gamma}{\gamma,a(\gamma),b(\gamma)}$, because of 4. of \autoref{giocare}.
    
    By \autoref{there are context identities}, by \autoref{semi destiny} and by \autoref{small reindexing}, the semantic type: $$\tp[\judge{\gamma}{\gamma,a(\gamma),b(\gamma)}](\id_{[\judge{\gamma}{A(\gamma):\type}]})=\tp[\judge{\gamma}{\gamma,a(\gamma),b(\gamma)}]([\judge{\gamma,x,y}{x=y:\type}])$$ in semantic context $[\gamma]$ admits the presentation: $$\judge{\gamma}{a(\gamma)=b(\gamma):\type}.$$ With this presentation, if we are given a semantic term of this semantic type in context $[\gamma]$, i.e. a section $[\gamma]\to[\gamma,p':a(\gamma)=b(\gamma)]$ of the corresponding display map $[\judge{\gamma,p':a(\gamma)=b(\gamma)}{\gamma}]$, then by \autoref{nuovo gioco} it admits a representative of the form $\judge{\gamma}{\gamma:\Gamma,p(\gamma):a(\gamma)=b(\gamma)}$ for some judgement $\judge{\gamma}{p(\gamma):a(\gamma)=b(\gamma)}$. Hence $[\judge{\gamma}{\gamma,a(\gamma)}]=[\judge{\gamma}{\gamma,b(\gamma)}]$, by 1. of \autoref{giocare}. We are left to verify that: $$[\judge{\gamma}{\gamma:\Gamma,p(\gamma):a(\gamma)=b(\gamma)}]=\refl_{[\judge{\gamma}{A(\gamma):\type}]}^{[\judge{\gamma}{\gamma:\Gamma,a(\gamma):A(\gamma)}]}.$$ As the diagram (where $p'':a(\gamma)=a(\gamma)$): \[\begin{tikzcd}
	\gamma &&&& {\gamma,x} \\
	\\
	{\gamma,p''} &&&& {\gamma,x,y,p}
	\arrow["{\judge{\gamma}{\gamma,a(\gamma)}}", from=1-1, to=1-5]
	\arrow["{\judge{\gamma}{\gamma,\refl(a(\gamma))}}"', from=1-1, to=3-1]
	\arrow["{\judge{\gamma,x}{\gamma,x,x,\refl(x)}}", from=1-5, to=3-5]
	\arrow["{\judge{\gamma,p''}{\gamma,a(\gamma),a(\gamma),p''}}", from=3-1, to=3-5]
	\end{tikzcd}\] commutes even judgementally and ${\judge{\gamma,p''}{\gamma,a(\gamma),a(\gamma),p''}}$ represents the extension: $$(\;[\gamma]\xrightarrow{[\judge{\gamma}{\gamma,a(\gamma):A(\gamma)}]}[\gamma,x]\xrightarrow{[\judge{\gamma}{\gamma:\Gamma,a(\gamma):A(\gamma)}].[\judge{\gamma}{A(\gamma):\type}]}[\gamma,x,y]\;).\id_{[\judge{\gamma}{A(\gamma):\type}]}$$ then $\refl_{[\judge{\gamma}{A(\gamma):\type}]}^{[\judge{\gamma}{\gamma:\Gamma,a(\gamma):A(\gamma)}]}$ is represented by $\judge{\gamma}{\gamma,\refl(a(\gamma))}$ and we are left to verify that: $$[\judge{\gamma}{\gamma,p(\gamma):a(\gamma)=b(\gamma)}]=[\judge{\gamma}{\gamma,\refl(a(\gamma))}].$$
    
    Let us consider a canonical context homotopy equivalence $\equivpair{\phi}{\psi}$ between $\gamma$ and itself such that $\judge{\gamma}{\alpha(\gamma):\ci(\gamma)=\gamma}$ and let $\equivpair{\phi}{\psi}$ be a canonical homotopy equivalence between $A(\gamma)$ and itself relative to $\equivpair{\ci}{\di}$ and such that $\judge{\gamma,x:A(\gamma)}{\alpha_1(\gamma,x):\phi(\gamma,x)=\alpha(\gamma)^*x}$ (see \autoref{there are context identities} and \autoref{there are identities}). Now, if: \[\begin{aligned}
    \judge{\gamma&}{\refl_1(\gamma)\equiv \alpha_1(\gamma,a(\gamma))\bullet a(\alpha(\gamma))^{-1}:\phi(\gamma,a(\gamma))=a(\ci(\gamma))}\\
    \judge{\gamma&}{\refl_2(\gamma)\equiv \alpha_1(\gamma,a(\gamma))\bullet a(\alpha(\gamma))^{-1}}\bullet p(\ci(\gamma))
    :\phi(\gamma,a(\gamma))=b(\ci(\gamma))
    \end{aligned}\] then the corresponding $\equivpair{\phi^=}{\psi^=}$ of \autoref{id equiv} happens to be a canonical homotopy equivalence between $a(\gamma)=a(\gamma)$ and $a(\gamma)=b(\gamma)$ relative to $\equivpair{\ci}{\di}$ (see the notion of canonical homotopy equivalence in \autoref{section2.3.VI}, as well as \autoref{crucial example}) and by propositional groupoidality: $$\judge{\gamma}{\phi^=(\gamma,\refl(a(\gamma)))=p(\ci(\gamma))}.$$ Therefore, the diagram: \[\begin{tikzcd}
	\gamma &&& {\gamma,p''} \\
	\\
	\gamma &&& {\gamma,p'}
	\arrow["{\judge{\gamma}{\gamma,\refl(a(\gamma))}}", from=1-1, to=1-4]
	\arrow["{\judge{\gamma}{\gamma,p(\gamma)}}", from=3-1, to=3-4]
	\arrow["\ci"{description}, from=1-1, to=3-1]
	\arrow["{\ci,\phi^=}"{description}, from=1-4, to=3-4]
	\end{tikzcd}\] commutes propositionally and, by the notion of canonical context homotopy equivalence, the morphism of contexts $\ci,\phi^=$ is the first component of a canonical context homotopy equivalence $\gamma,p''\to\gamma,p'$ (see \autoref{context homotopy equivalence expansion}). Therefore $\judge{\gamma}{\gamma,\refl(a(\gamma))}$ and $\judge{\gamma}{\gamma,p(\gamma)}$ represent the same morphism of semantic contexts i.e. the same arrow of $\mathcal{C}$ and we are done. 
 
    \item\textit{Compatibility with the substitution}. If we are given a morphism of contexts: $$[\delta:\Delta]\xrightarrow{[\judge{\delta}{f(\delta):\Gamma}]}[\gamma:\Gamma]$$ (such a representative exists by 3. of \autoref{giocare}) then the semantic context extension $[\judge{\delta}{f(\delta):\Gamma}].[\judge{\gamma}{A(\gamma):\type}]$ admits the presentation: $$[\delta,x':A(f(\delta))]\xrightarrow{[\judge{\delta,x'}{f(\delta),x'}]}[\gamma,x:A(\gamma)]$$ and the further semantic context extension: $$([\judge{\delta}{f(\delta):\Gamma}].[\judge{\gamma}{A(\gamma):\type}]).[\judge{\gamma,x:A(\gamma)}{A(\gamma):\type}]$$ admits the presentation: $$[\delta,x',y':A(f(\delta))]\xrightarrow{[\judge{\delta,x',y'}{f(\delta),x',y'}]}[\gamma,x,y:A(\gamma)].$$ Hence, by \autoref{there are context identities}, by \autoref{semi destiny} and by \autoref{small reindexing}, the type: $$\tp(([\judge{\delta}{f(\delta):\Gamma}].[\judge{\gamma}{A(\gamma):\type}]).[\judge{\gamma,x:A(\gamma)}{A(\gamma):\type}])\id_{[\judge{\gamma}{A(\gamma):\type}]}$$ in context $[\delta,x',y']$ is represented by $\judge{\delta,x',y'}{x'=y'}$, which is a representative of $\id_{[\judge{\delta}{A(f(\delta)):\type}]}$, and $\judge{\delta}{A(f(\delta)):\type}$ represents the type: $$(\tp[\judge{\delta}{f(\delta):\Gamma}])[\judge{\gamma}{A(\gamma):\type}]$$ in context $\delta$, again by \autoref{there are context identities}, by \autoref{semi destiny} and by \autoref{small reindexing}. We conclude that: \begin{center}$\tp(([\judge{\delta}{f(\delta):\Gamma}].[\judge{\gamma}{A(\gamma):\type}]).[\judge{\gamma,x:A(\gamma)}{A(\gamma):\type}])\id_{[\judge{\gamma}{A(\gamma):\type}]}=$\\$=\id_{(\tp[\judge{\delta}{f(\delta):\Gamma}])[\judge{\gamma}{A(\gamma):\type}]}.$\end{center} Now, the semantic context extension: $$(([\judge{\delta}{f(\delta):\Gamma}].[\judge{\gamma}{A(\gamma):\type}]).[\judge{\gamma,x:A(\gamma)}{A(\gamma):\type}]).\id_{[\judge{\gamma}{A(\gamma):\type}]}$$ admits the presentation: $$[\delta,x',y',p':x'=y']\xrightarrow{[\judge{\delta,x',y',p'}{f(\delta),x',y',p'}]}[\gamma,x,y,p:x=y]$$ hence we are left to verify that: \[\begin{tikzcd}
	{[\delta,x']} &&& {[\gamma,x]} \\
	\\
	{[\delta,x',y',p']} &&& {[\gamma,x,y,p]}
	\arrow["{\refl_{[\judge{\delta}{A(f(\delta))}]}}"{description}, from=1-1, to=3-1]
	\arrow["{\refl_{[\judge{\gamma}{A(\gamma)}]}}"{description}, from=1-4, to=3-4]
	\arrow["{[f(\delta),x']}", from=1-1, to=1-4]
	\arrow["{[f(\delta),x',y',p']}", from=3-1, to=3-4] \end{tikzcd}\] commutes in $\mathcal{C}$. By 1. and 4. of \autoref{giocare}, we are done since the diagram of the representatives: \[\begin{tikzcd}
	{\delta,x'} &&& {\gamma,x} \\
	\\
	{\delta,x',y',p'} &&& {\gamma,x,y,p}
	\arrow["{\delta,x',x',\refl(x')}"{description}, from=1-1, to=3-1]
	\arrow["{\gamma,x,x,\refl(x)}"{description}, from=1-4, to=3-4]
	\arrow["{f(\delta),x'}", from=1-1, to=1-4]
	\arrow["{f(\delta),x',y',p'}", from=3-1, to=3-4]
	\arrow[shift right=0.7, shorten <=40pt, shorten >=45pt, no head, from=3-1, to=1-4]
	\arrow[shorten <=40pt, shorten >=45pt, no head, from=3-1, to=1-4]
	\arrow[shift left=0.7, shorten <=40pt, shorten >=45pt, no head, from=3-1, to=1-4]
	\end{tikzcd}\] commutes even judgmentally.\qedhere
\end{itemize}
\endproof
\end{prop}

\subsection{Semantic dependent product and sum types}\label{section2.6.III}

Let $[\judge{\gamma}{A(\gamma):\type}]$ be a semantic type in semantic context $[\gamma:\Gamma]$ and let: $$[\judge{\gamma,x:A(\gamma)}{B(\gamma,x):\type}]$$ be a semantic type in semantic context $[\gamma].[\judge{\gamma}{A(\gamma):\type}]=[\gamma,x]$, where, as usual, we refer to see 1. of \autoref{simpatico} for the presentation we use for the semantic types. Below, we define the choices of the semantic types and the semantic terms that we need so that $\C$ is endowed with semantic dependent products and semantic dependent sums respectively. \begin{itemize}[align=left]
    \item[$\pid$] \textit{Formation}. We define the semantic type $\pid_{[\judge{\gamma}{A(\gamma):\type}]}^{[\judge{\gamma,x}{B(\gamma,x):\type}]}$ in context $[\gamma]$ as the semantic type of $\C$ represented by the type judgement $\judge{\gamma}{\Pi_{x:A(\gamma)}B(\gamma,x):\type}$.

    \textit{Introduction}. If we are given a semantic term: $$[\gamma,x]\xrightarrow{[\judge{\gamma,x}{\gamma,x,b(\gamma,x)}]}[\gamma,x,y:B(\gamma,x)]$$ (see \autoref{nuovo gioco} to justify this presentation of the semantic terms) of semantic type: $$[\judge{\gamma,x:A(\gamma)}{B(\gamma,x):\type}]$$ we define $\lambda [\judge{\gamma,x}{\gamma,x,b(\gamma,x)}]$ to be the semantic term $[\gamma]\to[\gamma,z:\Pi_{x:A(\gamma)}B(\gamma,x)]$ of sematic type $[\judge{\gamma}{\Pi_{x:A(\gamma)}B(\gamma,x):\type}]$ represented by the morphism of contexts: $$\judge{\gamma}{\gamma,\lambda x.b(x):\Pi_{x:A(\gamma)}B(\gamma,x)}$$ of $\ptt$.

    \textit{Elimination}. If we are given semantic terms: $$[\gamma]\xrightarrow{[\judge{\gamma}{\gamma,z(\gamma)}]}[\gamma,z]\text{ \, and \, }[\gamma]\xrightarrow{[\judge{\gamma}{\gamma,a(\gamma)}]}[\gamma,x]$$ of semantic type $[\judge{\gamma}{\Pi_{x:A(\gamma)}B(\gamma,x):\type}]$ and $[\judge{\gamma}{A(\gamma):\type}]$ respectively, then we define the semantic term: $$[\gamma]\xrightarrow{\ev_{[\judge{\gamma}{\gamma,z(\gamma)}]}^{[\judge{\gamma}{\gamma,a(\gamma)}]}}[\gamma,y':B(\gamma,a(\gamma))]$$ of semantic type: $$(\tp[\judge{\gamma}{\gamma,a(\gamma)}])[\judge{\gamma,x:A(\gamma)}{B(\gamma,x):\type}]=[\judge{\gamma}{B(\gamma,a(\gamma))}]$$ as the one represented by the morphism of contexts: $$\judge{\gamma}{\gamma,\ev(z(\gamma),a(\gamma)):B(\gamma,a(\gamma))}$$ of $\ptt$ (the presentation we use for $(\tp[\judge{\gamma}{\gamma,a(\gamma)}])[\judge{\gamma,x:A(\gamma)}{B(\gamma,x):\type}]$ is justified as usual by \autoref{there are context identities}, by \autoref{semi destiny} and by \autoref{small reindexing}).

    One can prove that:

    \begin{prop} With the above choices the properties of \textnormal{compatibility with the substitution}, \textnormal{computation}, and \textnormal{expansion} of \autoref{semantic dependent product types} are satisfied, hence $\C$ is equipped with semantic dependent product types.
    \end{prop}
    
    \item[$\sigmad$] \textit{Formation}. We define the semantic type $\sigmad_{[\judge{\gamma}{A(\gamma):\type}]}^{[\judge{\gamma,x}{B(\gamma,x):\type}]}$ in context $[\gamma]$ as the semantic type of $\C$ represented by the type judgement $\judge{\gamma}{\Sigma_{x:A(\gamma)}B(\gamma,x):\type}$.

    \textit{Introduction}. We define the morphism $p_{[\judge{\gamma}{A(\gamma):\type}]}^{[\judge{\gamma,x}{B(\gamma,x):\type}]}$ making the diagram: \[\begin{tikzcd}
	{[\gamma,x,y]} &&& {[\gamma,u]} \\
	\\
	\\
	{[\gamma,x]} &&& {[\gamma]}
	\arrow["{P_{\sigmad_{[\judge{\gamma}{A(\gamma):\type}]}^{[\judge{\gamma,x}{B(\gamma,x):\type}]}}}"{description}, from=1-4, to=4-4]
	\arrow["{P_{[\judge{\gamma}{A(\gamma):\type}]}}", from=4-1, to=4-4]
	\arrow["{P_{[\judge{\gamma,x}{B(\gamma,x):\type}]}}"{description}, from=1-1, to=4-1]
	\arrow["{p_{[\judge{\gamma}{A(\gamma):\type}]}^{[\judge{\gamma,x}{B(\gamma,x):\type}]}}", from=1-1, to=1-4] \end{tikzcd}\] commute (where $u:\Sigma_{x:A(\gamma)}B(\gamma,x):\type$) as the one represented by the morphism of contexts: $$\judge{\gamma,x,y}{\gamma,\langle x,y\rangle}$$ of $\ptt$.

    \textit{Elimination}. If we are given a semantic type $[\judge{\gamma,u}{C(\gamma,u):\type}]$ in semantic context $[\gamma,u]$ and a semantic term $[\judge{\gamma,x,y}{\gamma,x,y,c(\gamma,x,y):C(\gamma,\langle x,y\rangle)}]$ of semantic type: $$(\tp[\judge{\gamma,x,y}{\gamma,\langle x,y\rangle}])[\judge{\gamma,u}{C(\gamma,u):\type}]=[\judge{\gamma,x,y}{C(\gamma,\langle x,y\rangle)}]$$---as usual, use 1. of \autoref{simpatico} to justify the presentation we use for the semantic types; use \autoref{nuovo gioco} for the one of the semantic terms; use \autoref{there are context identities}, \autoref{semi destiny} and \autoref{small reindexing} for the one of the semantic type $(\tp[\judge{\gamma,x,y}{\gamma,\langle x,y\rangle}])[\judge{\gamma,u}{C(\gamma,u):\type}]$---then we define the semantic term: $$[\gamma,u] \xrightarrow{\splitt_{[\judge{\gamma,x,y}{\gamma,x,y,c(\gamma,x,y):C(\gamma,\langle x,y\rangle)}]}}[\gamma,u,c]$$ of semantic type $[\judge{\gamma,u}{C(\gamma,u):\type}]$ as the one represented by the morphism of contexts: $$\judge{\gamma,u}{\gamma,u,\splitt(c,\gamma,u):C(\gamma,u)}$$ of $\ptt$.
    
    Again, one can prove that: \begin{prop}
        With the above choices the properties of \textnormal{computation} and \textnormal{compatibility with the substitution} of \autoref{semantic dependent sum types} are satisfied, hence $\C$ is equipped with semantic dependent sum types.
    \end{prop}
    
\end{itemize}

\bigskip

We might summarise the results of the current section into the following:

\begin{thm}\label{modelofett}
    The category with attributes $\C$ is a model of $\ett$.
\end{thm}

\subsection{Conservativity result}\label{section2.6.IV}

We remind that $\sptt$ indicates the sub-theory of $\ptt$ generated by the atomic h-sets of $\ptt$. In other words, $\sptt$ is the propositional type theory whose contexts are the h-elementary contexts of $\ptt$, whose type judgements are the h-elementary type judgements of $\ptt$ and whose term judgements are the term judgements of $\ptt$ of an h-elementary type in h-elementary context. We refer to \autoref{section2.2.III} for further details.

We proved that the category with attributes $\C$ is a model of $\ett$ (\autoref{modelofett}). By \autoref{soundness i.e. initiality} and \autoref{extensional is propositional} there are: \begin{itemize}
    \item a unique morphism $\{\cdot\}:\ettm\to\C$ of $\ett$;
    \item a unique morphism $\{\cdot\}':\spttm\to\C$ of $\sptt$.
\end{itemize} The diagram: \[\begin{tikzcd}
	\spttm \\
	&& \C \\
	\ettm
	\arrow["{\{\cdot\}'}", from=1-1, to=2-3]
	\arrow["{\{\cdot\}}"', from=3-1, to=2-3]
	\arrow["{|\cdot|}"', from=1-1, to=3-1]
\end{tikzcd}\] where $|\cdot|$ is the canonical interpretation of \autoref{canonical interpretation}, commutes by \autoref{extensional is propositional} and by \autoref{soundness i.e. initiality}. However, the quotient mapping $[\cdot]$ defined by the equivalence relation of \autoref{section2.5} happens to be a morphism $\spttm\to\C$ of $\sptt$. Since $\{\cdot\}'=[\cdot]$ by \autoref{soundness i.e. initiality}, the diagram: \[\begin{tikzcd}
	\spttm \\
	&& \C \\
	\ettm
	\arrow["{[\cdot]}", from=1-1, to=2-3]
	\arrow["{\{\cdot\}}"', from=3-1, to=2-3]
	\arrow["{|\cdot|}"', from=1-1, to=3-1]
\end{tikzcd}\] commutes. Moreover, one can verify that the morphism of semantic types $\C\to\ettm$ of \autoref{[]=>||morphismofsemantictypes} is a morphism of $\ettm$. By \autoref{soundness i.e. initiality} it must be a retraction of $\{\cdot\}$, hence $\{\cdot\}$ is injective on semantic contexts, morphisms between them, and semantic types.

Now, let $\gamma:\Gamma$ be an h-elementary context of $\ptt$ and let $\judge{\gamma:\Gamma}{T(\gamma):\type}$ be an h-elementary judgement of $\ptt$. Let us suppose that: $$\judgext{|\gamma:\Gamma|}{t(|\gamma:\Gamma|):|T(\gamma)|}$$ in $\ett$. Here we recall that $|T(\gamma)|$ denotes the type of $\ett$ in context $|\gamma|$ constituting the type judgement $|\judge{\gamma}{T(\gamma):\type}|$ in $\ett$ (remind that $|\cdot|$ needs to map semantic types in semantic context $\gamma$---i.e. h-elementary type judgements of $\ptt$ in context $\gamma$---to semantic types in semantic context $|\gamma|$---i.e. type judgements of $\ett$ in context $|\gamma|$---), see \autoref{canonical interpretation} for more details. Hence, in $\ettm$ the judgement $\judgext{|\gamma|}{|\gamma|,t(|\gamma|):|T(\gamma)|}$ is a morphism of semantic contexts: $$|\gamma|\to|\gamma|.|\judge{\gamma}{T(\gamma):\type}|$$ happening to be a section of the display map $|\gamma|.|\judge{\gamma}{T(\gamma):\type}|\xrightarrow{P_{|\judge{\gamma}{T(\gamma):\type}|}}|\gamma|$. Therefore, being $\{\cdot\}$ a morphism of $\ett$, we get a section: $$\begin{aligned}\relax [\gamma]=\{|\gamma|\}\xrightarrow{\{\judgext{|\gamma|}{|\gamma|,t(|\gamma|):|T(\gamma)|}\}}\{|\gamma|.|\judge{\gamma}{T(\gamma):\type}|\}&=\{|\gamma|\}.\{|\judge{\gamma}{T(\gamma):\type}|\}\\
&=[\gamma].[\judge{\gamma}{T(\gamma):\type}]\\
&=[\gamma:\Gamma,t:T(\gamma)]\end{aligned}$$ of: $$[\gamma:\Gamma,t:T(\gamma)]\xrightarrow{\;\;P_{\{|\judge{\gamma}{T(\gamma):\type}|\}}=P_{[\judge{\gamma}{T(\gamma):\type}]}\;\;}[\gamma]$$ in $\C$, where the last equality follows by \autoref{there are context identities} and \autoref{semi destiny}. By \autoref{nuovo gioco}, the arrow $\{\judgext{|\gamma:\Gamma|}{|\gamma|,t(|\gamma|):|T(\gamma)|}\}$ happens to be represented by a morphism of contexts: $$\judge{\gamma:\Gamma}{\gamma:\Gamma,\tilde{t}(\gamma):T(\gamma)}$$ for  some judgement $\judge{\gamma:\Gamma}{\tilde{t}(\gamma):T(\gamma)}$. In particular: $$\begin{aligned}\{\judgext{|\gamma|}{|\gamma|,|\tilde{t}(\gamma)|:|T(\gamma)|}\}&=\{|\judge{\gamma}{\gamma,\;\tilde{t}(\gamma):T(\gamma)}|\}\\&=[\judge{\gamma}{\gamma,\;\tilde{t}(\gamma):T(\gamma)}]\\&=\{\judgext{|\gamma|}{|\gamma|,t(|\gamma|):|T(\gamma)|}\}\end{aligned}$$ and, by injectivity of $\{\cdot\}$ on morphisms of semantic contexts, we conclude that: $$\judgext{|\gamma|}{|\tilde{t}(\gamma)|\equiv t(|\gamma|)}.$$ Hence we have just proven the following:

\begin{thm}[Conservativity]\label{main result}
Let $\gamma:\Gamma$ be an h-elementary context of $\ptt$ and let $\judge{\gamma:\Gamma}{T(\gamma:\Gamma):\type}$ be an h-elementary type judgement of $\ptt$. Whenever $\ett$ infers $\judgext{|\gamma:\Gamma|}{t(|\gamma:\Gamma|):|T(\gamma)|}$, then $\ptt$ infers $\judge{\gamma:\Gamma}{\tilde{t}(\gamma:\Gamma):T(\gamma:\Gamma)}$ and $\ett$ infers $\judgext{|\gamma:\Gamma|}{|\tilde{t}(\gamma:\Gamma)|\equiv t(|\gamma:\Gamma|)}$.
\end{thm}

In particular, we infer that:

\begin{thm}
The theory $\ett$ is conservative over the theory $\ptt + \uip$.
\end{thm}

\newpage

\section{Inference rules for dependent type theories}
\label{inferencerules}

In this section we enumerate the inference rules that we consider in this paper, particularly for what we call \textit{extensional}, \textit{intensional}, and \textit{propositional} type theories with three type constructors: identities, dependent products and dependent sums. We omit the usual structural rules of a strict dependent type theory---that every theory considered in this paper is assumed to satisfy. For an enumeration of the structural rules, see \cite[Chapter III---context formation, context equality, judgement formation, judgement equality]{MR1134134}, \cite{MR1674451}, \cite{Hofmann1997}.

\bigskip

\subsection{Identity types}\label{sectionA.1}

\medskip

\text{ }

\begin{figure}[b]

\bigskip
\bigskip

\[\begin{aligned}
\text{Formation \quad}&\inferrule {\quad\quad\quad\quad\quad\judge{\textcolor{green}{\gamma : \Gamma}}{A\textcolor{green}{(\gamma)} : \type}\quad\quad\quad\quad\quad}{\judge{\textcolor{green}{\gamma:\Gamma;\;}x,x' : A\textcolor{green}{(\gamma)}}{x=x' : \type}}\\
\\
\\
\\
\text{Introduction \quad}&\inferrule {\quad\quad\quad\quad\quad\judge{\textcolor{green}{\gamma : \Gamma}}{A\textcolor{green}{(\gamma)} : \type}\quad\quad\quad\quad\quad}{\judge{\textcolor{green}{\gamma:\Gamma;\;}x : A\textcolor{green}{(\gamma)}}{\refl(x) : x=x}}\\
\\
\\
\\
\text{Extensionality \quad}&\inferrule{\quad\quad\quad\quad\quad\judge{\textcolor{green}{\gamma : \Gamma}}{A\textcolor{green}{(\gamma)} : \type}\quad\quad\quad\quad\quad}{\judge{\textcolor{green}{\gamma:\Gamma;\; }x, x' : A\textcolor{green}{(\gamma)};\; p : x=x'}{x\equiv x'}}\\
\\
\\
\\
\text{Identity proof irrelevance \quad}&\inferrule{\quad\quad\quad\quad\quad\judge{\textcolor{green}{\gamma : \Gamma}}{A\textcolor{green}{(\gamma)} : \type}\quad\quad\quad\quad\quad}{\judge{\textcolor{green}{\gamma:\Gamma;\;} x, x' : A\textcolor{green}{(\gamma)};\; p : x=x'}{p\equiv \refl(x)} }\\
\end{aligned}\]

\bigskip
\bigskip
\bigskip
\bigskip

\caption{\;\;\;Extensional identity types}
\label{figureA.1.1}\label{extensional id}

\end{figure}

\text{ }

\newpage

\text{ }

\begin{figure}[b]

\bigskip
\bigskip

\[\begin{aligned}
\text{Formation \quad}&\inferrule{\quad\quad\quad\quad\quad\quad\quad\quad\quad\judge{\textcolor{green}{\gamma : \Gamma}}{A\textcolor{green}{(\gamma)} : \type}\quad\quad\quad\quad\quad\quad\quad\quad\quad}{ \\ \judge{\textcolor{green}{\gamma : \Gamma;\;}x,x' : A\textcolor{green}{(\gamma)}}{x=x' : \type} \\ }\\
\\
\\
\\
\\
\text{Introduction \quad}&\inferrule{\quad\quad\quad\quad\quad\quad\quad\quad\quad\judge{\textcolor{green}{\gamma : \Gamma}}{A\textcolor{green}{(\gamma)} : \type}\quad\quad\quad\quad\quad\quad\quad\quad\quad}{ \\ \judge{\textcolor{green}{\gamma : \Gamma;\;}x : A\textcolor{green}{(\gamma)}}{\refl(x) : x=x} \\ }\\
\\
\\
\\
\\
\text{Elimination \quad}&\inferrule{\quad\quad\quad\quad\quad\quad\quad\quad\quad\judge{\textcolor{green}{\gamma : \Gamma}}{A\textcolor{green}{(\gamma)} : \type}\quad\quad\quad\quad\quad\quad\quad\quad\quad \\\\ \\\\ \judge{\textcolor{green}{\gamma : \Gamma;\;} x, x' : A\textcolor{green}{(\gamma)};\; p : x=x'}{C(\textcolor{green}{\gamma,\;}x,x',p) : \type} \\\\ \\\\ \judge{\textcolor{green}{\gamma : \Gamma;\;} x : A\textcolor{green}{(\gamma)}}{q(\textcolor{green}{\gamma,\;}x) : C(\textcolor{green}{\gamma,\;}x,x,\refl(x))}}{\text{ } \\\\ \judge{\textcolor{green}{\gamma : \Gamma;\;}x, x' : A\textcolor{green}{(\gamma)};\; p : x=x'}{\j(q,\textcolor{green}{\gamma,\;}x,x',p) : C(\textcolor{green}{\gamma,\;}x,x',p)}}\\
\\
\\
\\
\\
\text{Computation \quad}&\inferrule{\quad\quad\quad\quad\quad\quad\quad\quad\quad\judge{\textcolor{green}{\gamma : \Gamma}}{A\textcolor{green}{(\gamma)} : \type}\quad\quad\quad\quad\quad\quad\quad\quad\quad \\\\ \\\\ \judge{\textcolor{green}{\gamma : \Gamma;\;} x, x' : A\textcolor{green}{(\gamma)};\; p : x=x'}{C(\textcolor{green}{\gamma,\;}x,x',p) : \type} \\\\ \\\\ \judge{\textcolor{green}{\gamma : \Gamma;\;} x : A\textcolor{green}{(\gamma)}}{q(\textcolor{green}{\gamma,\;}x) : C(\textcolor{green}{\gamma,\;}x,x,\refl(x))}}{\text{ } \\\\ \judge{\textcolor{green}{\gamma : \Gamma;\;}x : A\textcolor{green}{(\gamma)}}{\j(q,\textcolor{green}{\gamma,\;}x,x,\refl(x))\equiv q(\textcolor{green}{\gamma,\;}x)}}\\
\end{aligned}\]

\bigskip
\bigskip
\bigskip
\bigskip

\caption{\;\;\;Intensional identity types}
\label{figureA.1.2}\label{intensional id}

\end{figure}

\text{ }

\newpage

\text{ }

\begin{figure}[b]

\bigskip
\bigskip

\[\begin{aligned}
\text{Formation \quad}&\inferrule{\quad\quad\quad\quad\quad\quad\quad\quad\quad\judge{\textcolor{green}{\gamma : \Gamma}}{A\textcolor{green}{(\gamma)} : \type}\quad\quad\quad\quad\quad\quad\quad\quad\quad}{ \\ \judge{\textcolor{green}{\gamma : \Gamma;\;}x,x' : A\textcolor{green}{(\gamma)}}{x=x' : \type} \\ }\\
\\
\\
\\
\\
\text{Introduction \quad}&\inferrule{\quad\quad\quad\quad\quad\quad\quad\quad\quad\judge{\textcolor{green}{\gamma : \Gamma}}{A\textcolor{green}{(\gamma)} : \type}\quad\quad\quad\quad\quad\quad\quad\quad\quad}{ \\ \judge{\textcolor{green}{\gamma : \Gamma;\;}x : A\textcolor{green}{(\gamma)}}{\refl(x) : x=x} \\ }\\
\\
\\
\\
\\
\text{Elimination \quad}&\inferrule{\quad\quad\quad\quad\quad\quad\quad\quad\quad\judge{\textcolor{green}{\gamma : \Gamma}}{A\textcolor{green}{(\gamma)} : \type}\quad\quad\quad\quad\quad\quad\quad\quad\quad \\\\ \\\\ \judge{\textcolor{green}{\gamma : \Gamma;\;} x, x' : A\textcolor{green}{(\gamma)};\; p : x=x'}{C(\textcolor{green}{\gamma,\;}x,x',p) : \type} \\\\ \\\\ \judge{\textcolor{green}{\gamma : \Gamma;\;} x : A\textcolor{green}{(\gamma)}}{q(\textcolor{green}{\gamma,\;}x) : C(\textcolor{green}{\gamma,\;}x,x,\refl(x))}}{\text{ } \\\\ \judge{\textcolor{green}{\gamma : \Gamma;\;}x, x' : A\textcolor{green}{(\gamma)};\; p : x=x'}{\j(q,\textcolor{green}{\gamma,\;}x,x',p) : C(\textcolor{green}{\gamma,\;}x,x',p)}}\\
\\
\\
\\
\\
\text{Prop computation \quad}&\inferrule{\quad\quad\quad\quad\quad\quad\quad\quad\quad\judge{\textcolor{green}{\gamma : \Gamma}}{A\textcolor{green}{(\gamma)} : \type}\quad\quad\quad\quad\quad\quad\quad\quad\quad \\\\ \\\\ \judge{\textcolor{green}{\gamma : \Gamma;\;} x, x' : A\textcolor{green}{(\gamma)};\; p : x=x'}{C(\textcolor{green}{\gamma,\;}x,x',p) : \type} \\\\ \\\\ \judge{\textcolor{green}{\gamma : \Gamma;\;} x : A\textcolor{green}{(\gamma)}}{q(\textcolor{green}{\gamma,\;}x) : C(\textcolor{green}{\gamma,\;}x,x,\refl(x))}}{\text{ } \\\\ \judge{\textcolor{green}{\gamma : \Gamma;\;}x : A\textcolor{green}{(\gamma)}}{\h(q,\textcolor{green}{\gamma,\;}x):\j(q,\textcolor{green}{\gamma,\;}x,x,\refl(x))= q(\textcolor{green}{\gamma,\;}x)}}\\
\end{aligned}\]

\bigskip
\bigskip
\bigskip
\bigskip

\caption{\;\;\;Propositional identity types}
\label{figureA.1.3}\label{propositional id}

\end{figure}

\text{ }

\newpage

\subsection{Dependent product types}\label{sectionA.2}

\medskip

\text{ }

\begin{figure}[b]

\bigskip
\bigskip

\[\begin{aligned}
\text{Formation \quad}&\inferrule{\judge{\textcolor{green}{\gamma : \Gamma}}{A\textcolor{green}{(\gamma)} : \type} \\ \judge{\textcolor{green}{\gamma : \Gamma;\;}x:A\textcolor{green}{(\gamma)}}{B(\textcolor{green}{\gamma,\;}x):\type}}{\judge{\textcolor{green}{\gamma : \Gamma}}{\Pi_{x:A\textcolor{green}{(\gamma)}}B(\textcolor{green}{\gamma,\;}x) : \type}}\\
\\
\\
\\
\text{Elimination \quad}&\inferrule{\judge{\textcolor{green}{\gamma : \Gamma}}{A\textcolor{green}{(\gamma)} : \type} \\ \judge{\textcolor{green}{\gamma : \Gamma;\;}x:A\textcolor{green}{(\gamma)}}{B(\textcolor{green}{\gamma,\;}x):\type}}{\judge{\textcolor{green}{\gamma : \Gamma;\;}z:\Pi_{x:A\textcolor{green}{(\gamma)}}B(\textcolor{green}{\gamma,\;}x);\;x:A\textcolor{green}{(\gamma)}}{\ev(z,x):B(\textcolor{green}{\gamma,\;}x)}}\\
\\
\\
\\
\text{Introduction \quad}&\inferrule{\judge{\textcolor{green}{\gamma : \Gamma}}{A\textcolor{green}{(\gamma)} : \type} \\ \judge{\textcolor{green}{\gamma : \Gamma;\;}x:A\textcolor{green}{(\gamma)}}{B(\textcolor{green}{\gamma,\;}x):\type} \\\\ \judge{\textcolor{green}{\gamma : \Gamma;\;}x:A\textcolor{green}{(\gamma)}}{y(\textcolor{green}{\gamma,\;}x):B(\textcolor{green}{\gamma,\;}x)}}{\judge{\textcolor{green}{\gamma : \Gamma}}{\lambda x. y(\textcolor{green}{\gamma,\;}x):\Pi_{x:A\textcolor{green}{(\gamma)}}B(\textcolor{green}{\gamma,\;}x)}}\\
\\
\\
\\
\text{Computation \quad}&\inferrule{\judge{\textcolor{green}{\gamma : \Gamma}}{A\textcolor{green}{(\gamma)} : \type} \\ \judge{\textcolor{green}{\gamma : \Gamma;\;}x:A\textcolor{green}{(\gamma)}}{B(\textcolor{green}{\gamma,\;}x):\type} \\\\ \judge{\textcolor{green}{\gamma : \Gamma;\;}x:A\textcolor{green}{(\gamma)}}{y(\textcolor{green}{\gamma,\;}x):B(\textcolor{green}{\gamma,\;}x)}}{\judge{\textcolor{green}{\gamma : \Gamma}}{\ev(\lambda x. y(\textcolor{green}{\gamma,\;}x),x)\equiv y(\textcolor{green}{\gamma,\;}x)}}\\
\end{aligned}\]

\bigskip
\bigskip
\bigskip
\bigskip

\caption{\;\;\;Dependent product types}
\label{dependent pi}
\label{figureA.2.1}

\end{figure}

\text{ }

\newpage

\text{ }

\begin{figure}[b]

\bigskip
\bigskip

\[\begin{aligned}
\text{Formation \quad}&\inferrule{\quad\judge{\textcolor{green}{\gamma : \Gamma}}{A\textcolor{green}{(\gamma)} : \type} \\ \judge{\textcolor{green}{\gamma : \Gamma;\;}x:A\textcolor{green}{(\gamma)}}{B(\textcolor{green}{\gamma,\;}x):\type}\quad}{\judge{\textcolor{green}{\gamma : \Gamma}}{\Pi_{x:A\textcolor{green}{(\gamma)}}B(\textcolor{green}{\gamma,\;}x) : \type}}\\
\\
\text{Elimination \quad}&\inferrule{\quad\judge{\textcolor{green}{\gamma : \Gamma}}{A\textcolor{green}{(\gamma)} : \type} \\ \judge{\textcolor{green}{\gamma : \Gamma;\;}x:A\textcolor{green}{(\gamma)}}{B(\textcolor{green}{\gamma,\;}x):\type}\quad}{\judge{\textcolor{green}{\gamma : \Gamma;\;}z:\Pi_{x:A\textcolor{green}{(\gamma)}}B(\textcolor{green}{\gamma,\;}x);\;x:A\textcolor{green}{(\gamma)}}{\ev(z,x):B(\textcolor{green}{\gamma,\;}x)}}\\
\\
\text{Introduction \quad}&\inferrule{\quad\judge{\textcolor{green}{\gamma : \Gamma}}{A\textcolor{green}{(\gamma)} : \type} \\ \judge{\textcolor{green}{\gamma : \Gamma;\;}x:A\textcolor{green}{(\gamma)}}{B(\textcolor{green}{\gamma,\;}x):\type}\quad\\\\ \judge{\textcolor{green}{\gamma : \Gamma;\;}x:A\textcolor{green}{(\gamma)}}{y(\textcolor{green}{\gamma,\;}x):B(\textcolor{green}{\gamma,\;}x)}}{\judge{\textcolor{green}{\gamma : \Gamma}}{\lambda x. y(\textcolor{green}{\gamma,\;}x):\Pi_{x:A\textcolor{green}{(\gamma)}}B(\textcolor{green}{\gamma,\;}x)}}\\
\\
\text{Prop computation \quad}&\inferrule{\quad\judge{\textcolor{green}{\gamma : \Gamma}}{A\textcolor{green}{(\gamma)} : \type} \\ \judge{\textcolor{green}{\gamma : \Gamma;\;}x:A\textcolor{green}{(\gamma)}}{B(\textcolor{green}{\gamma,\;}x):\type}\quad\\\\ \judge{\textcolor{green}{\gamma : \Gamma;\;}x:A\textcolor{green}{(\gamma)}}{y(\textcolor{green}{\gamma,\;}x):B(\textcolor{green}{\gamma,\;}x)}}{\judge{\textcolor{green}{\gamma : \Gamma}}{\beta^{\Pi}(y,\textcolor{green}{\gamma,\;}x):\ev(\lambda x. y(\textcolor{green}{\gamma,\;}x),x)= y(\textcolor{green}{\gamma,\;}x)}}\\
\\
\\
\\
\text{Introduction \quad}&\inferrule{\quad\judge{\textcolor{green}{\gamma : \Gamma}}{A\textcolor{green}{(\gamma)} : \type} \\ \judge{\textcolor{green}{\gamma : \Gamma;\;}x:A\textcolor{green}{(\gamma)}}{B(\textcolor{green}{\gamma,\;}x):\type}\quad}{\judge{\textcolor{green}{\gamma : \Gamma;\;}\;z,z':\Pi_{x:A\textcolor{green}{(\gamma)}}B(\textcolor{green}{\gamma,\;}x);\;p:\Pi_{x:A\textcolor{green}{(\gamma)}}\ev(z,x)=\ev(z',x)}\\\\{\funext(\textcolor{green}{\gamma,\;}z,z',p):z=z'}}\\
\\
\text{Prop computation \quad}&\inferrule{\quad\judge{\textcolor{green}{\gamma : \Gamma}}{A\textcolor{green}{(\gamma)} : \type} \\ \judge{\textcolor{green}{\gamma : \Gamma;\;}x:A\textcolor{green}{(\gamma)}}{B(\textcolor{green}{\gamma,\;}x):\type}\quad}{\judge{\textcolor{green}{\gamma : \Gamma;\;}\;z,z':\Pi_{x:A\textcolor{green}{(\gamma)}}B(\textcolor{green}{\gamma,\;}x);\;p:\Pi_{x:A\textcolor{green}{(\gamma)}}\ev(z,x)=\ev(z',x)}\\\\{\beta_\funext(\textcolor{green}{\gamma,\;}z,z',p):\lambda x.\ev(\funext(\textcolor{green}{\gamma,\;}z,z',p),x)=p}}\\
\\
\text{Prop expansion \quad}&\inferrule{\quad\judge{\textcolor{green}{\gamma : \Gamma}}{A\textcolor{green}{(\gamma)} : \type} \\ \judge{\textcolor{green}{\gamma : \Gamma;\;}x:A\textcolor{green}{(\gamma)}}{B(\textcolor{green}{\gamma,\;}x):\type}\quad}{\judge{\textcolor{green}{\gamma : \Gamma;\;}\;z,z':\Pi_{x:A\textcolor{green}{(\gamma)}}B(\textcolor{green}{\gamma,\;}x);\;q:z=z'}\\\\{\eta_\funext(\textcolor{green}{\gamma,\;}z,z',q):q=\funext(\textcolor{green}{\gamma,\;}z,z',\lambda x.\ev(q,x))}}\\
\end{aligned}\]

\bigskip
\bigskip
\bigskip
\bigskip

\caption{\;\;\;Propositional dependent product types}
\label{propositional pi with propositional function extensionality}\label{propositional pi}
\label{figureA.2.2}

\end{figure}

\text{ }

\newpage

\subsection{Dependent sum types}\label{sectionA.3}

\medskip

\text{ }

\begin{figure}[b]

\bigskip
\bigskip

\[\begin{aligned}
\text{Formation \quad}&\inferrule{\judge{\textcolor{green}{\gamma : \Gamma}}{A\textcolor{green}{(\gamma)} : \type} \\ \judge{\textcolor{green}{\gamma : \Gamma;\;}x:A\textcolor{green}{(\gamma)}}{B(\textcolor{green}{\gamma,\;}x):\type}}{\judge{\textcolor{green}{\gamma : \Gamma}}{\Sigma_{x:A\textcolor{green}{(\gamma)}}B(\textcolor{green}{\gamma,\;}x) : \type}}\\
\\
\\
\\
\text{Introduction \quad}&\inferrule{\judge{\textcolor{green}{\gamma : \Gamma}}{A\textcolor{green}{(\gamma)} : \type} \\ \judge{\textcolor{green}{\gamma : \Gamma;\;}x:A\textcolor{green}{(\gamma)}}{B(\textcolor{green}{\gamma,\;}x):\type}}{\judge{\textcolor{green}{\gamma : \Gamma;\;}x:A\textcolor{green}{(\gamma)},y:B(\textcolor{green}{\gamma,\;}x)}{\langle x,y\rangle:\Sigma_{x:A\textcolor{green}{(\gamma)}}B(\textcolor{green}{\gamma,\;}x)}}\\
\\
\\
\\
\\
\text{Elimination \quad}&\inferrule{\judge{\textcolor{green}{\gamma : \Gamma}}{A\textcolor{green}{(\gamma)} : \type} \\ \judge{\textcolor{green}{\gamma : \Gamma;\;}x:A\textcolor{green}{(\gamma)}}{B(\textcolor{green}{\gamma,\;}x):\type} \\\\ \\\\ \judge{\textcolor{green}{\gamma : \Gamma;\;}u:\Sigma_{x:A\textcolor{green}{(\gamma)}}B(\textcolor{green}{\gamma,\;}x)}{C(\textcolor{green}{\gamma,\;}u):\type} \\\\ \\\\ \judge{\textcolor{green}{\gamma : \Gamma;\;}x:A\textcolor{green}{(\gamma)};\;y:B(\textcolor{green}{\gamma,\;}x)}{c(\textcolor{green}{\gamma,\;}x,y):C(\textcolor{green}{\gamma,\;}\langle x,y\rangle)}}{\text{ } \\\\ \judge{\textcolor{green}{\gamma : \Gamma;\;}u:\Sigma_{x:A\textcolor{green}{(\gamma)}}B(\textcolor{green}{\gamma,\;}x)}{\splitt(c,\textcolor{green}{\gamma,\;}u):C(\textcolor{green}{\gamma,\;}u)}}\\
\\
\\
\\
\\
\text{Computation \quad}&\inferrule{\judge{\textcolor{green}{\gamma : \Gamma}}{A\textcolor{green}{(\gamma)} : \type} \\ \judge{\textcolor{green}{\gamma : \Gamma;\;}x:A\textcolor{green}{(\gamma)}}{B(\textcolor{green}{\gamma,\;}x):\type} \\\\ \\\\ \judge{\textcolor{green}{\gamma : \Gamma;\;}u:\Sigma_{x:A\textcolor{green}{(\gamma)}}B(\textcolor{green}{\gamma,\;}x)}{C(\textcolor{green}{\gamma,\;}u):\type} \\\\ \\\\ \judge{\textcolor{green}{\gamma : \Gamma;\;}x:A\textcolor{green}{(\gamma)};\;y:B(\textcolor{green}{\gamma,\;}x)}{c(\textcolor{green}{\gamma,\;}x,y):C(\textcolor{green}{\gamma,\;}\langle x,y\rangle)}}{\text{ } \\\\ \judge{\textcolor{green}{\gamma : \Gamma;\;}x:A\textcolor{green}{(\gamma)};\;y:B(\textcolor{green}{\gamma,\;}x)}{\splitt(c,\textcolor{green}{\gamma,\;}\langle x,y\rangle)\equiv c(\textcolor{green}{\gamma,\;}x,y)}}\\
\end{aligned}\]

\bigskip
\bigskip
\bigskip
\bigskip

\caption{\;\;\;Dependent sum types}
\label{figureA.3.1}\label{dependent sigma}

\end{figure}

\text{ }

\newpage

\text{ }

\begin{figure}[b]

\bigskip
\bigskip

\[\begin{aligned}
\text{Formation \quad}&\inferrule{\judge{\textcolor{green}{\gamma : \Gamma}}{A\textcolor{green}{(\gamma)} : \type} \\ \judge{\textcolor{green}{\gamma : \Gamma;\;}x:A\textcolor{green}{(\gamma)}}{B(\textcolor{green}{\gamma,\;}x):\type}}{\judge{\textcolor{green}{\gamma : \Gamma}}{\Sigma_{x:A\textcolor{green}{(\gamma)}}B(\textcolor{green}{\gamma,\;}x) : \type}}\\
\\
\\
\\
\\
\text{Introduction \quad}&\inferrule{\judge{\textcolor{green}{\gamma : \Gamma}}{A\textcolor{green}{(\gamma)} : \type} \\ \judge{\textcolor{green}{\gamma : \Gamma;\;}x:A\textcolor{green}{(\gamma)}}{B(\textcolor{green}{\gamma,\;}x):\type}}{\judge{\textcolor{green}{\gamma : \Gamma;\;}x:A\textcolor{green}{(\gamma)},y:B(\textcolor{green}{\gamma,\;}x)}{\langle x,y\rangle:\Sigma_{x:A\textcolor{green}{(\gamma)}}B(\textcolor{green}{\gamma,\;}x)}}\\
\\
\\
\\
\\
\text{Elimination \quad}&\inferrule{\judge{\textcolor{green}{\gamma : \Gamma}}{A\textcolor{green}{(\gamma)} : \type} \\ \judge{\textcolor{green}{\gamma : \Gamma;\;}x:A\textcolor{green}{(\gamma)}}{B(\textcolor{green}{\gamma,\;}x):\type} \\\\ \\\\ \judge{\textcolor{green}{\gamma : \Gamma;\;}u:\Sigma_{x:A\textcolor{green}{(\gamma)}}B(\textcolor{green}{\gamma,\;}x)}{C(\textcolor{green}{\gamma,\;}u):\type} \\\\ \\\\ \judge{\textcolor{green}{\gamma : \Gamma;\;}x:A\textcolor{green}{(\gamma)};\;y:B(\textcolor{green}{\gamma,\;}x)}{c(\textcolor{green}{\gamma,\;}x,y):C(\textcolor{green}{\gamma,\;}\langle x,y\rangle)}}{\text{ } \\\\ \judge{\textcolor{green}{\gamma : \Gamma;\;}u:\Sigma_{x:A\textcolor{green}{(\gamma)}}B(\textcolor{green}{\gamma,\;}x)}{\splitt(c,\textcolor{green}{\gamma,\;}u):C(\textcolor{green}{\gamma,\;}u)}}\\
\\
\\
\\
\\
\text{Prop computation \quad}&\inferrule{\judge{\textcolor{green}{\gamma : \Gamma}}{A\textcolor{green}{(\gamma)} : \type} \\ \judge{\textcolor{green}{\gamma : \Gamma;\;}x:A\textcolor{green}{(\gamma)}}{B(\textcolor{green}{\gamma,\;}x):\type} \\\\ \\\\ \judge{\textcolor{green}{\gamma : \Gamma;\;}u:\Sigma_{x:A\textcolor{green}{(\gamma)}}B(\textcolor{green}{\gamma,\;}x)}{C(\textcolor{green}{\gamma,\;}u):\type} \\\\ \\\\ \judge{\textcolor{green}{\gamma : \Gamma;\;}x:A\textcolor{green}{(\gamma)};\;y:B(\textcolor{green}{\gamma,\;}x)}{c(\textcolor{green}{\gamma,\;}x,y):C(\textcolor{green}{\gamma,\;}\langle x,y\rangle)}}{\text{ } \\\\ \judge{\textcolor{green}{\gamma : \Gamma;\;}x:A\textcolor{green}{(\gamma)};\;y:B(\textcolor{green}{\gamma,\;}x)}\\\\{\sigma(c,\textcolor{green}{\gamma,\;}x,y):\splitt(c,\textcolor{green}{\gamma,\;}\langle x,y\rangle)= c(\textcolor{green}{\gamma,\;}x,y)}}\\
\end{aligned}\]

\bigskip
\bigskip
\bigskip
\bigskip

\caption{\;\;\;Propositional dependent sum types}
\label{figureA.3.2}\label{propositional sigma}

\end{figure}

\text{ }

\newpage

\noindent In the next section we recall a characterisations of dependent sum types and propositional sum types \textit{inside} $\ett$ and $\ptt$, respectively. In fact we recall that:
\begin{itemize}
    \item a dependent type theory with extensional identity types has the rules of \autoref{dependent sigma} if and only if it has the rules of \autoref{dependent sigma char};
    \item a dependent type theory with propositional identity types has the rules of \autoref{propositional sigma} if and only if it has the rules of \autoref{propositional sigma char}.
\end{itemize} For further details, we refer the reader to \cite{MR1674451, MR0727078} and to \cite{spadettopropositional}, respectively.

\text{ }

\newpage

\subsection{Characterisation of dependent sum types}\label{sectionA.4}

\medskip

\text{ }

\begin{figure}[b]

\bigskip
\bigskip

\[\begin{aligned}
\text{Formation \quad}&\inferrule{\judge{\textcolor{green}{\gamma : \Gamma}}{A\textcolor{green}{(\gamma)} : \type} \\ \judge{\textcolor{green}{\gamma : \Gamma;\;}x:A\textcolor{green}{(\gamma)}}{B(\textcolor{green}{\gamma,\;}x):\type}}{\judge{\textcolor{green}{\gamma : \Gamma}}{\Sigma_{x:A\textcolor{green}{(\gamma)}}B(\textcolor{green}{\gamma,\;}x) : \type}}\\
\\
\\
\\
\\
\text{Introduction \quad}&\inferrule{\judge{\textcolor{green}{\gamma : \Gamma}}{A\textcolor{green}{(\gamma)} : \type} \\ \judge{\textcolor{green}{\gamma : \Gamma;\;}x:A\textcolor{green}{(\gamma)}}{B(\textcolor{green}{\gamma,\;}x):\type}}{\judge{\textcolor{green}{\gamma : \Gamma;\;}x:A\textcolor{green}{(\gamma)},y:B(\textcolor{green}{\gamma,\;}x)}{\langle x,y\rangle:\Sigma_{x:A\textcolor{green}{(\gamma)}}B(\textcolor{green}{\gamma,\;}x)}}\\
\\
\\
\\
\\
\text{Projection \quad}&\inferrule{\judge{\textcolor{green}{\gamma : \Gamma}}{A\textcolor{green}{(\gamma)} : \type} \\ \judge{\textcolor{green}{\gamma : \Gamma;\;}x:A\textcolor{green}{(\gamma)}}{B(\textcolor{green}{\gamma,\;}x):\type}}{\judge{\textcolor{green}{\gamma : \Gamma;\;}u:\Sigma_{x:A\textcolor{green}{(\gamma)}}B(\textcolor{green}{\gamma,\;}x)}{\pi_1u:A\textcolor{green}{(\gamma)}} \\\\ \\\\ \judge{\textcolor{green}{\gamma : \Gamma;\;}u:\Sigma_{x:A\textcolor{green}{(\gamma)}}B(\textcolor{green}{\gamma,\;}x)}{\pi_2u:B(\textcolor{green}{\gamma,\;}\pi_1(u))}}\\
\\
\\
\\
\\
\text{$\beta$-reduction \quad}&\inferrule{\judge{\textcolor{green}{\gamma : \Gamma}}{A\textcolor{green}{(\gamma)} : \type} \\ \judge{\textcolor{green}{\gamma : \Gamma;\;}x:A\textcolor{green}{(\gamma)}}{B(\textcolor{green}{\gamma,\;}x):\type}}{\judge{\textcolor{green}{\gamma : \Gamma;\;}x:A\textcolor{green}{(\gamma)};\;y:B(\textcolor{green}{\gamma,\;}x)}{\pi_1\langle x,y\rangle\equiv x}  \\\\ \\\\ \judge{\textcolor{green}{\gamma : \Gamma;\;}x:A\textcolor{green}{(\gamma)};\;y:B(\textcolor{green}{\gamma,\;}x)}{\pi_2\langle x,y\rangle\equiv y}}\\
\\
\\
\\
\\
\text{$\eta$-expansion \quad}&\inferrule{\judge{\textcolor{green}{\gamma : \Gamma}}{A\textcolor{green}{(\gamma)} : \type} \\ \judge{\textcolor{green}{\gamma : \Gamma;\;}x:A\textcolor{green}{(\gamma)}}{B(\textcolor{green}{\gamma,\;}x):\type}}{\judge{\textcolor{green}{\gamma : \Gamma;\;}u:\Sigma_{x:A\textcolor{green}{(\gamma)}}B(\textcolor{green}{\gamma,\;}x)}{u\equiv \langle \pi_1u,\pi_2u\rangle}}\\
\end{aligned}\]

\bigskip
\bigskip
\bigskip
\bigskip

\caption{\;\;\;Dependent sum types---negative form}
\label{figureA.4.1}\label{dependent sigma char}

\end{figure}

\text{ }

\newpage

\text{ }

\begin{figure}[b]

\bigskip
\bigskip

\[\begin{aligned}
\text{Formation \quad}&\inferrule{\makebox[\widthof{$\judge{\textcolor{green}{\gamma : \Gamma;\;}x:A\textcolor{green}{(\gamma)};\;y:B(\textcolor{green}{\gamma,\;}x)}{\beta_2(\textcolor{green}{\gamma,\;}x,y):\pi_2\langle x,y\rangle=\beta_1(\textcolor{green}{\gamma,\;}x,y)^*y}$}][c]{$\judge{\textcolor{green}{\gamma : \Gamma}}{A\textcolor{green}{(\gamma)} : \type} \quad\quad \judge{\textcolor{green}{\gamma : \Gamma;\;}x:A\textcolor{green}{(\gamma)}}{B(\textcolor{green}{\gamma,\;}x):\type}$}}{\judge{\textcolor{green}{\gamma : \Gamma}}{\Sigma_{x:A\textcolor{green}{(\gamma)}}B(\textcolor{green}{\gamma,\;}x) : \type}}\\
\\
\\
\\
\\
\text{Introduction \quad}&\inferrule{\makebox[\widthof{$\judge{\textcolor{green}{\gamma : \Gamma;\;}x:A\textcolor{green}{(\gamma)};\;y:B(\textcolor{green}{\gamma,\;}x)}{\beta_2(\textcolor{green}{\gamma,\;}x,y):\pi_2\langle x,y\rangle=\beta_1(\textcolor{green}{\gamma,\;}x,y)^*y}$}][c]{$\judge{\textcolor{green}{\gamma : \Gamma}}{A\textcolor{green}{(\gamma)} : \type} \quad\quad \judge{\textcolor{green}{\gamma : \Gamma;\;}x:A\textcolor{green}{(\gamma)}}{B(\textcolor{green}{\gamma,\;}x):\type}$}}{\judge{\textcolor{green}{\gamma : \Gamma;\;}x:A\textcolor{green}{(\gamma)},y:B(\textcolor{green}{\gamma,\;}x)}{\langle x,y\rangle:\Sigma_{x:A\textcolor{green}{(\gamma)}}B(\textcolor{green}{\gamma,\;}x)}}\\
\\
\\
\\
\\
\text{Projection \quad}&\inferrule{\makebox[\widthof{$\judge{\textcolor{green}{\gamma : \Gamma;\;}x:A\textcolor{green}{(\gamma)};\;y:B(\textcolor{green}{\gamma,\;}x)}{\beta_2(\textcolor{green}{\gamma,\;}x,y):\pi_2\langle x,y\rangle=\beta_1(\textcolor{green}{\gamma,\;}x,y)^*y}$}][c]{$\judge{\textcolor{green}{\gamma : \Gamma}}{A\textcolor{green}{(\gamma)} : \type} \quad\quad \judge{\textcolor{green}{\gamma : \Gamma;\;}x:A\textcolor{green}{(\gamma)}}{B(\textcolor{green}{\gamma,\;}x):\type}$}}{\judge{\textcolor{green}{\gamma : \Gamma;\;}u:\Sigma_{x:A\textcolor{green}{(\gamma)}}B(\textcolor{green}{\gamma,\;}x)}{\pi_1u:A\textcolor{green}{(\gamma)}} \\\\ \\\\ \judge{\textcolor{green}{\gamma : \Gamma;\;}u:\Sigma_{x:A\textcolor{green}{(\gamma)}}B(\textcolor{green}{\gamma,\;}x)}{\pi_2u:B(\textcolor{green}{\gamma,\;}\pi_1(u))}}\\
\\
\\
\\
\\
\text{Prop $\beta$-reduction \quad}&\inferrule{\makebox[\widthof{$\judge{\textcolor{green}{\gamma : \Gamma;\;}x:A\textcolor{green}{(\gamma)};\;y:B(\textcolor{green}{\gamma,\;}x)}{\beta_2(\textcolor{green}{\gamma,\;}x,y):\pi_2\langle x,y\rangle=\beta_1(\textcolor{green}{\gamma,\;}x,y)^*y}$}][c]{$\judge{\textcolor{green}{\gamma : \Gamma}}{A\textcolor{green}{(\gamma)} : \type} \quad\quad \judge{\textcolor{green}{\gamma : \Gamma;\;}x:A\textcolor{green}{(\gamma)}}{B(\textcolor{green}{\gamma,\;}x):\type}$}}{\judge{\textcolor{green}{\gamma : \Gamma;\;}x:A\textcolor{green}{(\gamma)};\;y:B(\textcolor{green}{\gamma,\;}x)}{\beta_1(\textcolor{green}{\gamma,\;}x,y):\pi_1\langle x,y\rangle = x}  \\\\ \\\\ \judge{\textcolor{green}{\gamma : \Gamma;\;}x:A\textcolor{green}{(\gamma)};\;y:B(\textcolor{green}{\gamma,\;}x)}{\beta_2(\textcolor{green}{\gamma,\;}x,y):\pi_2\langle x,y\rangle=\beta_1(\textcolor{green}{\gamma,\;}x,y)^*y}}\\
\\
\\
\\
\\
\text{Prop $\eta$-expansion \quad}&\inferrule{\makebox[\widthof{$\judge{\textcolor{green}{\gamma : \Gamma;\;}x:A\textcolor{green}{(\gamma)};\;y:B(\textcolor{green}{\gamma,\;}x)}{\beta_2(\textcolor{green}{\gamma,\;}x,y):\pi_2\langle x,y\rangle=\beta_1(\textcolor{green}{\gamma,\;}x,y)^*y}$}][c]{$\judge{\textcolor{green}{\gamma : \Gamma}}{A\textcolor{green}{(\gamma)} : \type} \quad\quad \judge{\textcolor{green}{\gamma : \Gamma;\;}x:A\textcolor{green}{(\gamma)}}{B(\textcolor{green}{\gamma,\;}x):\type}$}}{\judge{\textcolor{green}{\gamma : \Gamma;\;}u:\Sigma_{x:A\textcolor{green}{(\gamma)}}B(\textcolor{green}{\gamma,\;}x)}{\eta(\textcolor{green}{\gamma,\;}u):u= \langle \pi_1u,\pi_2u\rangle}}\\
\end{aligned}\]

\bigskip
\bigskip
\bigskip
\bigskip

\caption{\;\;\;Propositional dependent sum types---negative form}
\label{figureA.4.2}\label{propositional sigma char}

\end{figure}

\text{ }

\newpage

\section{Conclusion}\label{section2.7}

In this paper we compared a propositional dependent type theory to an extensional one, and observed that, despite the non-negligible weakening of the former with respect to the latter, there is actually an interesting family of judgements where the two theories have the same deductive power.

This result was obtained by adapting the argument introduced by Hofmann \cite{hofmannconservativity}. In spite of the amount of work in the purely syntactic part of this research, namely in the analysis of the family of canonical equivalences, we underline the central and fundamental role of the soundness property of the semantics induced by the class of categories with attributes. The effect is that, for a given term judgement $t:|T|$ built by $\ett$, where $T$ is an h-elementary type of $\ptt$, we do not know how a corresponding term $\tilde{t}:T$ of $\ptt$ (such that $|\tilde{t}|\equiv t$) is defined. A possible future research direction therefore concerns asking whether there is the possibility of making our argument more constructive from this point of view.

Our argument also applies in the case where $\ptt$ and $\ett$ are extended with a (weak) Tarski universe, requiring that it is provably an h-set in $\ptt$ and that the types associated with its terms are also provably h-sets in $\ptt$. Clearly, in this setup, the universe cannot be univalent in $\ptt$. To address the issue of accommodating a univalent universe, it would be interesting to explore whether the solution proposed by Winterhalter, Sozeau, and Tabareau \cite{winterhalter:hal-01849166}, as applied to two-level type theories, can be used to obtain a corresponding adaptation of our argument.

Our proof, strongly based on canonical equivalences, requires that the family of these be restricted, in such a way that further properties are satisfied. The restriction from general contexts to contexts with h-propositional identities and the one from contexts with h-propositional identities to h-elementary contexts are explicitly applied only at very specific points in the proof. However, they are fundamental for defining---by quotienting---an actual category with attributes, i.e. a model of a dependent type theory. From this perspective, it would then be worthwhile to consider what can be learned by replicating our argument without these restrictions. In detail, one might therefore ask whether the argument can be extended, without restrictions on the family of canonical equivalences, in order to deduce a conservativity result for those \textit{generalised type theories}---with less structural rules---that are modelled by the categorical structure obtained---by quotienting---from the general family of canonical equivalences. More generally, we are interested in looking for other ways of applying Hofmann's argument in order to get similar results or generalised versions of the ones of this paper, depending e.g. on the strictness of the structural rules that we allow in a theory of dependent types.

\section*{Acknowledgment}
  \noindent The author thanks his doctoral supervisors Nicola Gambino and Federico Olimpieri for their helpful suggestions during the development of this research. The author is also grateful to (in alphabetic order) Benedikt Ahrens, Ivan di Liberti, Jacopo Emmenegger, Maria Emilia Maietti, Taichi Uemura for useful discussions on the subject. The author would also like to thank the anonymous referees for their useful comments and suggestions. The research presented in this paper was conducted while the author was affiliated with the School of Mathematics of the University of Leeds, United Kingdom.

\bibliographystyle{alphaurl}
\bibliography{references}

\end{document}